\documentclass[a4paper,12pt,twoside]{article}
\usepackage{graphicx}
\usepackage{a4wide}
\usepackage{sidecap}
\usepackage[font=small,labelfont=bf]{caption}
\usepackage[export]{adjustbox}

\usepackage{setspace}

\usepackage[caption = false]{subfig}
\usepackage{amssymb,amsmath,amsthm}
\usepackage[only,llbracket,rrbracket,interleave]{stmaryrd}
\usepackage{gensymb}

\usepackage[sort,compress]{cite}

\usepackage{booktabs,array,dcolumn}
\newcommand{\PreserveBackslash}[1]{\let\temp=\\#1\let\\=\temp}
\newcolumntype{C}[1]{>{\PreserveBackslash\centering}p{#1}}
\newcolumntype{R}[1]{>{\PreserveBackslash\raggedleft}p{#1}}
\newcolumntype{L}[1]{>{\PreserveBackslash\raggedright}p{#1}}

\usepackage{xcolor}

\newcommand{\bm}[1]{\text{\boldmath $#1$\unboldmath}}
\newcommand{\abs}[1]{\lvert#1\rvert}
\newcommand{\norm}[1]{\lVert#1\rVert}
\DeclareMathOperator{\sgn}{sgn}
\newcommand{\sign}[1]{{\sgn \left( #1 \right)}}
\newcommand{\sprod}[2][\Omega]{\left(#2\right)_{#1}}
\newcommand{\dprod}[2][\Gamma]{\left\langle#2\right\rangle_{#1}}
\newcommand{\bcdot}{\operatorname{\bm{\cdot}}}
\newcommand{\vect}[1]{\bm{#1}}
\newcommand{\mat}[1]{\bm{#1}}
\newcommand{\bvect}[1]{\mathbf{#1}}
\newcommand{\bmat}[1]{\mathbf{#1}}

\newcommand{\Div}{{\bm{\nabla}\bcdot\,}}
\newcommand{\Grad}{\bm{\nabla}}
\newcommand{\Sym}{\bm{\nabla^{\texttt{S}}}}

\newcommand{\pd}[2]{\frac{\partial{#1}}{\partial{#2}}}

\newcommand{\RR}{\mathbb{R}}
\newcommand{\testh}{\ensuremath{\mathcal{W}^h}}

\newcommand{\testmh}{\ensuremath{\widehat{\mathcal{W}}^h}}

\newcommand{\eltwo}{\ensuremath{\mathcal{L}_2}}
\newcommand{\Sobo}{\ensuremath{\mathcal{H}}}
\newcommand{\Czero}{\ensuremath{\mathcal{C}^0}}
\newcommand{\Polyk}{\ensuremath{\mathcal{P}^{k}}}

\newcommand{\nsd}{\ensuremath{\texttt{n}_{\texttt{sd}}}}
\newcommand{\msd}{\ensuremath{\texttt{m}_{\texttt{sd}}}}
\newcommand{\numel}{\ensuremath{\texttt{n}_{\texttt{el}}}}
\newcommand{\nen}{\ensuremath{\texttt{n}_{\texttt{en}}}}

\newcommand{\nlayers}{\ensuremath{\texttt{n}_{\texttt{lay}}}}
\newcommand{\ndiv}{\ensuremath{\texttt{n}_{\texttt{div}}}}
\newcommand{\tras}{^{{T}}}
\newcommand{\tr}{\text{tr}}
\newcommand{\bu}{\bm{U}}
\newcommand{\bv}{\bm{v}}
\newcommand{\bw}{\bm{W}}
\newcommand{\bn}{\bm{n}}

\newcommand{\bbAll}{\bm{B}}
\newcommand{\bb}{\widehat{\bbAll}}
\newcommand{\bhu}{\widehat{\bu}}
\newcommand{\bhv}{\widehat{\bv}}

\newcommand{\bF}{\bm{F}}
\newcommand{\bG}{\bm{G}}
\newcommand{\bV}{\bm{V}}
\newcommand{\bW}{\bm{W}}

\newcommand{\btau}{\boldsymbol{\tau}}
\newcommand{\bq}{\bm{q}}
\newcommand{\bsigma}{\boldsymbol{\sigma}^d}
\newcommand{\beps}{\boldsymbol{\varepsilon}^{\! d}}
\newcommand{\bphi}{\boldsymbol{\phi}}
\newcommand{\buE}{\bu_{\! e}}
\newcommand{\bepsE}{\beps_{\! e}}
\newcommand{\bphiE}{\bphi_{\! e}}
\newcommand{\bUinf}{\bm{U}_{\! \infty}}
\newcommand{\Tend} {\textrm{T}_{\texttt{end}}}
\newcommand{\lamax} {\lambda_{\texttt{max}}}
\newcommand{\hlamax} {\widehat{\lambda}_{\texttt{max}}}
\newcommand{\Id}[1]{\bmat{I}_{#1}}

\newcommand{\hu}{\widehat{\bvect{U}}}
\newcommand{\re}{\bvect{R}_e}
\newcommand{\hre}{\widehat{\bvect{R}}_e}

\newcommand{\bznode}{\bm{Z}_{\! e}}
\newcommand{\bznodeAll}{\bm{Z}}

\newcommand{\bzE}{\bvect{Z}_{e}}
\newcommand{\AmatE}[2]{\bmat{A}_{{#1}{#2}}^e}
\newcommand{\FvectE}[1]{\bvect{F}_{{#1}}^e}
\newcommand{\Kmat}{\bmat{K}}
\newcommand{\Fvect}{\bvect{F}}
\newcommand{\jump}[1]{\llbracket #1\rrbracket}

\newcommand{\Rey}{\small{\textsl{Re}}}
\newcommand{\Pra}{\small{\textsl{Pr}}}
\newcommand{\Ma}{\small{\textsl{M}}}
\newcommand{\vinf}{v_{\! \infty}}
\newcommand{\rhoinf}{\rho_{\! \infty}}
\newcommand{\muinf}{\mu_{\! \infty}}
\newcommand{\Minf}{\Ma_{\! \infty}}

\newcommand{\cp}{c_p}
\newcommand{\cv}{c_v}
\newcommand{\Tref}{\ensuremath{T_{\texttt{ref}}}}

\newtheorem{remark}{Remark}

\usepackage{scalerel}
\usepackage{stackengine}
\stackMath
\def\hatgap{0pt}
\def\subdown{-2pt}
\newcommand\reallywidehat[2][]{
	\renewcommand\stackalignment{l}
	\stackon[\hatgap]{#2}{
		\stretchto{\scalerel*[\widthof{$#2$}]{\kern-.6pt\bigwedge\kern-.6pt}{\rule[-\textheight/2]{1ex}{\textheight}}}{0.5ex}_{\smash{\belowbaseline[\subdown]{\scriptstyle#1}}}}}

\graphicspath{{figsPNG/}{figsPDF/}}

\begin{document}
\title{Hybridisable discontinuous Galerkin formulation of compressible flows}

\author{Jordi Vila-P\'erez, Matteo Giacomini, Ruben Sevilla and Antonio Huerta}

\author{
	\renewcommand{\thefootnote}{\arabic{footnote}}
	Jordi Vila-P\'erez\footnotemark[1]\textsuperscript{\, ,}\footnotemark[2] \ ,
	Matteo Giacomini\footnotemark[1]\textsuperscript{\, ,}\footnotemark[3], \
	Ruben Sevilla\footnotemark[4], \
	Antonio Huerta\footnotemark[1]\textsuperscript{\, ,}\footnotemark[3]
}

\date{\today}
\maketitle

\renewcommand{\thefootnote}{\arabic{footnote}}

\footnotetext[1]{Laboratori de C\`alcul Num\`eric (LaC\`aN), ETS de Ingenieros de Caminos, Canales y Puertos, Universitat Polit\`ecnica de Catalunya, Barcelona, Spain.}
\footnotetext[2]{Barcelona Graduate School of Mathematics (BGSMath), Barcelona, Spain.}
\footnotetext[3]{International Centre for Numerical Methods in Engineering (CIMNE), Barcelona, Spain.}
\footnotetext[4]{Zienkiewicz Centre for Computational Engineering, College of Engineering, Swansea University, Bay Campus, SA1 8EN, Wales, United Kingdom.
	\vspace{5pt}\\
	Corresponding author: Jordi Vila-P\'erez \textit{E-mail:} \texttt{jordi.vila.perez@upc.edu}
}

\begin{abstract}
	This work presents a review of high-order hybridisable discontinuous Galerkin (HDG) methods in the context of compressible flows.
	Moreover, an original unified framework for the derivation of Riemann solvers in hybridised formulations is proposed.
	This framework includes, for the first time in an HDG context, the HLL and HLLEM Riemann solvers as well as the traditional Lax-Friedrichs and Roe solvers.
	HLL-type Riemann solvers demonstrate their superiority with respect to Roe in supersonic cases due to their positivity preserving properties.
	In addition, HLLEM specifically outstands in the approximation of boundary layers because of its shear preservation, which confers it an increased accuracy with respect to HLL and Lax-Friedrichs.
	A comprehensive set of relevant numerical benchmarks of viscous and inviscid compressible flows is presented.
	The test cases are used to evaluate the competitiveness of the resulting high-order HDG scheme with the aforementioned Riemann solvers and equipped with a shock treatment technique based on artificial viscosity.
\end{abstract}


\section{Introduction}

High-order methods have experienced a growing interest within the computational fluid dynamics (CFD) community because of their increased accuracy when compared to low-order methods \cite{Kroll2009,Wang2013}.
However, low-order finite volume (FV) or stabilised finite element (FE) methods are still the most employed strategies in commercial, industrial and open source CFD solvers \cite{sorensen2003multigrid,jasak2009openfoam,gerhold2005overview,chalot2004industrial}.
Low-order FV and FE methods are robust, easy to implement and provide a competitive alternative for the computation of steady state CFD solutions.
Nevertheless, the higher diffusion and dispersion errors introduced by such discretisations when compared to high-order methods limit their performance in problems involving transient effects \cite{Ekaterinaris2005,Drikakis2003,Moura-MSP:2015}.
This has prompted the extension of FV and stabilised FE schemes to high-order~\cite{Chalot-CN:2010,chassaing2013accuracy,Sevilla-SHM:2013,2017hsupgALE}. 
The development of high-order schemes is, thus, claimed necessary in order to tackle a variety of complex flow phenomena arising in many practical aerodynamic problems, such as the resolution of shear layers or the propagation of vortices over long distances and for long times~\cite{NASAvision,Wang2013}.

In the context of high-order discretisations, discontinuous Galerkin (DG) methods have become one of the most adopted approaches within the computational engineering community~\cite{Cockburn-CKS:2000,Arnold-ABCM:2002,Bassi-BR:1997b}.
In particular, DG discretisations have been often seen as a methodology to combine the advantages of both FV and FE schemes.
Contrary to FV methods, DG methods allow to define high-order local approximations~\cite{Cockburn-CLS:1989,Cockburn-CS:1998,Bassi-BR:1997b,Bassi-BR:2002}.
In addition, in DG methods, the stabilisation term required for solving convection dominated problems is easier to define when compared to traditional stabilised FE methods~\cite{Sevilla-SHM:2013,Chalot-CN:2010}.
The DG framework allows to devise high-order numerical methods that enforce element-by-element conservation and provides a suitable discretisation on unstructured meshes~\cite{Cockburn-CKS:2000,Ern2011}.
In addition, it permits an efficient exploitation of parallel computing architectures~\cite{Roca-RNP:2013,Fernandez-FNP:2017} and an easy implementation of adaptive strategies for non-uniform degree approximations \cite{HartmannHouston:2003,Cangiani2017,Sanjay2020,Giorgiani-GFH:2014,Balan-BWM:2015,}.
However, the duplication of nodes at the interface of neighbouring elements has limited its application mostly to academic problems, see the discussion in~\cite{MG-GS:19} and references therein.

Accordingly, hybrid discretisation methods, e.g. the hybrid/hybridised DG method~\cite{Egger-ES-09,Egger-EW:2012a,Egger-EW:2012b}, the hybridisable discontinuous Galerkin (HDG) methods~\cite{Jay-CGL:09,Cockburn-CS1998a,Cockburn-CG:2004,Jay-CG:2009} and the hybrid high-order (HHO) method~\cite{Ern-DPE-15,DiPietro-DPEL:2014,CockburnErn-HHO:16}, obtained from the hybridisation of traditional DG schemes, have been devised as a significantly less expensive alternative~\cite{Huerta-HARP:2013,May-WBMS-14}.
The HDG approach reduces the number of globally coupled degrees of freedom via the introduction of a hybrid variable, namely the trace of the unknown on the mesh faces, and appropriately defined inter-element numerical fluxes.
Recently, special attention has been devoted to the HDG method which relies on a mixed formulation for second-order problems~\cite{Cockburn-CG:2004,Cockburn-CG:2005a,Cockburn-CG:2005b,Jay-CGL:09,Nguyen-NPC:2009lCD,Nguyen-NPC:2009nlCD,Nguyen-NPC:20011INS,Nguyen-NP:2012,Giorgiani-GFH:2014,Shi-QS-16,RS-SH:16,Tutorial-GSH:2020}.

In the context of compressible flows, different hybrid methods, such as HDG~\cite{Peraire-PNC:10,Nguyen-NP:2012,Williams-18,Sanjay2020,Jaust-JS:2014,Jaust-JSW:2015}, the embedded DG (EDG)~\cite{Peraire-PNC:11}, the interior embedded DG (IEDG)~\cite{Nguyen-NPC-EDG2015} or hybrid mixed methods~\cite{May-SM-13,Schutz-SM:2013,Schutz-SWM:2012}, have been devised for the formulation of the inviscid Euler and the laminar compressible Navier-Stokes equations.
The HDG formulation has also been extended to turbulent compressible flows, both solving the Reynolds-averaged Navier-Stokes equations combined with the Spalart-Allmaras~\cite{Moro-MNP:2011,Moro2017} or the $\kappa-\omega$~\cite{Woopen-WLM:2014} turbulence models, or by means of a large-eddy simulation approach~\cite{Fernandez-FNP:2017}.
It is worth noting that in the inviscid limit, i.e. for the Euler equations, HDG methods based on primal and mixed formulations are equivalent.

A salient feature of the HDG method stemming from~\cite{Peraire-PNC:10,Nguyen-NP:2012} is its associated optimal order of convergence for the viscous stress and the heat flux. 
It follows that the HDG method provides an increased accuracy in the computation of typical quantites of interest in aerodynamic applications, such as lift and drag.
The optimal accuracy properties of the method in the approximation of the stress and heat flux rely on the equal-order approximation of the primal, mixed and hybrid variables.
In addition, the resulting HDG discretisation is robust in the incompressible limit, circumventing the Ladyzhenskaya-Babu\v ska-Brezzi (LBB)~\cite{Jay-CG:2009}.
In this context, when the Cauchy stress tensor formulation is employed for the momentum equation, the appropriate choice of the discretisation space for the mixed variable is crucial to ensure the optimal convergence of the method. This has been achieved by means of the $M$-decomposition framework~\cite{CockburnFu-Mdecomposition1:2016,CockburnFu-Mdecomposition2:2016,CockburnFu-Mdecomposition3:2016,Cockburn-Mdecomposition:2016,Cockburn-Mdecomposition:2017}, the utilisation of the reduced stabilisation~\cite{Lehrenfeld2016,Shi-QS-16,Oikawa2015} or the employment of a pointwise symmetric formulation of the stress tensor~\cite{Sevilla-SGKH:2018,Giacomini-GKSH:2018,Tutorial-GSH:2020}.


When convection phenomena are considered, e.g. in the context of systems of conservation laws and nonlinear hyperbolic partial differential equations (PDEs), the definition of the numerical fluxes has a seminal importance in the accuracy and stability of the approximate solution. For this reason, it has been object of intensive study by means of Riemann solvers, both in the context of traditional DG \cite{Cockburn-CS:1998,Mengaldo2014,Qiu-QKS:2006} and in low-order FV methods, see for instance the monographs by Toro~\cite{Toro2009}, Leveque~\cite{Leveque2013} and Hesthaven~\cite{Hesthaven2017}.
In contrast, the definition of approximate Riemann solvers for HDG methods have received considerably less attention, and only the traditional Lax-Friedrichs and Roe solvers have been considered~\cite{Nguyen-NP:2012,Peraire-PNC:10,Peraire-PNC:11}. \\

This work presents a review of the high-order HDG formulation of compressible flows, including both the inviscid Euler and the viscous compressible Navier-Stokes equations.
Moreover, the study proposes a unified framework for the derivation of Riemann solvers in hybridised formulations.
The framework includes the existing Lax-Friedrichs and Roe solvers and formulates, for the first time in the context of HDG, the HLL~\cite{Harten-HLL:1983} and the HLLEM~\cite{Einfeldt1988,Einfeldt1991} Riemann solvers.
The use of Riemann solvers of the HLL family is especially important in the context of supersonic flows, where the Roe numerical flux may fail to provide physically admissible solutions because of a lack of dissipation~\cite{Peery1988,Quirk1994}, whereas Lax-Friedrichs produces over-dissipative approximations~\cite{Moura-MMPS:2017,Moura-MSP:2015report}. On the contrary, HLL-type Riemann solvers provide a robust framework to compute accurate solutions while guaranteeing positiveness of the approximate density and pressure fields~\cite{Fleischmann-FAHA:2019,Quirk1994}. 
Furthermore, the HLLEM Riemann solver is also robust in the preservation of shear and contact waves \cite{Einfeldt1988,Einfeldt1991,Dumbser2016}, likewise Roe, thus improving the Lax-Friedrichs and HLL approximation of such kind of waves.

Additionally, this work introduces a mixed formulation of the compressible Navier-Stokes equations with strongly enforced symmetry of the viscous stress tensor.
Such approach uses the same discrete spaces for the primal and mixed variables and retrieves optimal convergence properties of the stress tensor and the heat flux, with reduced computational cost.

Finally, this study presents an exhaustive set of numerical benchmarks, spanning from subsonic flows to supersonic inviscid and viscous cases with shocks, that allow to verify the capabilities of the HDG method while examining the properties of the presented Riemann solvers.

In this work, the main focus is in the HDG formulation for compressible flows and the presentation of a unified framework for the Riemann solvers in HDG.
To this end, the examples considered involve steady state flows.
When the steady state is computed using a pseudo-time approach, the backward Euler scheme is employed.
For transient flows, the HDG method has been combined with a variety of low and high-order time integrators~\cite{Nguyen-NP:2012,Jaust-JSW:2015,Jaust-JS:2014,Jaust-JSW:2014,Fernandez-FNP:2017,Sanjay2020}.
Although less explored, there are also works where the HDG method has been employed with explicit time-marching algorithms~\cite{Samii2018,Samii2019}.

The remainder of this paper is organised as follows.
The compressible Navier-Stokes equations, governing compressible flows, are described in section~\ref{sc:compressibleNS}.
In section~\ref{sc:HDGformulation}, the HDG formulation of the compressible Navier-Stokes equations is detailed.
Section~\ref{sc:RiemannSolvers} presents a unified description of the Riemann solvers in the context of high-order HDG methods. Specifically, the HLL and HLLEM Riemann solvers are proposed for hybrid discretisations.
In section~\ref{sc:Implementation}, the solution strategy of the HDG solver for the resulting nonlinear problem and the numerical treatment of solutions with discontinuities and sharp gradients is discussed.
Section~\ref{sc:Convergence} examines the optimal accuracy properties of the computational method in a pair of convergence studies for inviscid and viscous flows.
A set of numerical benchmarks for a variety of flow conditions is then presented in section~\ref{sc:Benchmarks} to test the performance and robustness of the high-order HDG solver.
Finally, section~\ref{sc:Conclusions} summarises the main results of this study.

\section{Compressible flow equations}
\label{sc:compressibleNS}

Let $\Omega \subset \RR^{\nsd}$ be an open bounded domain with boundary $\partial \Omega$, being $\nsd$ the number of spatial dimensions, and $\Tend>0$ the final time of interest. 
The Navier-Stokes equations, governing unsteady viscous compressible flows in absence of external body forces are expressed in nondimensional conservation form as
\begin{equation}\label{eq:NScompact}
\pd{\bu}{t} + \Div \left( \bF (\bu)- \bG (\bu, \Grad \bu) \right) = \vect{0}, \quad \text{in } \Omega \times (0,\Tend],
\end{equation}
where $\bu \in \RR^{\nsd + 2}$ is the vector of dimensionless conservative variables and $\bF$ and $\bG \in \RR^{(\nsd + 2)\times \nsd}$ are the advection and diffusion flux tensors, respectively, given by
\begin{equation}\label{eq:NSterms}
\bu = \begin{Bmatrix} \rho\\ \rho\bv\\ \rho E \end{Bmatrix}\!,
\quad
\bF (\bu) = \begin{bmatrix} \rho \bv\tras\\
\rho\bv \otimes\bv + p \Id{\nsd}\\
(\rho E + p) \bv\tras \end{bmatrix} \!,
\quad
\bG (\bu, \Grad \bu)= \begin{bmatrix} \vect{0}\\
  \bsigma \\
(\bsigma \bv + \bq )\tras \end{bmatrix}.
\end{equation}

In these expressions, $\rho$ denotes the density, $\bv$ is the velocity vector, $E$ is the total specific energy, $p$ is the pressure, $\bsigma$ is the viscous stress tensor and $\bq$ is the heat flux.

The flow is assumed to obey the ideal gas law $\gamma p = (\gamma - 1) \rho T$, where $T$ is the temperature, and $\gamma = \cp / \cv$ is the ratio of specific heats at constant pressure, $\cp$, and constant volume, $\cv$, and takes value $\gamma = 1.4$ for air.
Moreover, for a calorically perfect gas, it holds that $p = (\gamma - 1) \rho \left( E -  \norm{\bv}^2/2 \right)$.

Under Stokes' hypothesis, the viscous stress tensor is expressed as
\begin{equation}\label{eq:StressTensor}
\bsigma = \frac{\mu}{\Rey} \left(2 \Sym \bv - \frac{2}{3} (\Div \bv) \Id{\nsd} \right),
\end{equation}
where $\Sym := (\Grad + \Grad\tras)/2$ is the symmetric part of the gradient operator.

\begin{remark}[Cauchy stress tensor]
	The Cauchy stress tensor, $\boldsymbol{\sigma}$, which assembles the mechanical stresses of the fluid, is the combination of the viscous stress tensor $\bsigma$ and the thermodynamical pressure $p$, that is $\boldsymbol{\sigma} = \bsigma - p \Id{\nsd}$.
\end{remark}

In addition, the heat flux is modelled according to Fourier's law of heat conduction, that is
\begin{equation}\label{eq:HeatFlux}
\bq = \frac{\mu}{\Pra \Rey} \Grad T,
\end{equation}
and the nondimensional dynamic viscosity, $\mu$, depends on the temperature following Sutherland's law, i.e.
\begin{equation}\label{eq:Sutherland}
\mu = \left(\frac{T}{T_\infty}\right)^{3/2} \frac{T_\infty + S}{T + S},
\end{equation}
where the non-dimensional free-stream temperature and the Sutherland constant are expressed, respectively, as $T_\infty = 1/\left((\gamma - 1) \Minf^2\right)$ and $S = S_0/\left((\gamma - 1) \Tref \Minf^2\right)$, with $S_0 = 110K$ for a reference temperature of $\Tref = 273K$.

The nondimensional description of the problem is completed with the definition of the Reynolds, Prandtl and Mach numbers, i.e., $\Rey = \rhoinf \vinf L / \muinf$, $\Pra = \cp \muinf / \kappa$ and $\Minf = \vinf / c_\infty$, respectively, where $c = \sqrt{\gamma p/\rho}$ denotes the speed of sound, $L$ is a characteristic length and $\kappa$ stands for the thermal conductivity.
Such quantities are expressed in terms of reference free-stream values, indicated by the subscript $\infty$.
The Prandtl number is considered constant and equal to $\Pra = 0.71$ for air.

The problem is closed with the prescription of initial and boundary conditions, namely
\begin{equation}\label{eq:BC-In_NS}
\begin{aligned}
\bu &= \bu^0  &  \text{in } & \Omega \times \{0\}, \\
\bbAll (\bu, \Grad \bu) &= \bm{0}  & \text{on } &\partial \Omega \times (0,\Tend], 
\end{aligned}
\end{equation}
where $\bu^0$ stands for an initial state and the vector $\bbAll$ describes a boundary condition operator, imposing inflow, outflow or wall conditions with isothermal, adiabatic or symmetry properties as detailed in section~\ref{ssc:BoundaryConditions}.

\begin{remark}[Compressible Euler equations]
	The compressible Euler equations are recovered in the inviscid limit, that is when $\Rey \to \infty$. In such case, the set of conservation equations~\eqref{eq:NScompact} becomes a system of first-order hyperbolic PDEs, namely
	\begin{equation}\label{eq:EulerCompact}
	\pd{\bu}{t} + \Div \bF (\bu)= \vect{0}, \quad \text{in } \Omega \times (0,\Tend].
	\end{equation}
\end{remark}

\section{HDG formulation of the compressible Navier-Stokes equations}
\label{sc:HDGformulation}

Consider a partition of the domain $\Omega$ in $\numel$ disjoint subdomains $\Omega_e$ such that $\Omega = \bigcup_{e=1}^{\numel} \Omega_e$.
Let $\Gamma$ denote the \emph{mesh skeleton} or internal interface, namely
\begin{equation}
\Gamma:=\left[ \bigcup_{e=1}^{\numel}\partial\Omega_e \right] \setminus \partial \Omega.
\end{equation}

In addition, the notation for the \emph{jump} operator, $\jump{\circledcirc} = \circledcirc^+ + \circledcirc^-$, is introduced, defining the sum of the values in the elements $\Omega^+$ and $\Omega^-$ at each side of the internal interface $\Gamma$, respectively~\cite{AdM-MFH:08}.

\subsection{Mixed variables for the compressible Navier-Stokes equations}
\label{ssc:MixedFormulation}

One of the main features of the HDG mixed formulation is the introduction of mixed variables for the approximation of derivative terms in second-order problems~\cite{Jay-CGL:09,Jay-CG:2009,Nguyen-NP:2012,Cockburn2016}.
In the case of the compressible Navier-Stokes equations, the mixed variables are responsible for the description of the viscous stress tensor $\bsigma$ and the heat flux $\bq$ appearing in the viscous fluxes~\eqref{eq:NSterms}.

Usual mixed formulations of the compressible Navier-Stokes equations introduce the gradient of the primal variable, $\Grad \bu$, as mixed variable~\cite{Peraire-PNC:10,Nguyen-NP:2012,May-WBMS-14, Fernandez-FNP:2017}.
The advantage of using $\Grad \bu$ is its linear expression with respect to the primal variable $\bu$. Then, \eqref{eq:NScompact} is rewritten as a system of first-order PDEs with an additional linear equation, that is
\begin{equation} \label{eq:NSmixedGradient} 
\left\{\begin{aligned}
\bm{Q} - \Grad \bu &= \bm{0},  \\
\pd{\bu}{t} +\Div \left( \bF(\bu) - \bG(\bu,\bm{Q}) \right) &= \bm{0},  \\
\end{aligned}\right.
\end{equation}
where the viscous stress tensor and the heat flux appearing in $\bG(\bu,\bm{Q})$~\eqref{eq:NSterms} are given by
\begin{subequations} \label{eq:ViscousTermsGradient}
	\begin{align}
	& \bsigma = \frac{1}{\Rey} \frac{\mu}{\rho} \left[ \Sym (\rho \bv) - \frac{1}{\rho} \left( \rho \bv \otimes \Grad \rho +  \Grad \rho \otimes \rho \bv \right) - \frac{2}{3} \left( \Div (\rho \bv) - \frac{1}{\rho} \Grad \rho \cdot \rho \bv \right) \Id{\nsd} \right], \\
	& \bq = \frac{\gamma}{\Rey \Pra} \frac{\mu}{\rho} \left[ \Grad (\rho E) - \frac{\rho E}{\rho} \Grad \rho - \frac{1}{\rho} \left( \Grad (\rho \bv)\tras -  \frac{1}{\rho} \Grad \rho \otimes \rho \bv \right) \rho\bv \right].
	\end{align}
\end{subequations}
It is worth noticing that the viscous stresses and the heat flux are linear with respect to the mixed variable. However, their expression presents a number of non-linearities with respect to the conservation variables.

An alternative formulation, inspired by the mechanical description of the problem, employs the deviatoric strain rate tensor
\begin{equation}\label{eq:DeviatoricStrainRate}
\beps = 2 \Sym \bv - \frac{2}{3} (\Div \bv) \Id{\nsd},
\end{equation}
and the gradient of temperature $\bphi =  \Grad T$ as mixed variables for the HDG formulation. The resulting system of first-order PDEs is given by
\begin{equation} \label{eq:NSmixedPhysical} 
\left\{\begin{aligned}
\beps - 2 \Sym \bv - \frac{2}{3} (\Div \bv) \Id{\nsd} &= \bm{0},  \\
\bphi - \Grad T &= \bm{0},  \\
\pd{\bu}{t} +\Div \left( \bF(\bu) - \bG(\bu,\beps,\bphi) \right) &= \bm{0},  \\
\end{aligned}\right.
\end{equation}
where the viscous stress tensor and the heat flux in $\bG(\bu,\beps,\bphi)$~\eqref{eq:NSterms} can be expressed in a neat manner as
\begin{equation} \label{eq:ViscousTermsNewMixed} 
\bsigma = \frac{\mu}{\Rey} \beps, \qquad
\bq = \frac{\mu}{\Rey \Pra} \bphi.
\end{equation}

Note that, whereas such mixed variables are nonlinear with respect to the conservation variables, this choice vastly reduces the number of nonlinearities and simplifies the expression of the viscous fluxes, in contrast to~\eqref{eq:ViscousTermsGradient}.

\begin{remark}
	Such choice for the mixed variables, which resembles the mixed formulation proposed in~\cite{Williams-18}, involves a reduced number of degrees of freedom, when compared to $\bm{Q} = \Grad \bu$, thus decreasing the computational cost of the local problems.
\end{remark}

\subsection{Strong form of the local and global problems}
\label{ssc:StrongForm}

In this work, the deviatoric strain rate tensor $\beps$ and the temperature gradient $\bphi$ are considered as mixed variables. The problem is then written as a system of first-order PDEs, in mixed form, in the so-called \emph{broken computational domain}.

The HDG method solves the problem in two stages.
First, $\numel$ local problems, given by
\begin{equation} \label{eq:NSstrongLocal} 
\left\{\begin{aligned}
\bepsE - \left( 2 \Sym \bv_{\! e} - \frac{2}{3} (\Div \bv_{\! e}) \Id{\nsd}\right)  &= \bm{0} & \text{in } &\Omega_e \times (0,\Tend],  \\
\bphiE - \Grad T_{\! e} &= \bm{0} & \text{in } &\Omega_e \times (0,\Tend],  \\
\pd{\buE}{t} +\Div \left( \bF(\buE) - \bG(\buE,\bepsE,\bphiE) \right) &= \bm{0}  & \text{in } &\Omega_e \times (0,\Tend],  \\
\buE &= \bu^0  &  \text{in } & \Omega_e \times \{0\},  \\
\buE &= \bhu  & \text{on } &\partial \Omega_e \times (0,\Tend], \\
\end{aligned}\right.
\end{equation}
for $e = 1,\dotsc,\numel$, define the solution $(\buE, \bepsE, \bphiE)$ in each element as a function of an independent variable $\bhu$, representing the trace of the solution on $\Gamma \cup \partial \Omega$.

Then, $\bhu$ is computed as the solution of a global problem imposing boundary conditions on $\partial\Omega$ and enforcing inter-element continuity of the solution and of the normal fluxes on $\Gamma$ via the so-called \emph{transmission conditions}, namely 
\begin{equation} \label{eq:NSstrongGlobal} 
\left\{\begin{aligned}
\bb(\bu,\bhu,\beps,\bphi) &= \bm{0},  & \text{on } &\partial \Omega \times (0,\Tend], \\
\jump{\bu\otimes \bn} &= \mat{0} &\text{on } &\Gamma \times (0,\Tend],\\
\jump{\left( \bF(\bu) - \bG(\bu,\beps,\bphi) \right) \bn} &= \bm{0} &\text{on } &\Gamma \times (0,\Tend],
\end{aligned}\right.
\end{equation}
where $\bn$ is the outward unit normal vector and the boundary trace operator $\bb(\bu,\bhu,\beps,\bphi)$ imposes the boundary conditions along $\partial \Omega$ exploiting the hybrid variable.

Note that the second equation in \eqref{eq:NSstrongGlobal} is automatically satisfied due to the Dirichlet boundary condition $\buE = \bhu$ imposed in the local problems \eqref{eq:NSstrongLocal} and by the fact that the hybrid variable $\bhu$ is unique on each face of the mesh skeleton.

\subsubsection{Boundary conditions}
\label{ssc:BoundaryConditions}
The global system~\eqref{eq:NSstrongGlobal} involves the boundary trace operator $\bb(\bu,\bhu,\beps,\bphi)$, whose definition depends on the type of boundary under analysis.
Following the philosophy of \cite{Peraire-PNC:10,Nguyen-NP:2012,Fernandez-FNP:2017,Mengaldo2014}, different definitions of boundary conditions that commonly arise in simulation of compressible flow problems are presented in table~\ref{tb:boundaryConditions}.

To this effect, consider a partition of the boundary $\partial \Omega$ such that $\partial \Omega = \Gamma_\infty \cup \Gamma_{\text{out}} \cup \Gamma_{\text{ad}} \cup \Gamma_{\text{iso}} \cup \Gamma_{\text{inv}}$ and the subdomains $\Gamma_\infty$, $\Gamma_{\text{out}}$, $\Gamma_{\text{ad}}$, $\Gamma_{\text{iso}} $ and $\Gamma_{\text{inv}}$ are disjoint by pairs.
Here, $\Gamma_\infty$ accounts for a far-field boundary type, $\Gamma_{\text{out}}$ denotes a subsonic outflow with imposed pressure, $\Gamma_{\text{ad}}$ and $\Gamma_{\text{iso}}$ refer to adiabatic and isothermal walls, respectively, and $\Gamma_{\text{inv}}$ stands for an inviscid wall with slip conditions or a symmetry wall.

\begin{table} [htbp]
	\centering
	\caption{Definition of boundary conditions for compressible flow problems using a hybrid discretisation.}
	\makebox[\linewidth]{
		\begin{tabular}{ L{0.5cm} L{4.6cm} L{11cm}}
			\toprule
			& Boundary type & Boundary condition operator \\
			\midrule
			$\Gamma_\infty$ & Far-field, subsonic inflow, supersonic inflow/outflow & $\bb = \mat{A}_n^+(\bhu) (\buE - \bhu) + \mat{A}_n^-(\bhu) (\bUinf - \bhu)$, \\
			$\Gamma_{\text{out}}$ & Subsonic outflow (pressure outflow) & $\bb = \begin{Bmatrix} \rho_{\! e} - \widehat{\rho}, \, \left[\rho \bv_{\! e} - \widehat{\rho \bv} \right]\tras,\,  p_{\text{out}}/(\gamma - 1) + \rho_{\! e} \norm{\bv_{\! e}}^2/2 - \widehat{\rho E} \end{Bmatrix}\tras$, \\
			$\Gamma_{\text{ad}}$ & Adiabatic wall & $\bb = \begin{Bmatrix} \rho_{\! e} - \widehat{\rho},\,  \widehat{\rho \bv}\tras,\,  (\mu / \Rey \Pra) \bphiE \bn - \tau^d_{\rho E} \left(\rho E_{e} - \widehat{\rho E}\right) \end{Bmatrix}\tras$, \\
			$\Gamma_{\text{iso}}$ & Isothermal wall & $\bb = \begin{Bmatrix} \rho_{\! e} - \widehat{\rho},\,  \widehat{\rho \bv}\tras,\,  \rho_{\! e} T_w/\gamma - \widehat{\rho E} \end{Bmatrix}\tras$, \\
			$\Gamma_{\text{inv}}$ & Inviscid wall or symmetry surface & $\bb = \begin{Bmatrix} \rho_{\! e} - \widehat{\rho},\,  \left[ (\Id{\nsd} - \bn\otimes \bn) \rho \bv_{\! e} - \widehat{\rho \bv} \right]\tras,\,  \rho E_{e} - \widehat{\rho E} \end{Bmatrix}\tras$. \\
			\bottomrule
		\end{tabular}
	}
	\label{tb:boundaryConditions}
\end{table}

In the expressions in table~\ref{tb:boundaryConditions}, $p_{\text{out}}$ and $T_w$ stand for prescribed values of outflow pressure and wall temperature, respectively, and $\tau^d_{\rho E} = 1/\left[(\gamma - 1) \Minf^2 \Rey \Pra \right]$ is a diffusive stabilisation term for the heat flux~\cite{Peraire-PNC:11,Nguyen-NP:2012}.

Moreover, note that inflow and outflow boundary conditions on $\Gamma_\infty$ are imposed in a characteristics-based approach using the Jacobian matrix of the convective flux in the normal direction to the boundary, namely $\mat{A}_n(\bhu):=[\partial \bF(\bhu) / \partial \bhu] \cdot \bn$.
The spectral decomposition of the matrix, $\mat{A}_n(\bhu) = \mat{R} \mat{\Lambda} \mat{L}$ is then computed, where $\mat{\Lambda}$, $\mat{R}$ and $\mat{L}$ denote the matrices of eigenvalues, right eigenvectors and left eigenvectors, respectively.
Finally, the matrices $\mat{A}_n^-$ and $\mat{A}_n^+$ are defined as $\mat{A}_n^\pm := (\mat{A}_n \pm \abs{\mat{A}_n})/2$, where $\abs{\mat{A}_n(\bhu) } := \mat{R} \abs{\mat{\Lambda}} \mat{L}$ and the matrix $ \abs{\mat{\Lambda}}$ is a diagonal matrix containing the absolute value of the eigenvalues in $\mat{\Lambda}$. The expression of the matrices of eigenvectors and eigenvalues, $\mat{R}$, $\mat{L}$ and $\mat{\Lambda}$, can be found in~\cite{Rohde2001}.

\subsection{Weak form of the local and global problems}

Following the notation in~\cite{RS-SH:16,Tutorial-GSH:2020}, the following discrete functional spaces
\begin{subequations}
	\begin{align}
	\testh(\Omega)    &:=\bigl\{ w \in \eltwo(\Omega): w \vert_{\Omega_e} \in \Polyk(\Omega_e) \;\forall \Omega_e, \, e=1,\dots,\numel \bigr\},\\
	\testmh(S)    &:=\bigl\{ \widehat{w} \in \eltwo(S): \widehat{w} \vert_{\Gamma_i}\in \Polyk(\Gamma_i) \;\forall \Gamma_i \subset S \subseteq \Gamma \cup \partial \Omega \bigr\} ,
	\end{align}
\end{subequations}
are introduced, where $\mathcal{P}^{k}(\Omega_e)$ and $\mathcal{P}^{k}(\Gamma_i)$ denote the spaces of polynomial functions of complete degree at most $k$ in $\Omega_e$ and on $\Gamma_i$, respectively.
Moreover, let
\begin{subequations}
	\begin{align}
	\testh_t(\Omega)    &:= \eltwo \left( (0,\Tend]; \testh (\Omega) \right),\\
	\testmh_t(S)    &:= \eltwo \left( (0,\Tend]; \testmh (S) \right),
	\end{align}
\end{subequations}
denote the spaces of square-integrable functions on the time interval $(0,\Tend]$ with spatial approximation in $\testh (\Omega)$ and $\testmh (S)$, respectively.\\

Henceforth, the classical notation for $\eltwo$ inner products of vector and tensor-valued functions on a generic subdomain $D \subset\Omega$ is considered, that is
\begin{equation}
\sprod[D]{\bv,\bm{w}} := \int_{D} \bv \cdot \bm{w} \ d\Omega
\quad \text{and} \quad   
\sprod[D]{\bV,\bW} := \int_{D} \bV : \bW \ d\Omega.
\end{equation}
Analogously, the $\eltwo$ inner products on a surface $S \subset \Gamma \cup \partial \Omega$ are denoted by $\dprod[S]{\cdot,\cdot}$.\\

\begin{remark}
	It is worth noticing that the mixed variable $\beps$ requires the definition of an appropriate functional space. In particular, $\beps \in \left[ \Sobo(\text{div};D);\mathbb{S} \right]$, $D \subseteq \Omega$, that is, the space of $\eltwo(D)$ symmetric tensors $\mathbb{S}$ of order $\nsd$ with $\eltwo(D)$ row-wise divergence.
	Accordingly, its element-by-element approximation $\bepsE$ must be defined in an appropriate discrete space for symmetric second-order tensors of dimension $\nsd \times \nsd$. Several approaches have been proposed in the literature, see \cite{Cockburn-CS:2012,Cockburn-Mdecomposition:2017, Qiu-QSS:2017}. In this work, Voigt notation \cite{Fish-Belytschko:2007} is exploited to rearrange the diagonal and off-diagonal components of the tensor into an $\msd$-dimensional vector, being $\msd = \nsd(\nsd + 1)/2$ the number of non-redundant terms. This allows a simple construction of a pointwise symmetric mixed variable with reduced computational cost, while retrieving optimal convergence of the approximation, see~\cite{Giacomini-GKSH:2018,Sevilla-SGKH:2018}. For a detailed derivation of such approach, interested readers are referred to\cite{Tutorial-GSH:2020}.
\end{remark}

With the introduced notation, the discrete weak form associated to the local problems \eqref{eq:NSstrongLocal} is: for every element $\Omega_e$, $e=1,\dotsc, \numel$, find  an approximation $(\buE,\bepsE,\bphiE) \in [\testh_t(\Omega_e)]^{\nsd + 2} \times  [\testh_t(\Omega_e)]^{\msd} \times [\testh_t(\Omega_e)]^{\nsd}$, given $\bhu \in [\testmh_t(\Gamma \cup \partial \Omega)]^{\nsd + 2}$, such that
\begin{subequations} \label{eq:NSweakLocal}
	\begin{align}
	\sprod[\Omega_e]{\bm{\zeta},\bepsE} + \sprod[\Omega_e]{ \Div \left( 2 \bm{\zeta} - \frac{2}{3} \tr (\bm{\zeta}) \Id{\nsd} \right), \bv_{\! e}} - \dprod[\partial \Omega_e]{\left( 2 \bm{\zeta} - \frac{2}{3} \tr (\bm{\zeta}) \Id{\nsd} \right) \bn, \widehat{\bv}} &= 0, \label{eq:NSweakLocalEps}\\ 
	\sprod[\Omega_e]{\bm{\xi},\bphiE} + \sprod[\Omega_e]{\Div \bm{\xi}, T_{\! e}} - \dprod[\partial \Omega_e]{\bm{\xi},\widehat{T}\bn} &= 0, \label{eq:NSweakLocalPhi} \\ 
	\sprod[\Omega_e]{\bw,\pd{\buE}{t}} - \sprod[\Omega_e]{\Grad \bw,\bF - \bG} + \dprod[\partial \Omega_e]{\bw,\reallywidehat{\left( \bF(\buE) - \bG (\buE, \bepsE, \bphiE) \right) \bn}} &=  0, \label{eq:NSweakLocalPrimal}
	\end{align}
\end{subequations}
for all $(\bw,\bm{\zeta},\bm{\xi}) \in  [\testh_t(\Omega_e)]^{\nsd + 2} \times [\testh_t(\Omega_e)]^{\msd} \times [\testh_t(\Omega_e)]^{\nsd}$.

\begin{remark}
	Note that, rigorously, equation~\eqref{eq:NSweakLocalEps} should be derived under the assumption that $\beps$ belongs to the functional space $\left[ \Sobo(\text{div};D);\mathbb{S} \right]$.
	Nonetheless, in an abuse of notation, $\beps$ has been substituted by its discrete counterpart $\bepsE \in  [\testh_t(\Omega_e)]^{\msd}$.
	For further details on the functional spaces and the derivation of the discrete forms, interested readers are referred to~\cite{Tutorial-GSH:2020}.
\end{remark}

Similarly, the discrete weak formulation of the global problem in equation~\eqref{eq:NSstrongGlobal} is: find $\bhu \in [\testmh_t(\Gamma \cup \partial \Omega)]^{\nsd + 2}$ such that
\begin{equation} \label{eq:NSweakGlobal}
\sum_{e = 1}^{\numel} \left\lbrace \dprod[\partial \Omega_e \cap \Gamma]{\widehat{\bw},   \reallywidehat{\left( \bF(\buE) - \bG (\buE, \bepsE, \bphiE) \right) \bn}  } + \dprod[\partial \Omega_e \cap \partial \Omega]{\widehat{\bw}, \bb } \right\rbrace =  0,
\end{equation}
for all $\widehat{\bw} \in [\testmh_t(\Gamma \cup \partial \Omega)]^{\nsd + 2}$.

Equations~\eqref{eq:NSweakLocalPrimal} and~\eqref{eq:NSweakGlobal} introduce the traces of the numerical fluxes on the boundary,
\begin{equation} \label{eq:HDGfluxes}
\reallywidehat{\left( \bF(\buE) - \bG (\buE, \bepsE, \bphiE) \right) \bn} = \reallywidehat{\bF(\buE) \bn} - \reallywidehat{\bG (\buE, \bepsE, \bphiE)  \bn},
\end{equation}
where
\begin{subequations}
	\begin{align}
	&\reallywidehat{\bF(\buE) \bn} :=  \bF (\bhu)  \bn + \btau^a(\bhu) (\buE - \bhu) \quad \text{ and } \label{eq:convNumericalFlux} \\
	&\reallywidehat{\bG (\buE, \bepsE, \bphiE)  \bn} := \bG (\bhu, \bepsE, \bphiE)  \bn - \btau^d (\buE - \bhu) \label{eq:diffNumericalFlux}
	\end{align}
\end{subequations}
stand for the convective and the diffusive numerical fluxes, respectively, whose approximation is essential for the quality and accuracy of the HDG method.

On the one hand, the diffusive numerical fluxes, $\reallywidehat{\bG (\buE, \bepsE, \bphiE)  \bn}$, involve the diffusive stabilisation term $\btau^d$, selected as the diagonal matrix $\btau^d = \Rey^{-1}  \, \text{diag} \left(  0,  \vect{1}_{\nsd}, \left[(\gamma - 1) \Minf^2 \Pra\right]^{-1} \right)$, being $\vect{1}_{\nsd}$ a $\nsd$-dimensional vector of ones.
This approach follows the philosophy of \cite{Peraire-PNC:11,Nguyen-NP:2012} owing to dimensional consistency but considers different amounts of diffusive stabilisation for each of the three conservation equations, i.e., mass, momentum and energy. In particular, note that the continuity equation, which has a purely convective nature, does not include any diffusive stabilisation.\\
It is worth noting that the term $\bG (\bhu, \beps, \bphi)  \bn$ containing the physical flux in~\eqref{eq:diffNumericalFlux} can be approximated either using the interior state $\buE$ or the trace of the primal variable $\bhu$. In this work, the latter has been chosen, following the classical formulation in HDG~\cite{Peraire-PNC:10,Peraire-PNC:11,Nguyen-NP:2012}, which exploits the presence of an intermediate state, namely the trace of the conservation variables, $\bhu$.

\begin{remark}
	Note that in the incompressible limit and using the current choice for mixed variables, both alternatives lead to the same numerical flux.
	Indeed, the energy equation for which the tensor $\bG$ depends on the primal variable $\bu$, is decoupled from the system of conservation equations.
\end{remark}

On the other hand, the convective numerical fluxes, $\reallywidehat{\bF(\buE) \bn}$, are approximated using Riemann solvers \cite{Toro2009}. More precisely, they are introduced implicitly within the numerical fluxes by means of the convective stabilisation parameter, $\btau^a$.
Different definitions of such convective fluxes are detailed in section~\ref{sc:RiemannSolvers}, where a unified framework, including the newly proposed HLL and HLLEM Riemann solvers, is presented in the context of HDG.

\begin{remark}[Compressible Euler equations]
	The associated weak forms for the inviscid Euler equations reduce to:\\
	\textbf{Local problems:} given $\bhu \in [\testmh_t(\Gamma \cup \partial \Omega)]^{\nsd + 2}$ and for every element $\Omega_e$, $e=1,\dotsc, \numel$, find $\buE \in  [\testh_t(\Omega_e)]^{\nsd + 2}$ such that, for all $\bw \in  [\testh_t(\Omega_e)]^{\nsd + 2}$,
	\begin{equation} \label{eq:EulerWeakLocal}
		\sprod[\Omega_e]{\bw,\pd{\buE}{t}} - \sprod[\Omega_e]{\Grad \bw,\bF(\buE)} + \dprod[\partial \Omega_e]{\bw,\reallywidehat{ \bF(\buE)  \bn}} =  0.
	\end{equation}
	\textbf{Global problem:} for all $\widehat{\bw} \in [\testmh_t(\Gamma \cup \partial \Omega)]^{\nsd + 2}$, find $\bhu \in [\testmh_t(\Gamma \cup \partial \Omega)]^{\nsd + 2}$ such that
	\begin{equation} \label{eq:EulerWeakGlobal}
	\sum_{e = 1}^{\numel} \left\lbrace \dprod[\partial \Omega_e \cap \Gamma]{\widehat{\bw},   \reallywidehat{ \bF(\buE) \bn}  } + \dprod[\partial \Omega_e \cap \partial \Omega]{\widehat{\bw}, \bb } \right\rbrace =  0.
	\end{equation}
\end{remark}

\section{A unified framework for Riemann solvers in hybridised discontinuous Galerkin methods}
\label{sc:RiemannSolvers}

As mentioned above, the choice of the convective numerical fluxes $\reallywidehat{\bF(\buE) \bn}$ appearing in equations~\eqref{eq:NSweakLocalPrimal} and~\eqref{eq:NSweakGlobal} ---or in equations~\eqref{eq:EulerWeakLocal} and~\eqref{eq:EulerWeakGlobal} for the Euler equations--- has a critical influence on the accuracy and stability of the numerical solution.
More precisely, such numerical fluxes are responsible for encapsulating the information of the convective nature of the flow under analysis.
For this reason, the approximation of such interface fluxes has received great attention in the context of discontinuous Galerkin methods~\cite{Toro2009,Cockburn-CS:1998,Moura-MMPS:2017,Qiu-QKS:2006} and, more recently, of HDG~\cite{Nguyen-NP:2012,Peraire-PNC:11,Peraire-PNC:10} by means of Riemann solvers.

This section details the expression of numerical fluxes arising in DG discretisations with some of the most popular approximate Riemann solvers for compressible flows, namely Lax-Friedrichs, Roe, HLL and HLLEM. 
Then, a unified framework for the derivation of numerical fluxes in hybridised discretisations is presented. This framework allows to extend the aforementioned numerical fluxes to HDG, including the HLL and HLLEM Riemann solvers, devised for the first time in the context of hybridised formulations.

\subsection{Riemann solvers in standard DG methods}
\label{ssc:RiemannSolversDG}

Consider a pair of neighbouring elements, $\Omega_e^+$ and $\Omega_e^-$, with shared interface $\Gamma_i = \partial \Omega_e^+ \cap \partial \Omega_e^- \subset \Gamma$. 
The solution at each side of the interface is denoted by $\buE^\pm$, whereas $\bu^\star(\buE^+,\buE^-)$ represents an intermediate state between $\buE^+$ and $\buE^-$.
Following the monograph by Toro~\cite{Toro2009}, the definition of Lax-Friedrichs, Roe, HLL and HLLEM Riemann solvers is first recalled for standard DG formulations.

\subsubsection{Lax-Friedrichs Riemann solver}
\label{ssc:LFRiemannsolver}

The first option is represented by the Lax-Friedrichs numerical flux. This Riemann solver is obtained as an extrapolation of the result for a scalar convection equation~\cite{Lax1954} and defines the numerical flux as
\begin{equation} \label{eq:DG-LF}
\reallywidehat{\bF(\buE) \bn^\pm} = \frac{1}{2} \left[ \bF (\buE^+) + \bF (\buE^-) \right] \bn^\pm + \frac{\lamax^\star}{2} (\buE^\pm - \buE^\mp),
\end{equation}
where $\lamax^\star := \abs{\bv^\star\cdot \bn} + c^\star$ is the maximum eigenvalue of the matrix $\mat{A}_n(\bu^\star)$ evaluated at the intermediate state $\bu^\star$.
It is well-known that the Lax-Friedrichs numerical flux~\eqref{eq:DG-LF} is extremely robust but leads to over-diffusive solutions.

\subsubsection{Roe Riemann solver} 
\label{ssc:RoeRiemannsolver}

The Roe Riemann solver~\cite{Roe1981} approximates the complete wave structure of the Riemann problem \cite{Leveque2013,Toro2009} by means of the matrix $\abs{\mat{A}_n(\bu^\star)}$ that linearises the convective fluxes $\bF(\bu^\star)$. More precisely, the Roe numerical flux is given by
\begin{equation} \label{eq:DG-Roe}
\reallywidehat{\bF(\buE) \bn^\pm} = \frac{1}{2} \left[ \bF (\buE^+) + \bF (\buE^-) \right] \bn^\pm + \frac{1}{2}\abs{\mat{A}_n(\bu^\star)}  (\buE^\pm - \buE^\mp),
\end{equation}
where $\mat{A}_n(\bu^\star)$ and $\abs{\mat{A}_n(\bu^\star)}$ are the matrices introduced in section~\ref{ssc:BoundaryConditions} evaluated at the intermediate state $\bu^\star$.

Although more accurate than the Lax-Friedrichs flux, the Roe Riemann solver is not \emph{positivity preserving} and it may produce nonphysical solutions in transonic and supersonic cases due to the violation of entropy conditions~\cite{Quirk1994,Perthame-PS1996}.
In this context, the linearised Roe solver is modified via a so-called \emph{entropy fix} (EF) in order to recover the entropy conditions.
The entropy fix by Harten and Hyman (HH)~\cite{Harten-HH:1983} proposes the following modification of the Roe numerical flux
\begin{equation} \label{eq:DG-RoeEF}
\reallywidehat{\bF(\buE) \bn^\pm} = \frac{1}{2} \left[ \bF (\buE^+) + \bF (\buE^-) \right] \bn^\pm + \frac{1}{2}\abs{\mat{A}_n^\delta(\bu^\star)}   (\buE^\pm - \buE^\mp),
\end{equation}
where $\abs{\mat{A}_n^\delta(\bu^\star)}$ denotes a dissipation matrix.
The HH-EF dissipation matrix is defined as $\abs{\mat{A}_n^\delta(\bu^\star)} := \mat{R} \mat{\Phi} \mat{L}$, being $\mat{R}$ and $\mat{L}$ the right and left eigenvector matrices previously introduced and $\mat{\Phi}$ a diagonal matrix such that $\Phi_{ii} = \max \left( \abs{\lambda_{i}}, \delta \right)$, being $\lambda_i$ the $i$-th eigenvalue of the matrix $\mat{A}_n(\bu^\star)$ and $\delta>0$ a user-defined threshold for the entropy fix.

\begin{remark}
	In the expression of the dissipation matrix, a user-defined threshold parameter $\delta>0$ needs to be appropriately tuned to introduce the correct amount of extra diffusion for the problem under analysis. 
	Note that, generally, $\delta \ll \lamax$. Nonetheless, this value is problem-dependent and may require an empirical tuning to provide the best performance of the Roe solver.
\end{remark}

\subsubsection{Harten-Lax-van Leer (HLL) Riemann solver}

An alternative approach to remedy the entropy violation of the Roe solver is represented by the HLL Riemann solver~\cite{Harten-HLL:1983}.
Such approach relies on a weighted average of the information in two neighbouring elements $\Omega_e^+$ and $\Omega_e^-$ and leads to the following numerical flux
\begin{equation} \label{eq:DG-HLL}
\reallywidehat{\bF(\buE) \bn^\pm} = \left[ \frac{s^+\bF (\buE^+) - s^-\bF (\buE^-)}{s^+ - s^-}\right] \bn^\pm + \frac{s^+ s^-}{s^+ - s^-}(\buE^\pm - \buE^\mp),
\end{equation}
where, respectively, $s^+ := \max \left(0, \bv^\star \cdot \bn^+ + c^\star \right)$ and $s^- := \min \left(0, \bv^\star \cdot \bn^+ - c^\star \right)$ denote the estimates of the largest and smallest wave speeds, with the corresponding signs.

\subsubsection{HLLEM Riemann solver}

Finally, the HLLEM Riemann solver \cite{Einfeldt1988,Einfeldt1991} is introduced as a modification of the HLL Riemann solver which approximates the complete wave structure of the Riemann problem. More precisely, differently from the HLL method, it introduces a special treatment for middle waves, ensuring an accurate description of contact waves and shear layers \cite{Dumbser2016}. In addition, HLLEM inherits the positivity-preserving properties of HLL-type Riemann solvers, fulfilling entropy conditions without the need of the user defined entropy fix required by the Roe solver.

In particular, the HLLEM numerical flux is expressed as
\begin{equation} \label{eq:DG-HLLEM} 
\reallywidehat{\bF(\buE) \bn^\pm} = \left[ \frac{s^+\bF (\buE^+) - s^-\bF (\buE^-)}{s^+ - s^-}\right] \bn^\pm + \frac{s^+ s^-}{s^+ - s^-}\mat{\theta}(\bu^\star) (\buE^\pm - \buE^\mp),
\end{equation}
being $s^+$ and $s^-$ the HLL estimates of the largest and smallest wave speeds previously introduced. In addition, it holds that $\mat{\theta}(\bu^\star) = \mat{R} \mat{\Theta} \mat{L}$, where $\mat{\Theta}$ denotes the diagonal matrix $\mat{\Theta} =  \text{diag} \left(  1,   \theta^\star \vect{1}_{\nsd}, 1 \right)$ and $\theta^\star = \abs{\bv^\star\cdot\bn}/(\abs{\bv^\star\cdot\bn} + c^\star)$ is placed in the position of the eigenvalues corresponding to contact waves. For more details on such matrices, interested readers are referred to~\cite{Rohde2001}.

Note that, in contrast to HLL, the HLLEM flux reduces the amount of numerical dissipation associated to contact waves by means of the coefficient $\theta^\star < 1$. Moreover, it maintains an analogous treatment for shock waves and rarefactions, guarenteeing its entropy enforcement and positivity-preserving properties.

\subsection{Riemann solvers in hybridised DG methods}
\label{ssc:RiemannSolverHDG}

In this section, a unified framework for the formulation of the above introduced Riemann solvers in the context of HDG methods is proposed. The framework includes, for the first time, the formulation of the HLL and HLLEM Riemann solvers within an HDG formulation for compressible flows.
This derivation stems from the seminal work of Peraire and co-workers on linear and nonlinear convection-diffusion equations~\cite{Nguyen-NPC:2009lCD,Nguyen-NPC:2009nlCD} and on compressible flows~\cite{Peraire-PNC:10,Peraire-PNC:11,Nguyen-NP:2012}. The topic has also been studied in~\cite{Bui-Thanh2015}.

As described before, the general structure of the trace of the HDG convective numerical flux for a nonlinear problem is
\begin{equation} \label{eq:HDGConvectiveFlux}
\reallywidehat{\bF(\buE) \bn} = \bF (\bhu)  \bn + \btau^a(\bhu) (\buE - \bhu),
\end{equation}
where $\btau^a$ is the convective stabilisation matrix which encapsulates the information of the Riemann solvers. In order to ease readability, the superindex in $\btau^a$ to denote the convective stabilisation term will be dropped in the upcoming derivations along this section.

It is worth noting that in~\eqref{eq:HDGConvectiveFlux} the hybrid variable $\bhu$ defined on the interface $\Gamma_i$ between two neighboring elements $\Omega_e^+$ and $\Omega_e^-$ is utilised as the intermediate state $\bu^\star$ introduced in section~\ref{ssc:RiemannSolversDG}.

In order to derive the formulation of the Riemann solvers in the context of HDG methods, the inter-element continuity of the trace of the numerical fluxes is considered in the convective limit, namely $\jump{\reallywidehat{\bF(\buE) \bn}}=\bm{0}$. It follows that the sum of the contributions $\reallywidehat{\bF(\buE) \bn}$ from two neighbouring elements is set to zero.
Exploiting definition~\eqref{eq:HDGConvectiveFlux} and observing that $\jump{\bF(\bhu) \bn} = \bm{0}$ because of the uniqueness of $\bhu$ on the internal faces, the above transmission condition reduces pointwise to
\begin{equation}
(\btau^+ + \btau^-) \bhu = \btau^+ \buE^+  + \btau^- \buE^- ,
\end{equation}
where $\btau^+$ and $\btau^-$ denote stabilisation matrices seen from element $\Omega_e^+$ and $\Omega_e^-$, respectively.
Under the assumption of $(\btau^+ + \btau^-)$ being invertible, the intermediate state $\bhu$ is determined pointwise as
\begin{equation}
\bhu = (\btau^+ + \btau^-)^{-1}\left[\btau^+ \buE^+  + \btau^- \buE^-\right] .
\end{equation}

Hence, the convective numerical flux~\eqref{eq:HDGConvectiveFlux} is formulated as an explicit function of the left and right states $\buE^\pm$. 
From the framework above, two cases are analysed hereafter. On the one hand, a stabilisation matrix continuous across the interface is obtained by setting $\btau^+ = \btau^-$. On the other hand, a stabilisation matrix, discontinuous across the interface, is considered when $\btau^+ \neq \btau^-$.

\subsubsection{Continuous stabilisation across the interface: Lax-Friedrichs and Roe Riemann solvers} 

Consider a continuous definition of the stabilisation matrix across the interface, that is $\btau^+ = \btau^- = \overline{\btau}$. It follows
\begin{subequations} \label{eq:hDGfluxContinuous}
	\begin{align} 
	\bhu &= \frac{\buE^+  + \buE^-}{2} , 
	\label{eq:uHatStateContinuous} \\
	\reallywidehat{\bF(\buE) \bn^\pm} &= \bF \left( \frac{\buE^+  + \buE^- }{2} \right)  \bn^\pm + \frac{1}{2} \overline{\btau}  (\buE^\pm - \buE^\mp).
	\label{eq:FluxContinuous}
	\end{align}
\end{subequations}

By considering $\bhu$ as an intermediate state between $\buE^+$ and $\buE^-$ and under appropriate choices of the stabilisation matrix $\overline{\btau}$, a formulation that mimicks Lax-Friedrichs and Roe Riemann solvers for DG methods, see~\eqref{eq:DG-LF} and~\eqref{eq:DG-Roe}, is retrieved for HDG methods~\cite{Peraire-PNC:10,Peraire-PNC:11,Nguyen-NP:2012}.
More precisely, for each element $\Omega_e, \ e=1,\ldots,\numel$, setting $\overline{\btau} = \hlamax \bmat{I}_{\nsd + 2}$, with $\hlamax := \abs{\bhv\cdot \bn} + \widehat{c}$, the Lax-Friedrichs numerical flux is retrieved for the HDG method, namely
\begin{equation} \label{eq:LFRiemannsolver}
\reallywidehat{\bF(\buE) \bn} = \bF (\bhu) \bn + \hlamax  (\buE - \bhu) .
\end{equation}

Similarly, the intermediate state~\eqref{eq:uHatStateContinuous} and the stabilisation matrix $\overline{\btau} = \abs{\mat{A}_n(\bhu)}$ lead to the formulation of the Roe Riemann solver in the context of HDG methods, that is,
\begin{equation} \label{eq:RoeRiemannsolver}
\reallywidehat{\bF(\buE) \bn} = \bF (\bhu) \bn + \abs{\mat{A}_n(\bhu)}  (\buE - \bhu).
\end{equation}
Finally, the HH-EF variant of the Roe numerical flux is given by $\overline{\btau} = \abs{\mat{A}_n^\delta(\bhu)}$, according to the correction to matrix $\mat{A}_n(\bhu)$ introduced in~\eqref{eq:DG-RoeEF}.

\begin{remark}
	It is worth noting that the stabilisation matrix introduced in~\eqref{eq:LFRiemannsolver} for the Lax-Friedrichs Riemann solver is isotropic, whereas for the Roe numerical fluxes in~\eqref{eq:RoeRiemannsolver}, different values of the stabilisation term are introduced in the equations of conservation of mass, momentum and energy.
\end{remark}

\subsubsection{Discontinuous stabilisation across the interface: HLL-type Riemann solvers}

Consider a discontinuous stabilisation matrix across the interface, defined as $\btau^\pm = s^\pm \mat{\theta}$, with $s^+ \neq s^-$ and $\mat{\theta}$ a positive-definite square matrix of dimension $\nsd+2$. It follows
\begin{subequations} \label{eq:hDGfluxDiscontinuous}
	\begin{align} 
	\bhu &= \frac{s^+ \buE^+  + s^- \buE^- }{s^+ + s^-} , 
	\label{eq:uHatStateDiscontinuous}\\
	\reallywidehat{\bF(\buE) \bn^\pm} &= \bF \left( \frac{s^+ \buE^+  + s^- \buE^- }{s^+ + s^-} \right)  \bn^\pm + \frac{s^+ s^- }{s^+ + s^-}  \mat{\theta} (\buE^\pm - \buE^\mp).
	\end{align}
\end{subequations}

It is worth noting that the intermediate state in~\eqref{eq:uHatStateDiscontinuous} is obtained as a weighted average of the states $\buE^+$ and $\buE^-$.
From this framework, HLL-type numerical fluxes, mimicking the behaviour of HLL~\eqref{eq:DG-HLL} and HLLEM~\eqref{eq:DG-HLLEM} for DG approaches, are devised for the first time in the context of HDG methods. More precisely, the HLL Riemann solver is given by
\begin{equation} \label{eq:HLLRiemannsolver}
\reallywidehat{\bF(\buE) \bn} = \bF (\bhu) \bn + \left[ s^+ \bmat{I}_{\nsd+2}\right]  (\buE - \bhu),
\end{equation}
where $s^+ := \max (0,\bhv\cdot \bn + \widehat{c})$.

\begin{remark}
	A variant of the HLL Riemann solver in~\eqref{eq:HLLRiemannsolver}, the so-called Harten-Lax-van Leer-Einfeldt (HLLE) numerical flux~\cite{Einfeldt1988}, can be devised by simply modifying the term $s^+$ in the stabilisation parameter as
	\begin{equation} \label{eq:HLLEstabilisation}
	s^+ := \max (0,\bhv\cdot \bn + \widehat{c},\bv^+\cdot \bn + c^+,\bv^-\cdot \bn + c^-),
	\end{equation}
	being $\circledcirc^+$ and $\circledcirc^-$ the variables associated with the states $\buE^+$ and $\buE^-$, respectively, at each side of the interface under analysis.
	Numerical experiments have shown that, in the context of high-order discretisations, the practical difference between HLL and HLLE numerical fluxes is not significant since  the jumps across the interface are very small. Henceforth, the former choice is considered for simplicity.
\end{remark}

Following the same rationale, the HLLEM numerical flux can be devised as
\begin{equation} \label{eq:HLLEMRiemannsolver}
\reallywidehat{\bF(\buE) \bn} = \bF (\bhu) \bn +  \left[ s^+ \mat{\theta} (\bhu) \right]  (\buE - \bhu),
\end{equation}
where $s^+ := \max (0,\bhv\cdot \bn + \widehat{c})$ is the HLL estimate for the largest wave speed and $\mat{\theta} (\bhu) = \mat{R} \mat{\Theta} \mat{L}$, as defined in~\eqref{eq:DG-HLLEM}. It is worth noticing that the intermediate state is selected such that $\bu^\star = \bhu$. Therefore, $\mat{\Theta}$ employs $\theta^\star = \widehat{\theta} = \abs{\bhv\cdot\bn}/(\abs{\bhv\cdot\bn} + \widehat{c})$, where the hat quantities are evaluated using the hybrid variable $\bhu$.

\begin{remark}
	Because of the positive definition of the matrix $\mat{\theta}$ introduced here, the coefficient $\widehat{\theta}$ is not allowed to reach zero. This situation is experienced in flows that are perfectly aligned with the faces of the mesh. From a practical point of view, it may be useful to set a minimum threshold $0 < \theta_0 \ll 1$ to guarantee that $\widehat{\theta} > \theta_0$ and avoid a null stabilisation.
\end{remark}

\section{Implementation details of the high-order HDG solver}
\label{sc:Implementation}

In this section, some details on the implementation of the nonlinear solver in the high-order HDG method and on the numerical treatment of solutions with discontinuities and sharp gradients are provided.

\subsection{Solution strategy}
\label{ssc:solutionStrategy}

By introducing the numerical flux~\eqref{eq:HDGfluxes} and boundary conditions (detailed in table~\ref{tb:boundaryConditions}) in the weak forms of the local~\eqref{eq:NSweakLocal} and global~\eqref{eq:NSweakGlobal} problems, the complete form of the discrete problems is obtained.

It is worth recalling that the HDG solver features two stages. First, the local problems are devised.
Denote by $\bznode = ( \buE, \bepsE, \bphiE ) \in [\testh_t(\Omega_e)]^{\nsd + 2} \times  [\testh_t(\Omega_e)]^{\msd} \times [\testh_t(\Omega_e)]^{\nsd}$ the vector of local unknowns, which includes the primal and mixed variables.
By considering an isoparametric approximation in space for the local, $\bznodeAll$, and hybrid, $\bhu$, variables, the semi-discrete system of differential-algebraic equations resulting from the local problem at element $\Omega_e, \ e = 1,\dotsc,\numel$ reads
\begin{equation} \label{eq:LocalSemiDiscreteCompact} 
\bmat{M}_e \frac{d \bzE}{d t} + \re(\bzE,\hu) = \bmat{0}.
\end{equation}
where $\bzE$ and $\hu$ denote the vectors of nodal values of the local and hybrid variables, respectively, and $\bmat{M}_e$ and $\re$ are the mass matrix and nonlinear residual vector obtained from the spatial discretisation of the integral terms of the local problem~\eqref{eq:NSweakLocal} in element $\Omega_e$.

In a similar fashion, from the global problem~\eqref{eq:NSweakGlobal} it follows
\begin{equation} \label{eq:GlobalSemiDiscreteCompact} 
\sum_{e = 1}^{\numel} \hre(\hu,\bzE) = \bmat{0},
\end{equation}
where $\hre$ denotes the nonlinear residual vector involving the degrees of freedom associated with element $\Omega_e$.

Finally, upon temporal discretisation, the resulting nonlinear system is solved using a Newton-Raphson iterative method at each time step.
In particular, the linear system of equations arising at each time step and Newton-Raphson iteration for the local problems reads
\begin{equation} \label{eq:LocalDiscreteCompact} 
\AmatE{Z}{Z} \bzE + \AmatE{Z}{\widehat{U}} \hu = \FvectE{Z}
\end{equation}
for $e = 1,\dotsc,\numel$, where vectors $\FvectE{\diamond}$ and matrices $\AmatE{\diamond}{\circ}$ are obtained from Newton-Raphson linearisation of the system of equations~\eqref{eq:LocalSemiDiscreteCompact}.
Similarly, the linear system corresponding to the global problem~\eqref{eq:GlobalSemiDiscreteCompact} upon Newton-Raphson linearisation can be expressed as
\begin{equation} \label{eq:GlobalDiscreteCompact} 
\sum_{e = 1}^{\numel} \left\lbrace \AmatE{\widehat{U}}{\widehat{U}}\hu + \AmatE{\widehat{U}}{Z} \bzE - \FvectE{\widehat{U}}  \right\rbrace = \bmat{0}.
\end{equation}

Note that, owing to the hybridisation procedure, the elemental degrees of freedom of $\bzE$ can be rewritten in terms of the globally coupled degrees of freedom of $\hu$ via~\eqref{eq:LocalDiscreteCompact}, namely
\begin{equation} \label{eq:LocalSolution} 
\bzE = \left[ \AmatE{Z}{Z} \right]^{-1} \FvectE{Z} - \left[ \AmatE{Z}{Z} \right]^{-1} \AmatE{Z}{\widehat{U}} \hu,
\end{equation}
which just involves the inverse of matrix $\AmatE{Z}{Z}$, of dimension $\left((\nsd + 2 + \msd + \nsd) \nen \right)$, for each element of the mesh, being $\nen$ the number of element nodes of $\Omega_e$. This computation can be effectively parallelised and only involves the solution of small systems with limited computing effort.
Dimension of such local systems is displayed in table~\ref{tb:dimensionLocalSystem} for different degrees of approximation $k$ on simplexes and parallelepipeds in 2D and 3D.

\begin{table} [htbp]
	\caption{Dimension of the local problem.}
	\centering
	\makebox[\linewidth]{
		\begin{tabular}{ L{6.5cm}  L{1cm}  L{1cm} L{1cm} L{1cm} L{1cm} L{1cm}}
			\toprule
			Degree of approximation, $k$ & 1 & 2 & 3 & 4 & 5 & 6\\
			\midrule
			\multicolumn{7}{l}{Simplexes} \\ 	\hline
			2D & 27 & 54 & 90 & 135 & 189 & 252\\
			3D & 56 & 140 & 280 & 490 & 784 & 1,176\\
			\midrule
			\multicolumn{7}{l}{Parallelepipeds} \\ 	\hline
			2D & 36 & 81 & 144 & 225 & 324 & 441\\
			3D & 112 & 378 & 896 & 1,750 & 3,024 & 4,802\\	
			\bottomrule
		\end{tabular}
	}
	\label{tb:dimensionLocalSystem}
\end{table}

The hybridisation precedure~\eqref{eq:LocalSolution} permits to eliminate $\bzE$ in equation~\eqref{eq:GlobalDiscreteCompact}, giving rise to a linear system with a reduced number of degrees of freedom \cite{Cockburn2016,Guyan1965}. This global system is the one to be solved at each Newton-Raphson iteration and reads as
\begin{equation} \label{eq:GlobalSystem} 
\Kmat \hu = \Fvect,
\end{equation}
where the global matrix $\Kmat$ and the right-hand side vector $\Fvect$ are obtained by assembling the elemental contributions
\begin{subequations} \label{eq:GlobalSystemTerms}
	\begin{align}
	\Kmat^e &= \AmatE{\widehat{U}}{\widehat{U}} - \AmatE{\widehat{U}}{Z} \left[ \AmatE{Z}{Z} \right]^{-1} \AmatE{Z}{\widehat{U}} , \\
	\Fvect^e &= \FvectE{\widehat{U}} - \AmatE{\widehat{U}}{Z} \left[ \AmatE{Z}{Z} \right]^{-1} \FvectE{Z} .
	\end{align}
\end{subequations}

As the main purpose of this work is the HDG formulations of compressible flows, section~\ref{sc:Benchmarks} only considers steady state problems. In this context, the temporal discretisation in equation~\eqref{eq:LocalSemiDiscreteCompact} is used as a relaxation method to improve the convergence process in complex numerical examples, e.g. in presence of shocks. To this effect, the backward Euler method is considered in the simulations. 
However, the proposed methodology is applicable to other time discretisations, such as high-order time integrators like backward difference formulas (BDF) or diagonally implicit Runge-Kutta (DIRK) methods, especially suited for transient problems~\cite{Nguyen-NP:2012,Jaust-JSW:2015,Jaust-JS:2014}.

\subsection{Shock-capturing method}
\label{ssc:shockCapturing}

It is well-known that high-order methods experience an oscillatory behaviour in the vicinity of shocks and regions with sharp gradients, requiring an appropriate shock-capturing technique~\cite{Godunov1959,Donea2003}.
For this purpose, an artificial viscosity term is added to regularise the numerical approximation of the problem.

Different approaches can be adopted to introduce artificial dissipation. In this section, two different alternatives are presented. First, a physics-based shock capturing term, which is introduced within the viscous flux $\bG$, is detailed~\cite{Fernandez-FNP:2017}. Additionally, a Laplacian-based approach~\cite{Jaust-JSW:2015}, formulated in a discrete version to avoid the introduction of mixed variables, is considered for the Euler equations.

\subsubsection{Physics-based shock capturing}
\label{ssc:PhysicsBasedShockCapturing}

In this approach, shock waves are stabilised by correcting the diffusive flux in equation~\eqref{eq:NScompact} using the physics-based approach proposed in \cite{Fernandez-FNP:2018}. This methodology, stemming from the work of Von Neumann and Richtmyer \cite{VonNeumann1950} and later considered in \cite{Cook2005,Kawai2008,Kawai2010}, relies on defining the diffusive flux as a combination of the physical flux $\bG$ with an additional numerical contribution $\bG^*$. The latter is thus based on an artificial bulk viscosity $\beta^*$, namely
\begin{equation}\label{eq:BulkViscosityFlux}
\bG^* = \beta^*\begin{bmatrix} \vect{0}\\
(\Div \bv) \Id{\nsd} \\
\left[(\Div \bv) \bv + \Pra_{\beta}^{-1} \bphi \right]\tras \end{bmatrix},
\end{equation}
where $\Pra_{\beta}$ is an artificial Prandtl number.

First, a dilatation-based shock sensor~\cite{Moro-MNP:2016,Fernandez-FNP:2018}, which identifies the regions of high compression, is defined as
\begin{equation}\label{eq:sensorDilatation}
s_{\beta} = - \frac{h}{k} \frac{\Div \bv}{\tilde{c}},
\end{equation}
where $h$ is the element size, $k$ is the degree of polynomial approximation and $\tilde{c}$ is a reference speed of sound for non-dimensionalisation. Common choices for $\tilde{c}$ are the critical speed of sound $c^\star$, the speed of sound at the actual point $c$, or simply the reference free-stream value $c_\infty$. In the simulations presented in section~\ref{sc:Benchmarks}, the latter option is employed.

The shock sensor $s_{\beta}$ is thus utilised to define the artificial bulk viscosity $\beta^*$ as
\begin{equation}\label{eq:artBulkViscosity}
\beta^* = \Psi \left(  \varepsilon_0 \left[\rhoinf  \frac{h}{k} (\vinf^2 + c_{\infty}^2)^{1/2} \right] f_\beta(s_{\beta}) \right),
\end{equation}
where $\Psi$ denotes a smoothing operator consisting of a $\Czero$ reconstruction, see \cite{Persson2013}, $\varepsilon_0$ is a user-defined positive constant and $f_\beta(s_{\beta}) = \min \left\lbrace s_{\max}, \max \lbrace s_{\min}, s_{\beta} - s_0 \rbrace \right\rbrace$.
Following~\cite{Fernandez-FNP:2018} the values $\varepsilon_0 = 1.5$, $s_0 = 0.01$, $s_{\min} = 0$ and $s_{\max} = 2/\sqrt{\gamma^2 - 1}$ and the artificial Prandtl number $\Pra_{\beta} = 0.9$ are employed in the numerical simulations of section~\ref{sc:Benchmarks}.

\subsubsection{Laplacian-based shock capturing}
\label{ssc:LaplacianBasedShockCapturing}

The second alternative for the shock capturing detailed in this section consists of a discretised Laplace operator, applied in HDG discretisations~\cite{Jaust-JSW:2015,Jaust-JSW:2014} following standard approaches in the context of DG and SUPG methods~\cite{Bassi-BR:1995,Cockburn2001,Casoni-CPH:2012,Sevilla-SHM:2013}. Given the artificial viscosity $\varepsilon$, it relies on adding the term 
\begin{equation}
\sprod[\Omega_e]{\Grad \bw, \varepsilon \Grad \bu}
\end{equation}
to the left-hand side of the local equation~\eqref{eq:NSweakLocalPrimal}, or~\eqref{eq:EulerWeakLocal} for the Euler case.
This approach is especially suited for the inviscid case, where the second-order term $\bG$ vanishes and the mixed variables in~\eqref{eq:NSweakLocalEps} and~\eqref{eq:NSweakLocalPhi} are neglected.

The shock capturing technique is equipped with a discontinuity sensor $S_e$, introduced in~\cite{Persson-PP:2006} and expressed in terms of the density field according to~\cite{Persson2013}, namely
\begin{equation} \label{eq:sensorPersson}
S_e := \frac{\sprod[\Omega_e]{\rho_e - \widetilde{\rho}_e,\rho_e - \widetilde{\rho}_e}}{\sprod[\Omega_e]{\rho_e, \rho_e}}.
\end{equation}
The smoothness indicator $S_e$ is utilised to detect the regions with discontinuities. In~\eqref{eq:sensorPersson}, $\rho_e$ denotes the density in the element $\Omega_e$, computed using a polynomial approximation of degree $k$, and $\widetilde{\rho}_e$ is its truncation of order $k-1$.
The sensor measures the regularity of the approximate solution based on the rate of decay of its Fourier coefficients. More precisely, if $S_e > k^{-4}$, such approximation is expected to be at most $\Czero$, whereas smooth functions are expected to decay more rapidly~\cite{Casoni-CPH:2012}. 

Following~\cite{RSC-SFH:08,Huerta-HCP:2011}, the sensor~\eqref{eq:sensorPersson} is implemented using nodal basis functions. It follows that
\begin{equation}
S_e = \frac{\vect{\rho}_e\tras \bmat{V}^{-T} \bmat{P} \bmat{V}^{-1} \vect{\rho}_e}{\vect{\rho}_e\tras \bmat{V}^{-T} \bmat{V}^{-1} \vect{\rho}_e},
\end{equation}
where $\vect{\rho}_e$ is the vector containing the nodal values of the density field in the element $\Omega_e$, $\bmat{V}$ is the Vandermonde matrix whose inverse maps the Lagrange basis onto the orthonormal one and $\bmat{P}$ is the orthogonal projection matrix onto the space of monomials of degree $k$, namely
\begin{equation}
\bmat{P} := \text{diag} ( \overbrace{0, \dotsc, 0}^{\texttt{n}_{\texttt{L}}}, \overbrace{1, \dotsc, 1}^{\texttt{n}_{\texttt{H}}} ),
\end{equation}
being $\texttt{n}_{\texttt{L}}$ and $\texttt{n}_{\texttt{H}}$ the number of degrees of freedom for monomials of degree $k-1$ and $k$, respectively. 
In two dimensions, it holds $\texttt{n}_{\texttt{L}} := k + 1$ and $\texttt{n}_{\texttt{H}} := k(k+1)/2$.

The amount of artificial viscosity introduced in each element is determined according to
\begin{align} \label{eq:ArtificialViscosity}
\varepsilon_e =  \begin{cases}
0, & \text{if } s_e < s_0 - \xi,\\
\displaystyle \frac{\varepsilon_0}{2} \left( 1 + \sin \left( \frac{\pi (s_e-s_0)}{2 \xi} \right) \!\! \right), & \text{if } s_0 - \xi < s_e < s_0 + \xi,\\
\varepsilon_0, & \text{if } s_e > s_0 + \xi,
\end{cases}
\end{align}
where $s_e := \log_{10} S_e$, $\varepsilon_0 \sim h/k$ and $s_0$ and $\xi$ are selected such that $s_0 + \xi = -4 \log_{10} k$ and $s_0 - \xi$ is sufficiently large to detect the regions in which mild shock waves are present~\cite{Huerta-HCP:2011}. 
In particular, a value $s_0 - \xi = -11 \log_{10} k$ is considered in the numerical studies in section~\ref{sc:Benchmarks}.
Finally, the smoothing operator $\Psi$ is employed to perform a $\mathcal{C}^0$ reconstruction of the elemental artificial viscosity obtained in~\eqref{eq:ArtificialViscosity}, that is $\varepsilon = \Psi(\varepsilon_e)$.

\section{Numerical convergence studies}
\label{sc:Convergence}

The optimal convergence properties of the HDG method are tested both in inviscid and viscous cases. The accuracy of the approximation is examined using the four Riemann solvers presented in section~\ref{sc:RiemannSolvers} for different degrees of polynomial approximation.

\subsection{Convergence analysis for inviscid flows: Ringleb flow}
\label{ssc:Ringleb}

The Ringleb flow problem is considered to verify the optimal convergence of the HDG method for inviscid flows. It consists of a smooth transonic 2D solution of the Euler equations with analytical expression obtained via the hodograph method~\cite{Chiocchia1985}. For any given spatial coordinates $(x,y)$, the solution of the Ringleb flow can be computed by solving the following nonlinear implicit equation in terms of the speed of sound $c$,
\begin{equation}\label{eq:RinglebNLeq}
\left( x + \frac{J}{2} \right)^2 + y^2 = \frac{1}{4 \rho^2 V^4},
\end{equation}
where the following relationships for density $\rho$, radial velocity $V$ and $J$ hold
\begin{equation}
\rho = c^{2/(\gamma - 1)}, \quad V = \sqrt{\frac{2(1 - c^2)}{\gamma - 1}}, \quad J = \frac{1}{c} + \frac{1}{3c^3} + \frac{1}{5c^5} - \frac{1}{2}\log\left( \frac{1+c}{1-c}\right).
\end{equation}
The exact velocity and pressure fields are
\begin{equation}\label{eq:RinglebVandP}
\bv = \begin{pmatrix} -  \sign{y} V \sin \theta\\ V \cos \theta 	\end{pmatrix}
\quad \text{and} \quad
p = \frac{1}{\gamma} c^{2\gamma/(\gamma - 1)} ,
\end{equation}
where $\sign{\cdot}$ is the \emph{sign} operator, $\sin \theta := \Psi V$ and
\begin{equation}\label{eq:RinglebPsi}
\Psi := \sqrt{\frac{1}{2 V^2} + \rho \left( x + \frac{J}{2}\right)}.
\end{equation}

\begin{remark}[Computation of the Ringleb solution]
	\label{rk:RinglebPrecision}
	It is worth noting that the nonlinear equation driving the analytical solution of the Ringleb problem~\eqref{eq:RinglebNLeq} needs to be solved iteratively upon a certain tolerance, thus introducing an approximation error in the estimated analytical solution. Further operations in order to compute the rest of variables of the problem may be responsible for the propagation of such error, which may become critical in high-order convergence tests. In these cases, the error introduced in the exact solution may be of similar level or even higher than the error of the approximate solution. Then, the computed approximation error is no longer reliable, showing a stagnation in the levels of accuracy.
	
	Such numerical issues were circumvented in this study by avoiding the computation of $\Psi$ directly as in~\eqref{eq:RinglebPsi} but using trigonometric identities and algebraic manipulation of~\eqref{eq:RinglebNLeq} to compute the direction of the flow, namely	
	\begin{equation} 
		\sin (2\theta) = 2 \sin \theta \cos \theta = 2 \Psi V \sqrt{1 - \Psi^2 V^2} = 2 \sqrt{\frac{1}{4} - \rho^2 V^4 \left( x + J/2\right)^2} = 2 \rho V^2 y.
	\end{equation}
\end{remark}

\begin{remark}[Domain of the Ringleb solution]
	Classicaly, the Ringleb flow problem has been solved in a curvilinear domain symbolising a channel around a symmetric blunt obstacle, bounded by two streamlines of the flowfield, see~\cite{Bassi-BR:97,HartmannHouston:2003,WangLiu:2006,Dumbser-DKTT:2007,Vymazal2015}. In such domain, the flow is transonic, displaying a large supersonic region near the nose of the blunt body.
	Alternatively, this problem has also been studied in rectangular domains located at different regions, thus avoiding the introduction of geometric errors in the approximation of the curved boundaries. Numerical tests performed both in regions of subsonic~\cite{Nguyen-NP:2012} or transonic~\cite{Sanjay2020} speeds have been presented in the literature.
	It is worth mentioning that the numerical issues described in remark~\ref{rk:RinglebPrecision} may be more evident in those regions where the solution displays greater variations, namely those including supersonic speeds.
\end{remark}

\begin{figure}[htbp]
	\subfloat[Mesh 1 \label{fig:Ringleb_mesh1}]{\includegraphics[width=0.31\textwidth]{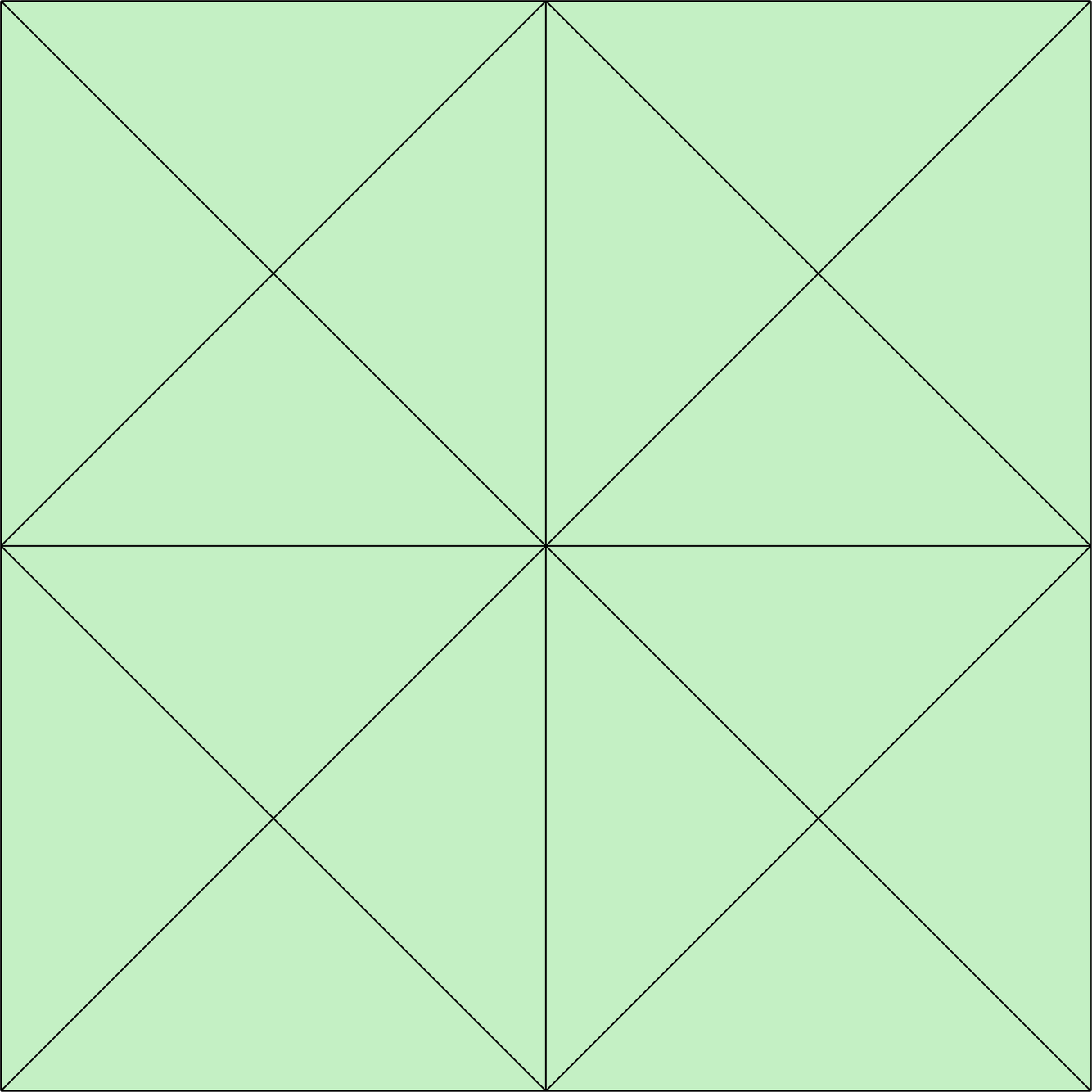}} \quad
	\subfloat[Mesh 2 \label{fig:Ringleb_mesh2}]{\includegraphics[width=0.31\textwidth]{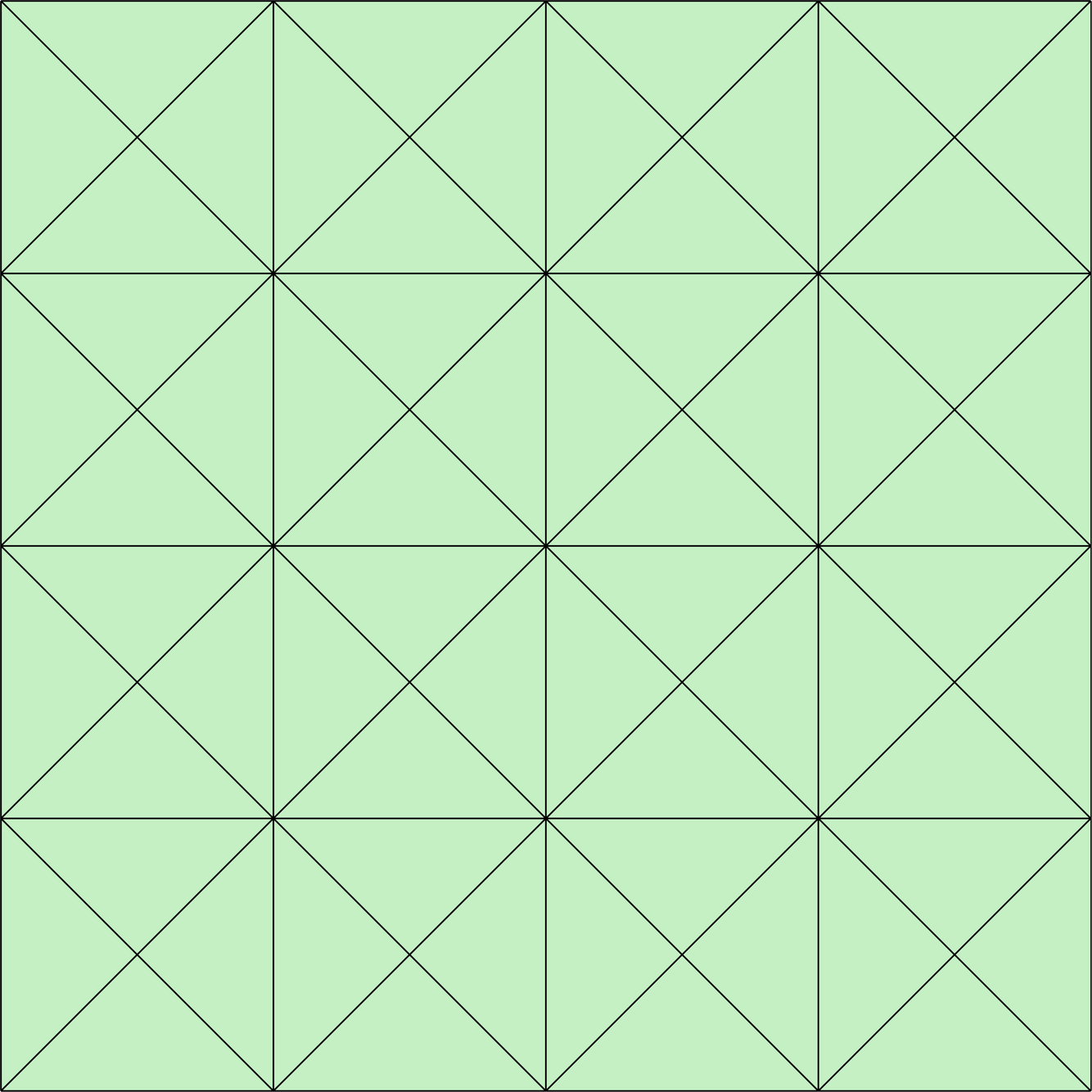}} \quad
	\subfloat[Mesh 3 \label{fig:Ringleb_mesh3}]{\includegraphics[width=0.31\textwidth]{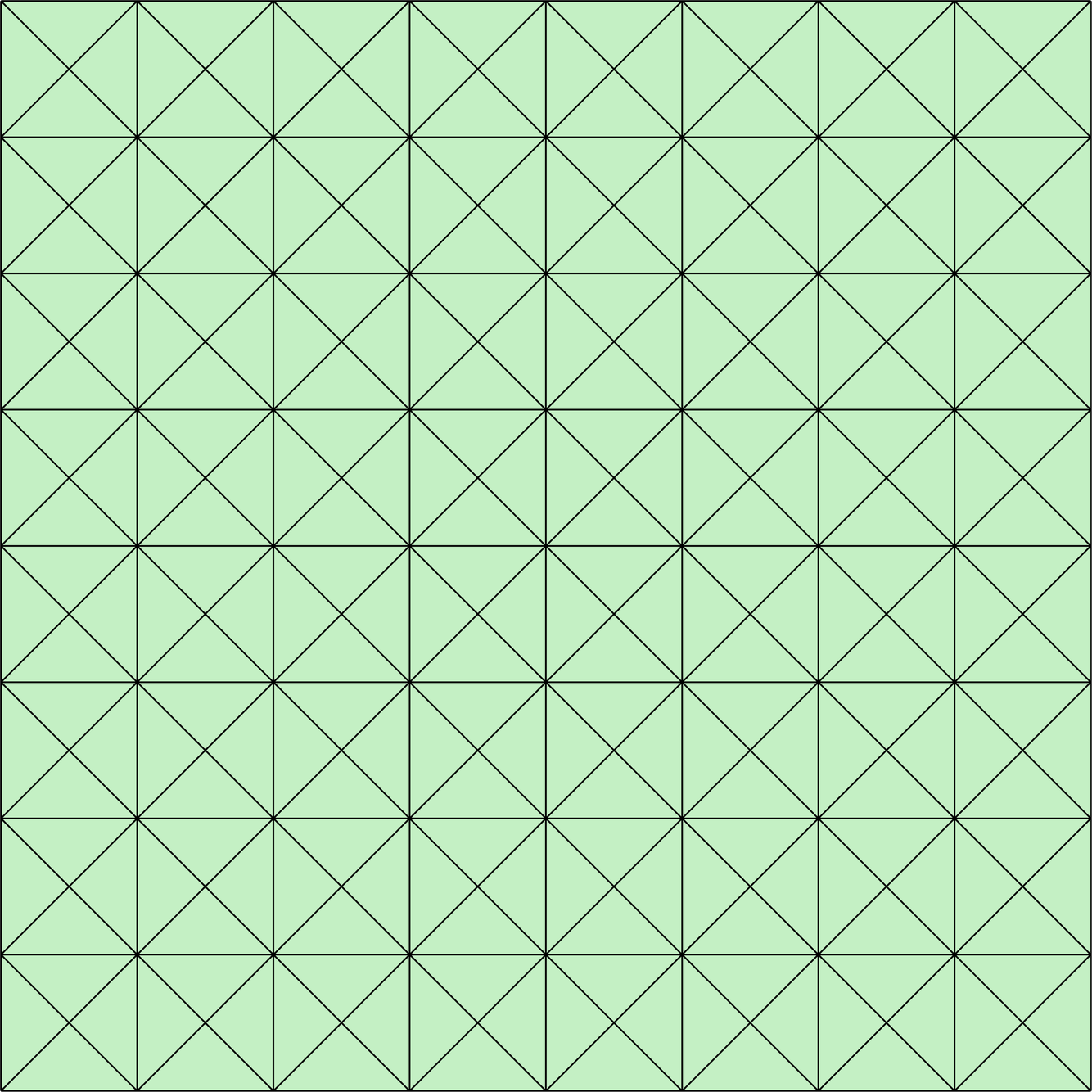}}
	\caption{Ringleb flow - Triangular meshes of $\Omega = [0,1]^2$ for the $h$-convergence analysis.}
	\label{fig:Ringleb_meshes}
\end{figure}

In this work, the Ringleb problem is solved in the domain of transonic flow $\Omega = [0,1]^2$, such as in~\cite{Sanjay2020}, with a far-field boundary condition imposed on $\partial \Omega$. The computational domain is discretised using uniform meshes of triangular elements. Figure~\ref{fig:Ringleb_meshes} displays the first three levels of refinement employed.

\begin{figure}[htbp]
	\subfloat[$k=1$]{\includegraphics[width=0.32\textwidth]{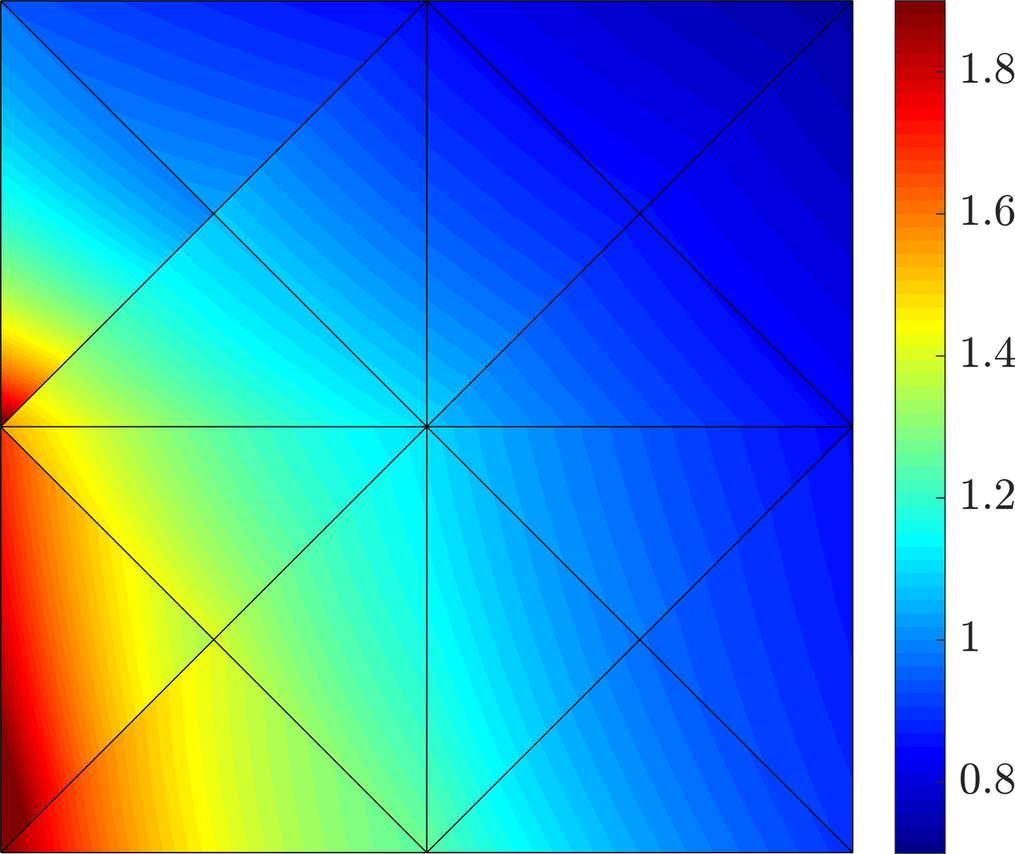}} \hfill
	\subfloat[$k=2$]{\includegraphics[width=0.32\textwidth]{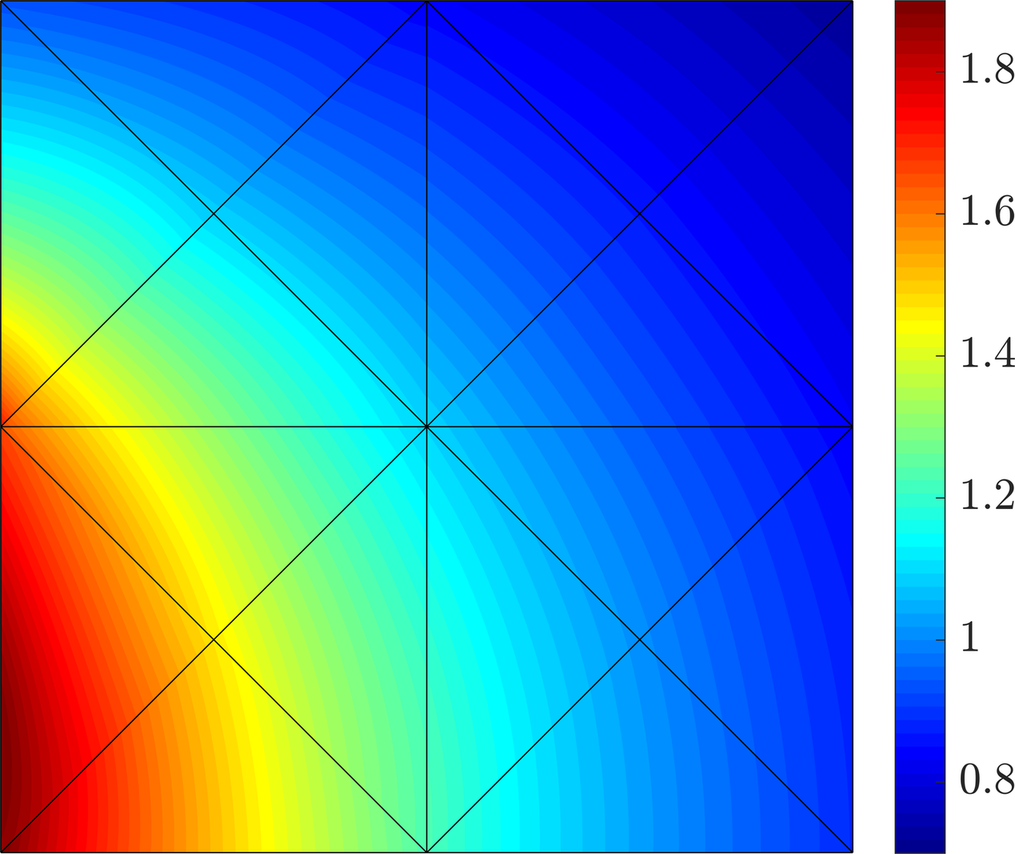}} \hfill
	\subfloat[$k=3$]{\includegraphics[width=0.32\textwidth]{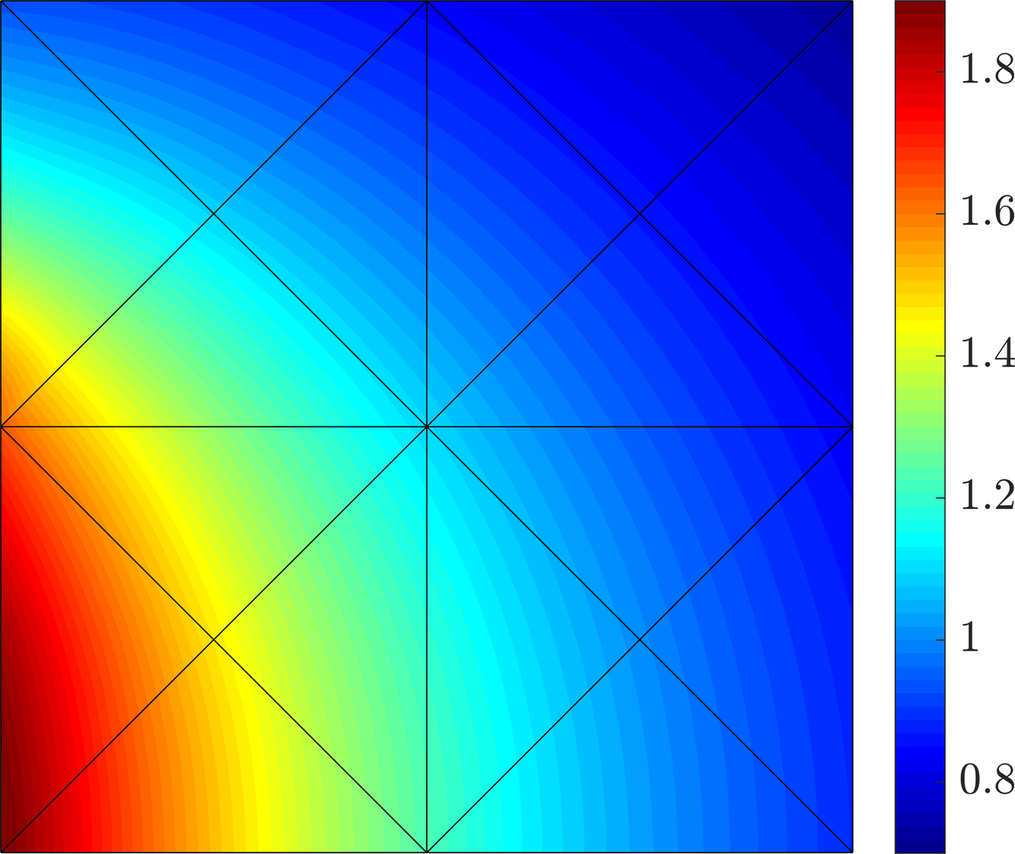}}
	\caption{Ringleb flow - Mach number distribution computed using the HLL Riemann solver on the first level of mesh refinement with polynomial degree $k=1,\ldots,3$.}
	\label{fig:Ringleb_Mach}
\end{figure}

The approximate solution of the Mach number distribution computed on the mesh in figure~\ref{fig:Ringleb_mesh1} using polynomial degree $k=1,\ldots,3$ is depicted in figure~\ref{fig:Ringleb_Mach}. The results clearly display the gain in accuracy obtained increasing the degree of the polynomial approximation, even in presence of extremely coarse meshes, motivating the interest in high-order discretisations.

\begin{figure}[htbp]
	\centering
	\subfloat[Density, $\rho$ \label{fig:Ringleb_density}] {\includegraphics[width=0.33\textwidth]{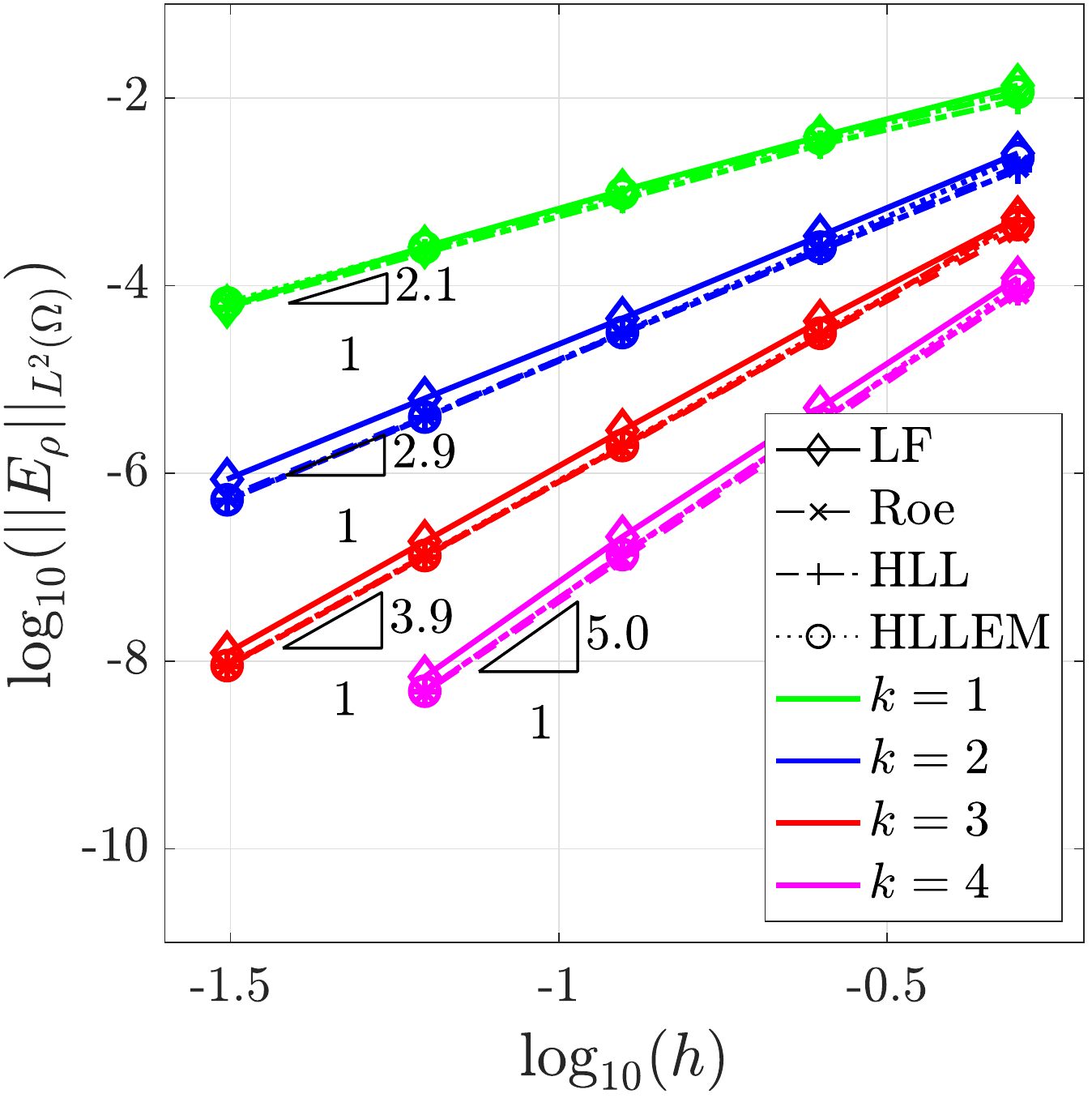}} \hfill
	\subfloat[Momentum, $\rho \bv$ \label{fig:Ringleb_momentum}] {\includegraphics[width=0.33\textwidth]{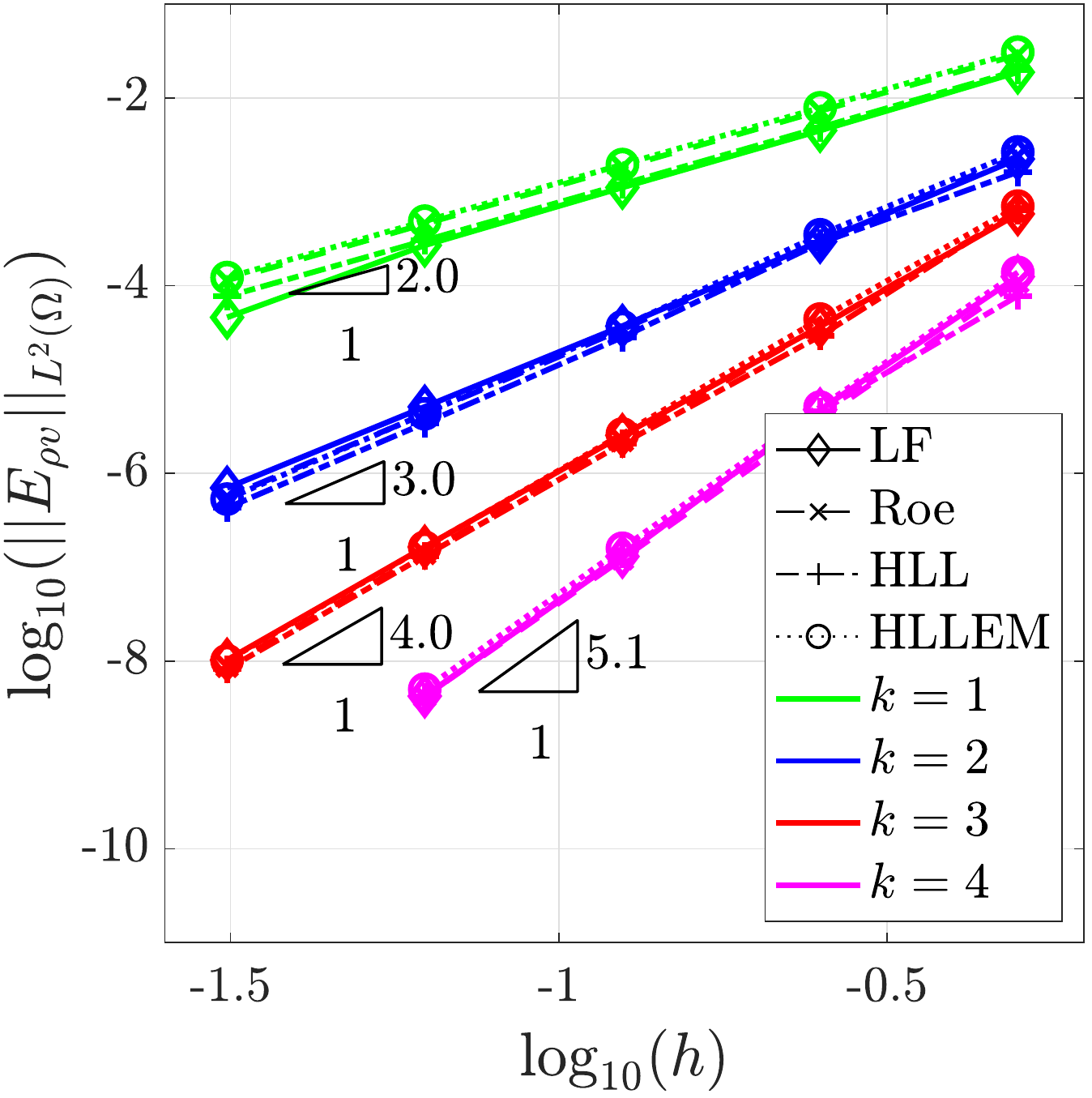}} \hfill
	\subfloat[Energy, $\rho E$ \label{fig:Ringleb_energy}] {\includegraphics[width=0.33\textwidth]{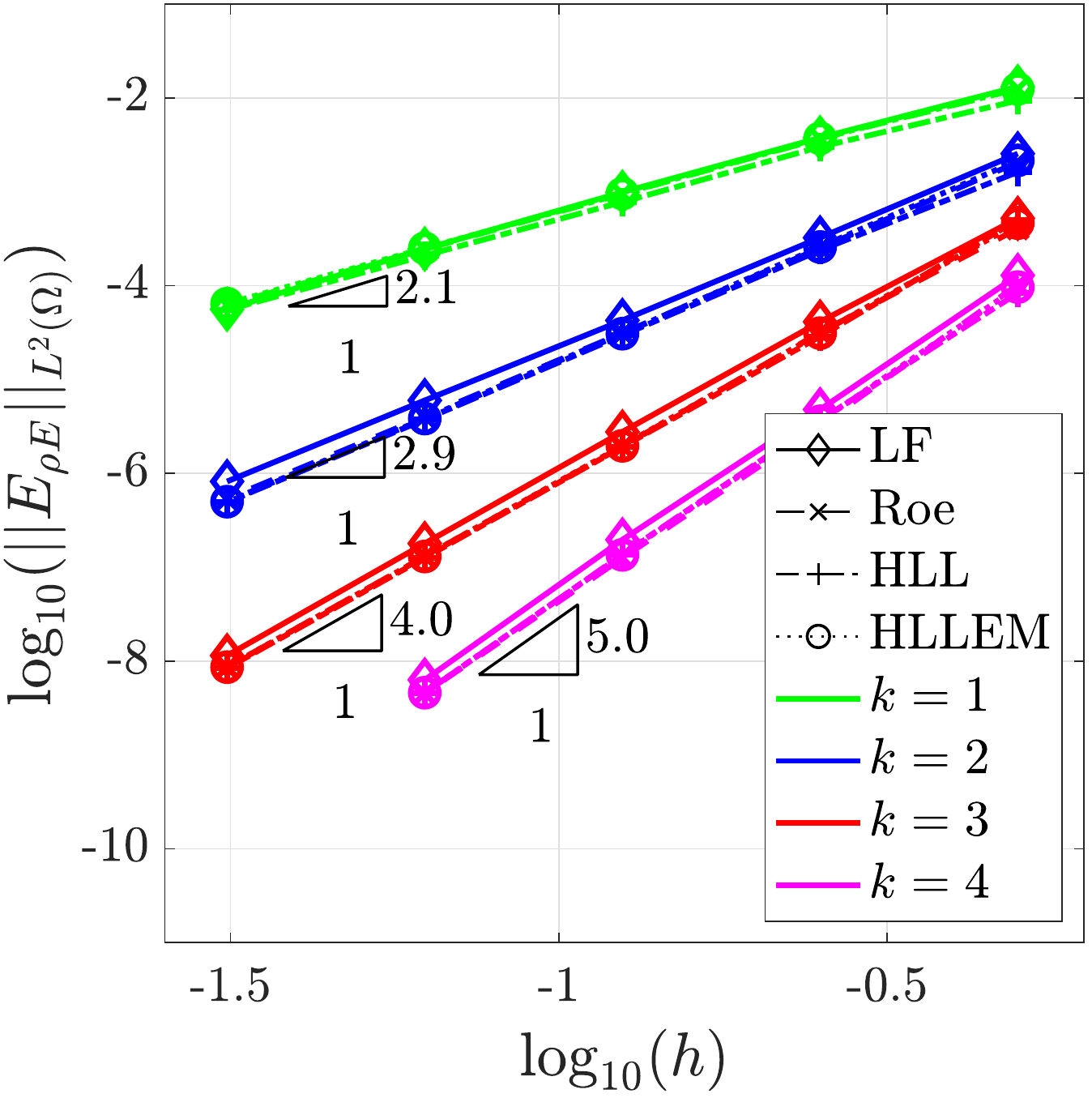}}
	\caption{Ringleb flow - Mesh convergence of the $\eltwo$ error of (a) density, (b) momentum and (c) energy, using Lax-Friedrichs (LF), Roe, HLL and HLLEM Riemann solvers and polynomial degree of approximation $k=1,\dotsc,4$.}
	\label{fig:Ringleb_Convergence}
\end{figure}

An $h$-convergence study is performed using a degree of approximation ranging from $k=1$ up to $k=4$ and for the four Riemann solvers presented in section~\ref{sc:RiemannSolvers}. 
Figure~\ref{fig:Ringleb_Convergence} displays the error for the conserved variables, i.e. $\rho$, $\rho \bv$ and $\rho E$, measured in the $\eltwo(\Omega)$ norm, as a function of the characteristic mesh size $h$.
It can be observed that the different Riemann solvers lead to an optimal rate of convergence $h^{k+1}$ and a comparable accuracy in all cases.

It is worth mentioning that the level of accuracy obtained in mesh 5 with a linear approximation $k=1$ (49,664 DOFs) is comparable to the one achieved on the coarsest mesh with polynomial degree of approximation $k=4$ (560 DOFs).
Hence, the results show the superiority of high-order discretisations, which allow to highly reduce the size of the HDG problem for a given level of accuracy.

\subsection{Convergence analysis for viscous laminar flows: Couette flow}
\label{ssc:Couette}

A compressible Couette flow with a source term~\cite{Nguyen-NP:2012,Schutz-SWM:2012} is considered to numerically verify the accuracy and convergence properties of the HDG method for the compressible Navier-Stokes equations using the different Riemann solvers presented in section~\ref{sc:RiemannSolvers}.\\
The analytical expression of the solution, defined on the square domain $\Omega = [0, 1]^2$, is
\begin{equation}\label{eq:Couette_analytical}
\begin{aligned}
\bv &= \begin{Bmatrix} y \log (1 + y) \\ 0  \end{Bmatrix}, \qquad p =\frac{ 1}{\gamma \Minf^2} \\
T &= \frac{1}{(\gamma - 1) \Minf^2} \left[ \alpha_c + y(\beta_c - \alpha_c) + \frac{(\gamma - 1)\Minf^2 \Pra}{2} y(1-y)\right],
\end{aligned}
\end{equation}
where $\alpha_c = 0.8$ and $\beta_c = 0.85$ are positive constants. The viscosity is assumed constant and the source term, which is determined from the exact solution, is given by
\begin{equation}\label{eq:Couette_source}
\bm{S} =  \frac{- 1}{\Rey} \left\lbrace 0, \,  \frac{2 + {y}}{(1+{y})^2}, \, 0, \,  \log^2(1 + {y}) + \frac{{y}\log(1+{y})}{1+{y}} + \frac{{y}(3+2{y})\log(1+{y})-2{y}-1}{(1+{y})^2} \right\rbrace  \tras.
\end{equation}

The exact solution is utilised to impose the boundary conditions on $\partial \Omega$ and the nondimensional quantities are set to $\Minf = 0.15$ and $\Rey = 1$ in order to replicate the case presented in~\cite{Nguyen-NP:2012,Schutz-SWM:2012}, taking a characteristic length $L=1$.

\begin{figure}[htbp]
	\subfloat[$k = 1$ \label{fig:Couette_Density11}]{\includegraphics[width=0.32\textwidth]{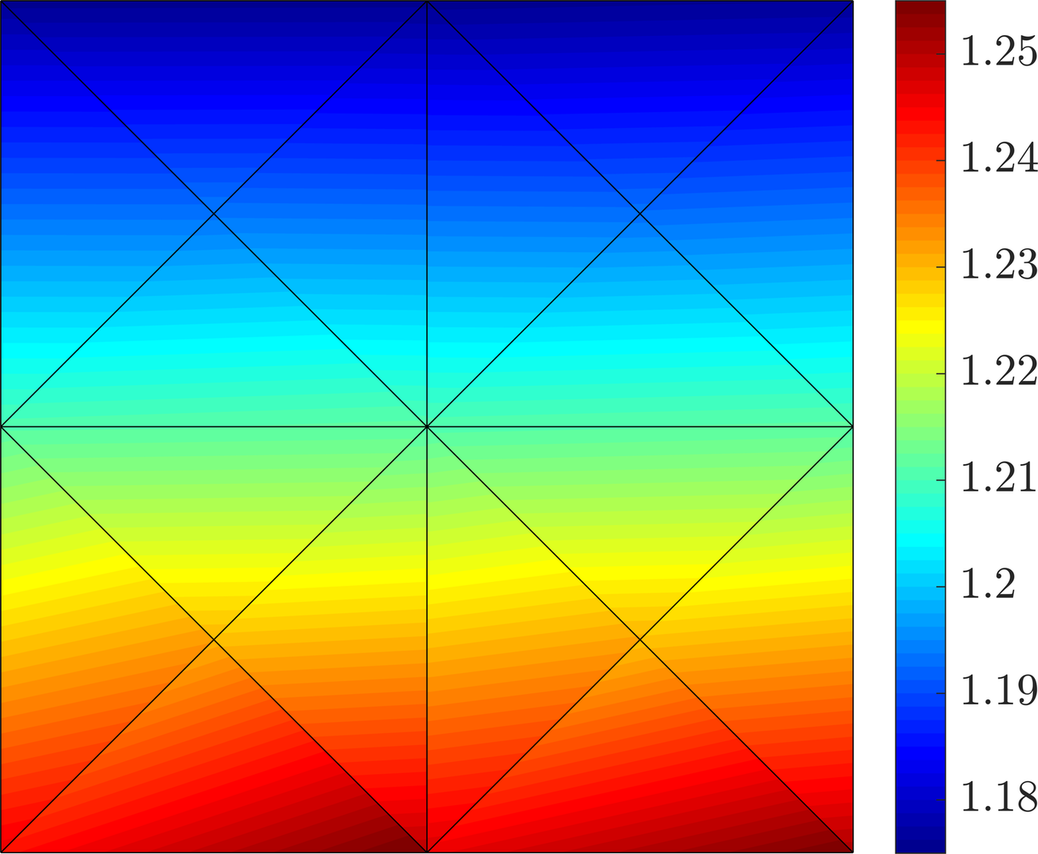}} \hfill
	\subfloat[$k = 2$ \label{fig:Couette_Density2}]{\includegraphics[width=0.32\textwidth]{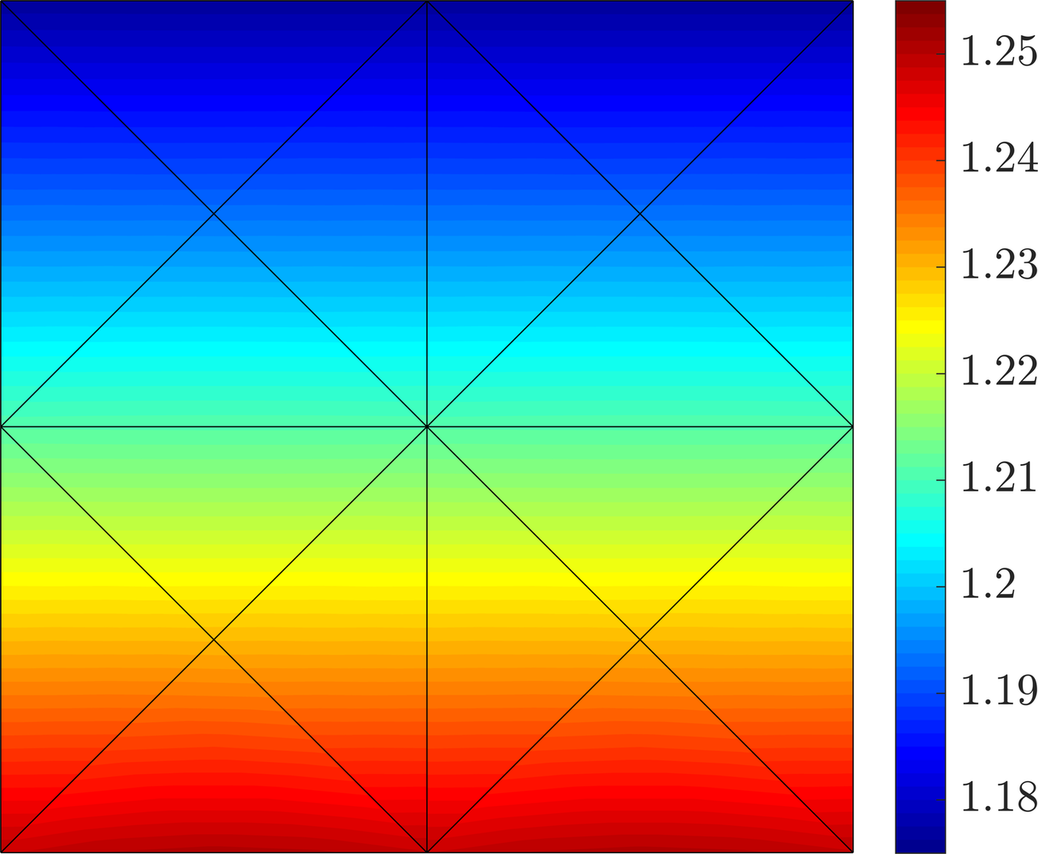}} \hfill
	\subfloat[$k = 3$ \label{fig:Couette_Density3}]{\includegraphics[width=0.32\textwidth]{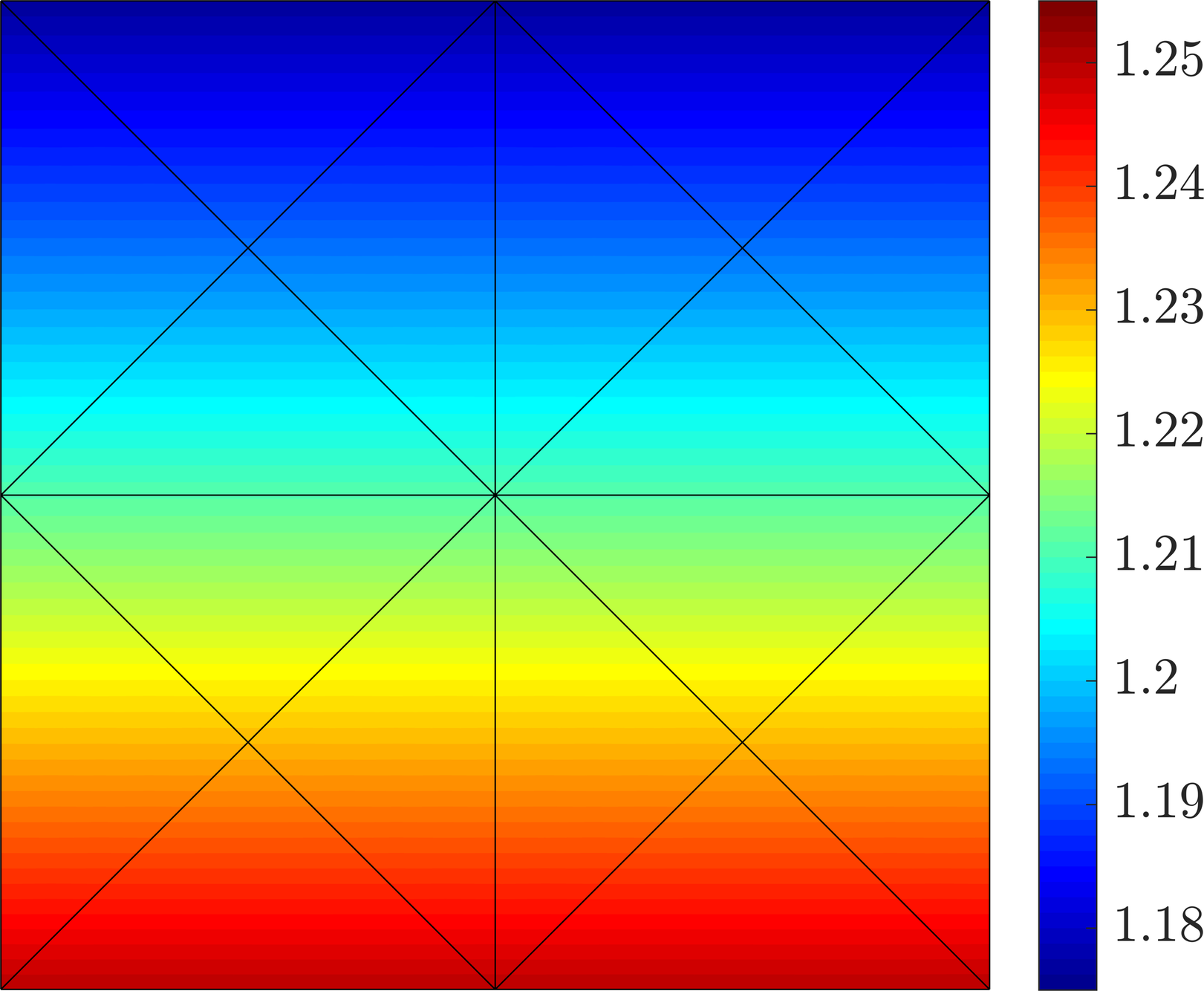}}
	\caption{Couette flow - Density distribution computed using the HLLEM Riemann solver on the first level of mesh refinement with polynomial degree $k = 1, \dotsc , 3$.}
	\label{fig:Couette_Density}
\end{figure}

The computational domain is discretised using the uniform meshes of triangular elements employed in the Ringleb example of section~\ref{ssc:Ringleb}. Figure~\ref{fig:Couette_Density} displays the approximate solution of the density field on the first mesh refinement for polynomial degrees $k = 1, \dotsc,3$.

\begin{figure}[htbp]
	\centering
	\subfloat[Density, $\rho$ \label{fig:Couette_ErrorDensity}]{\includegraphics[width=0.33\textwidth]{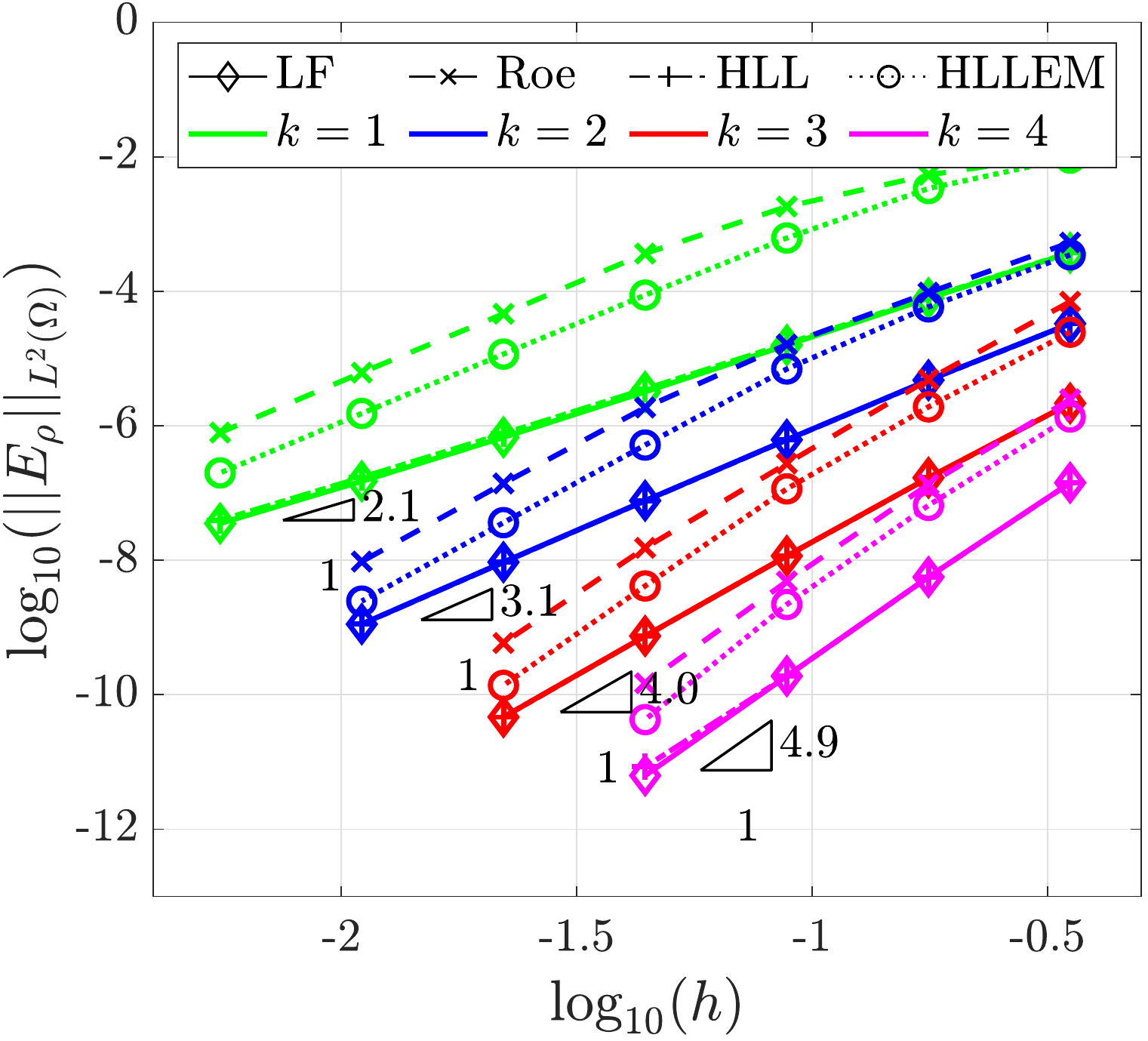}} \hfill
	\subfloat[Momentum, $\rho \bv$ \label{fig:Couette_ErrorMomentum}]{\includegraphics[width=0.33\textwidth]{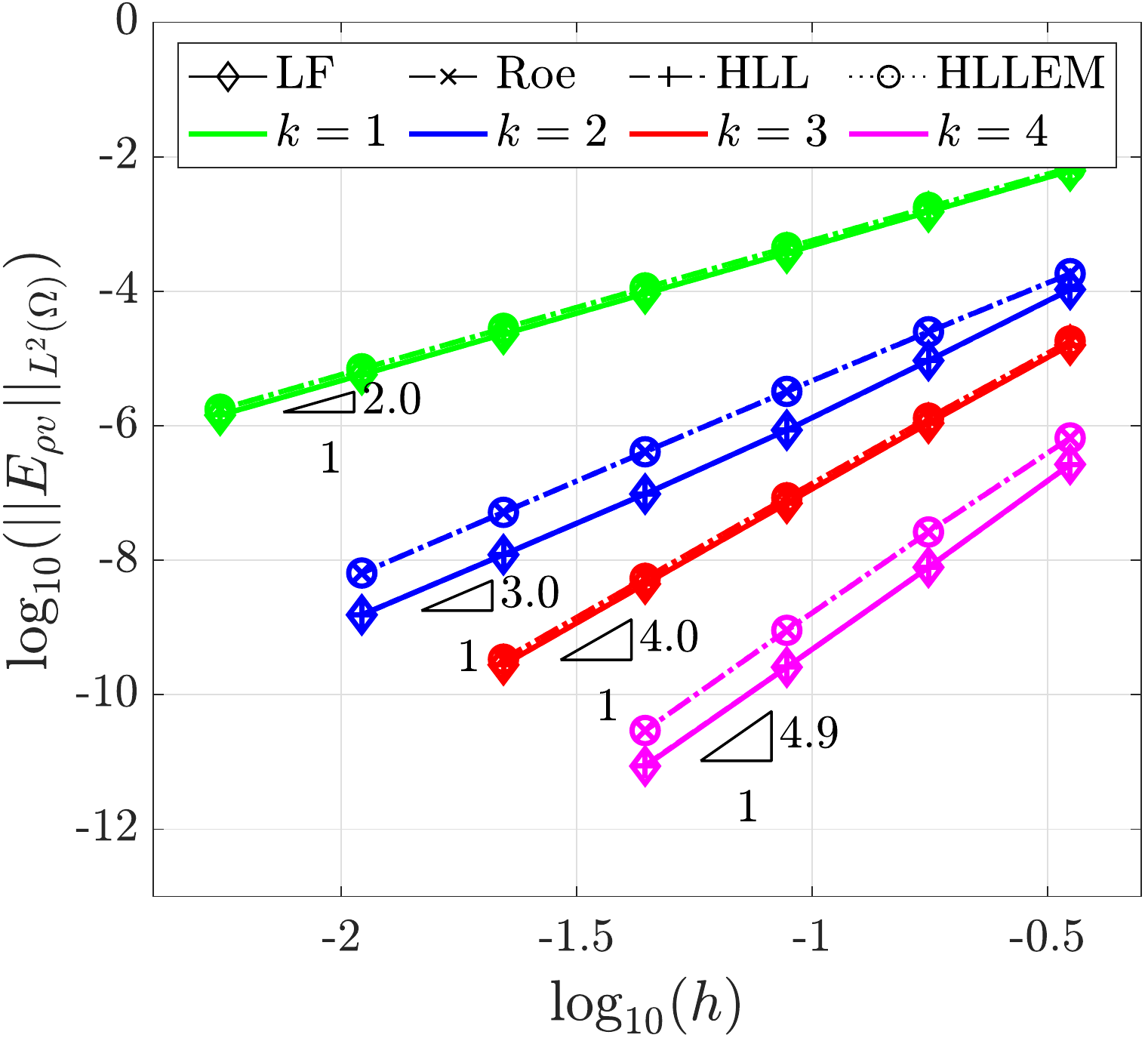}} \hfill
	\subfloat[Energy, $\rho E$ \label{fig:Couette_ErrorEnergy}]{\includegraphics[width=0.33\textwidth]{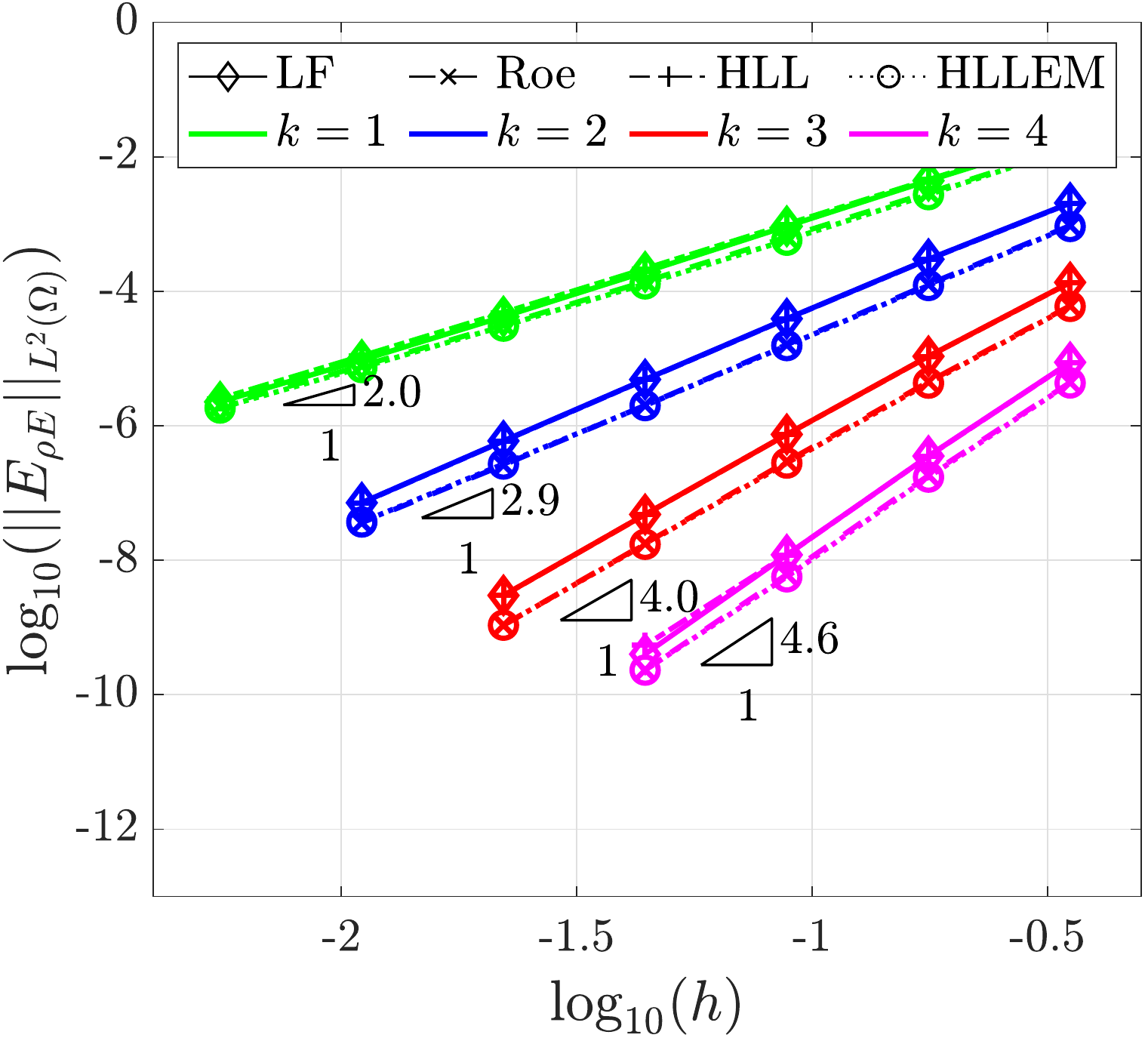}} \\
	\subfloat[Deviatoric strain rate, $\beps$ \label{fig:Couette_ErrorStrain}]{\includegraphics[width=0.33\textwidth]{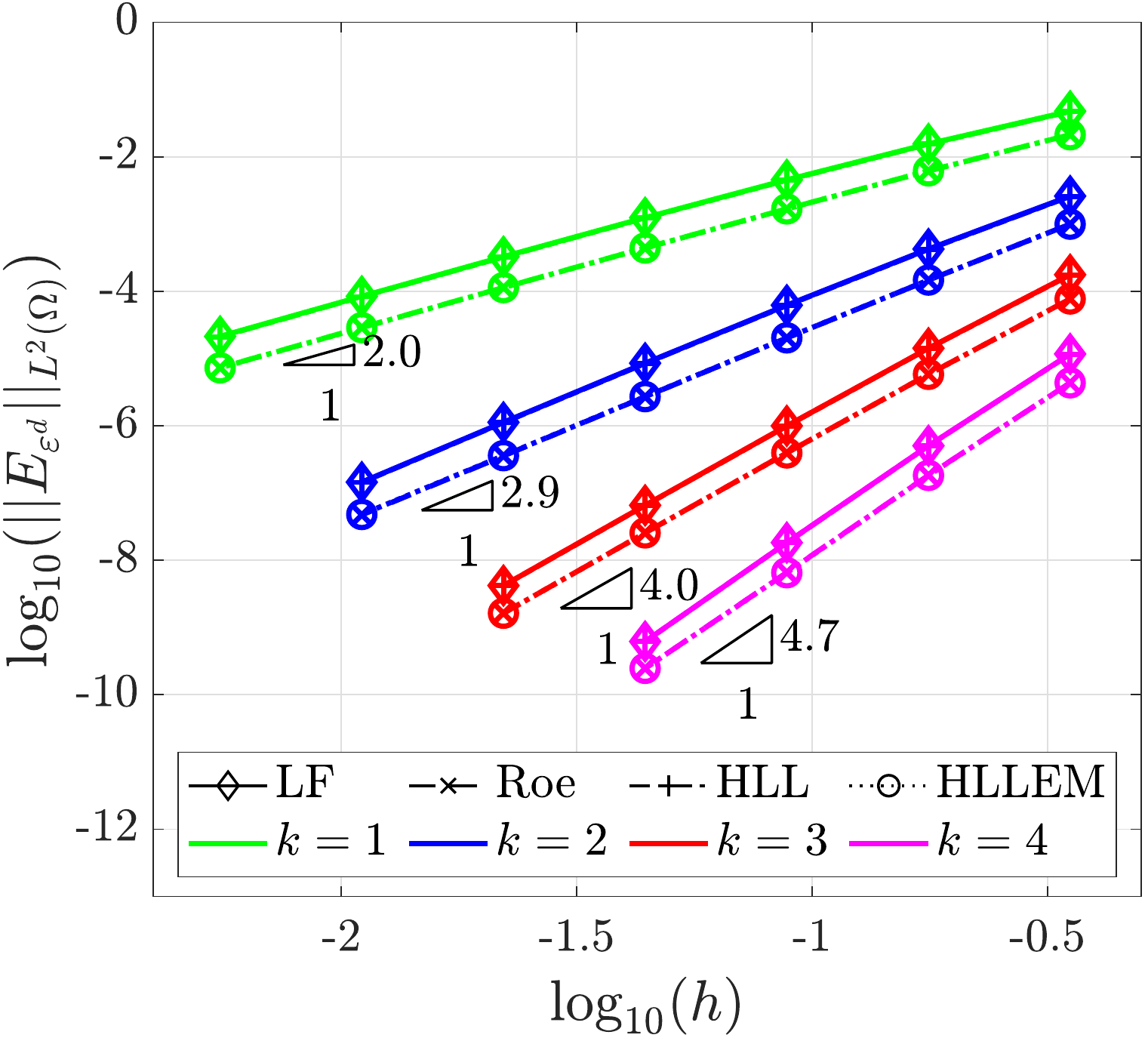}} \,
	\subfloat[Temperature gradient, $\bphi$ \label{fig:Couette_ErrorGradT}]{\includegraphics[width=0.33\textwidth]{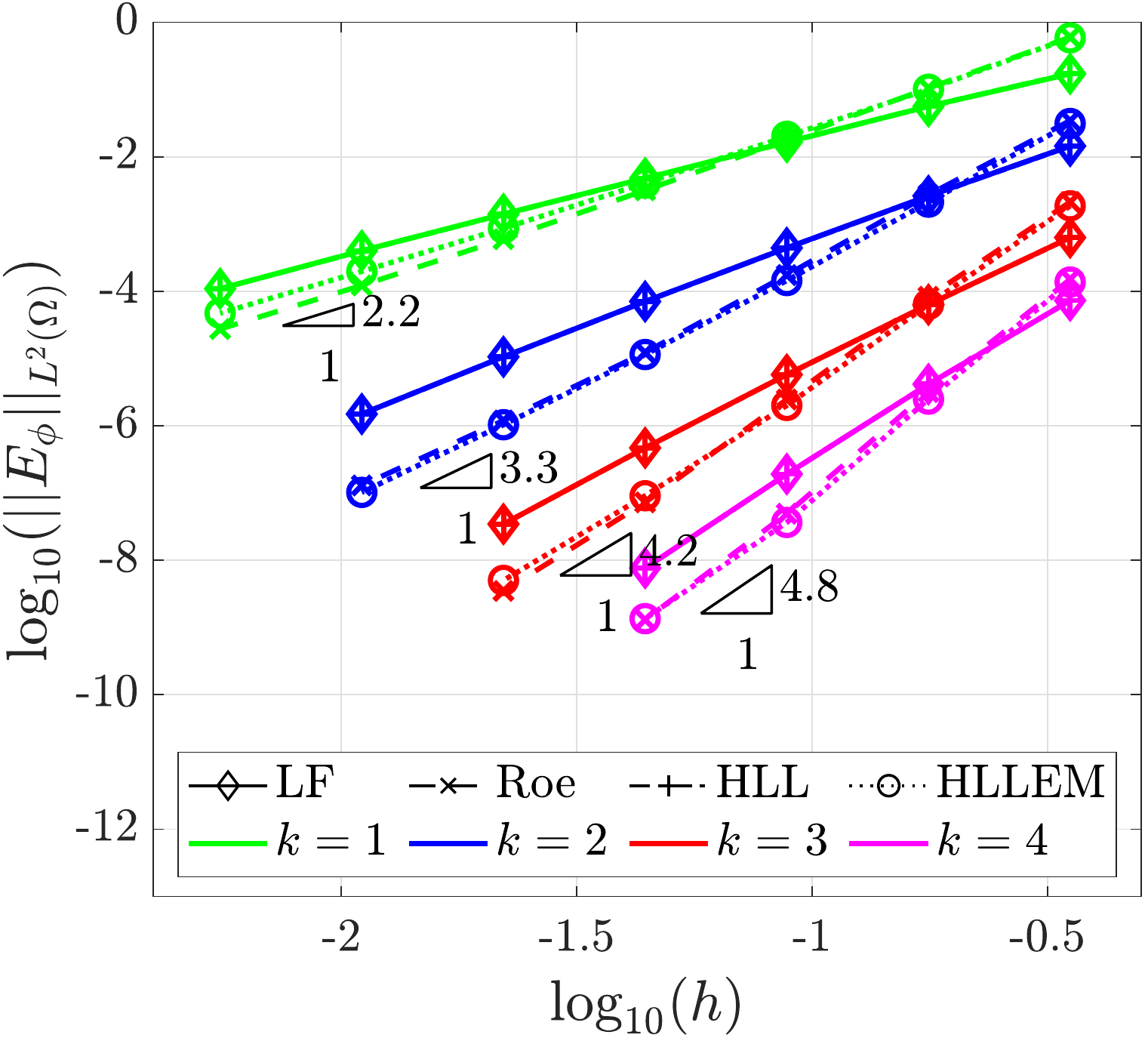}} \\
	\caption{Couette flow - Mesh convergence of the $\eltwo$ error of the (a--c) primal and (d--e) mixed variables of the discretisation, using Lax-Friedrichs (LF), Roe, HLL and HLLEM Riemann solvers and polynomial degree of approximation $k=1, \dotsc , 4$.}
	\label{fig:Couette_Error}
\end{figure}

The evolution of the error of the primal (conserved) and mixed variables measured in the $\eltwo(\Omega)$ norm is displayed in figure~\ref{fig:Couette_Error}, as a function of the characteristic element size $h$.
The $h$-convergence study compares the results of the Lax-Friedrichs, Roe, HLL and HLLEM Riemann solvers, using polynomial degrees of approximation from $k=1$ to $k=4$.
Optimal rates of convergence and comparable levels of accuracy are obtained for the approximation of the primal and mixed variables using the different Riemann solvers.

\begin{figure}[htbp]
	\centering
	\subfloat[Deviatoric strain rate, $\beps$ \label{fig:Couette_ConvergenceReStrain}]{\includegraphics[width=0.48\textwidth]{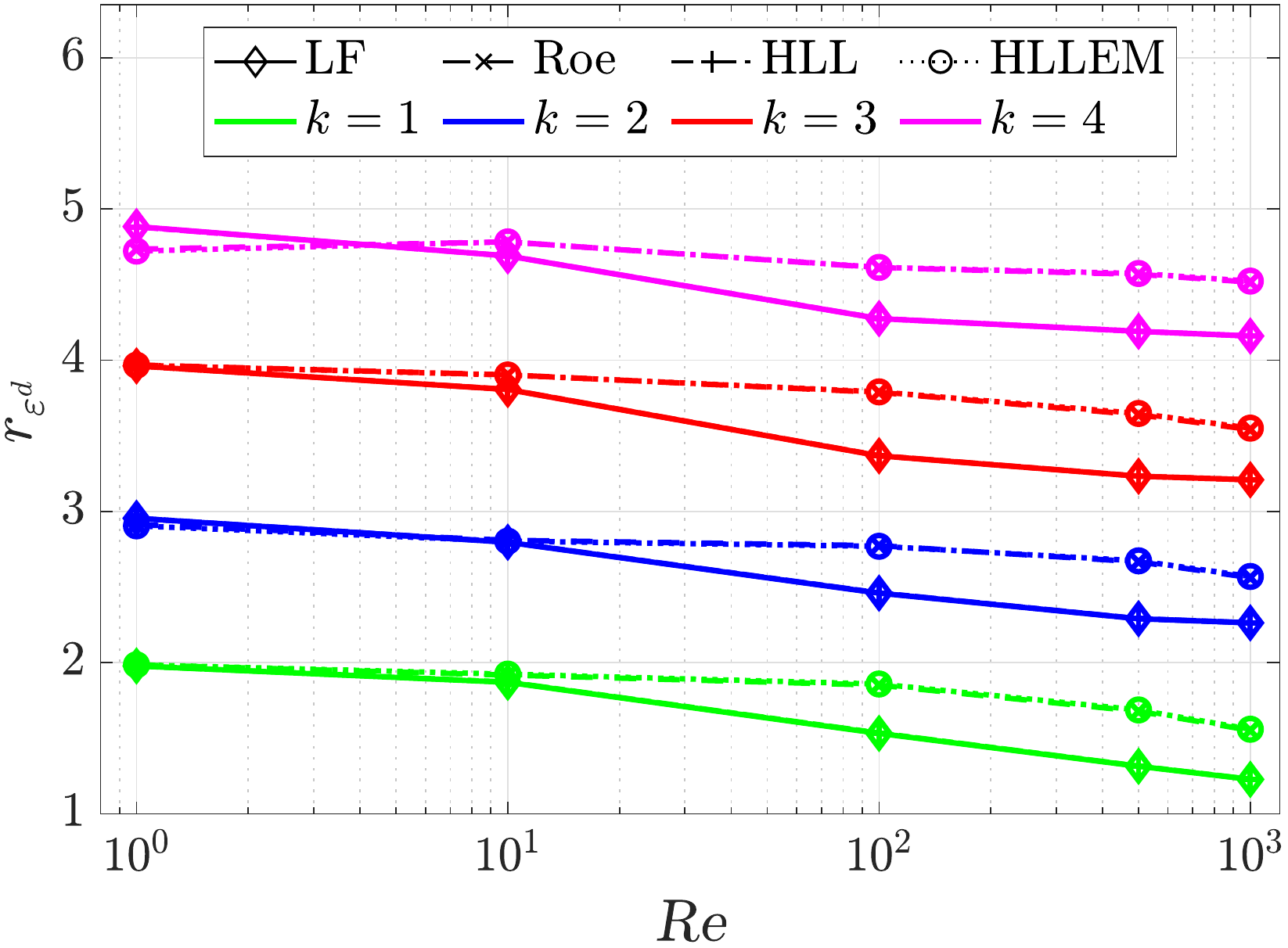}} \,
	\subfloat[Temperature gradient, $\bphi$ \label{fig:Couette_ConvergenceReGradT}]{\includegraphics[width=0.48\textwidth]{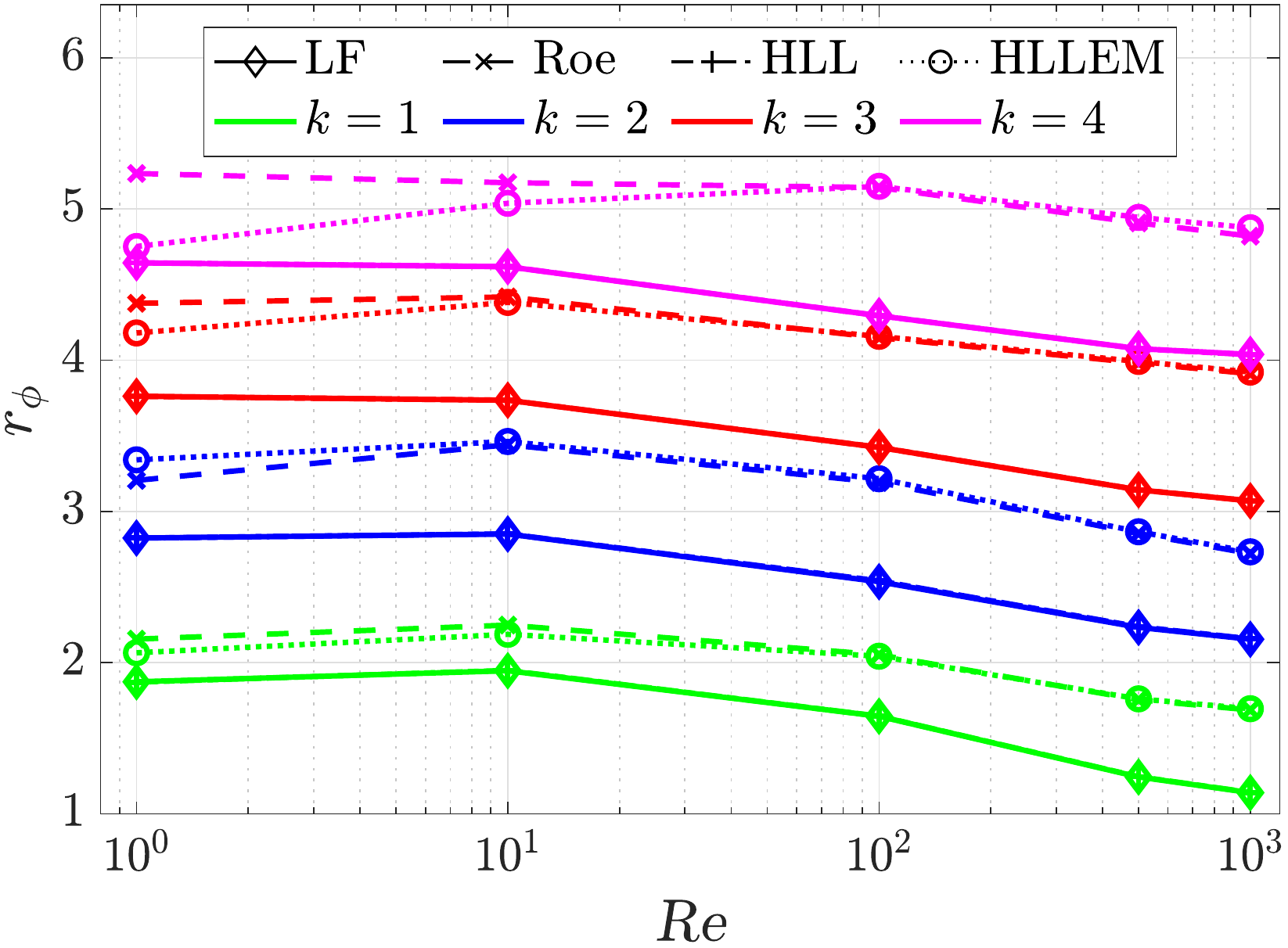}} \\
	\caption{Couette flow - Rate of convergence of the mixed variables for variable Reynolds, using Lax-Friedrichs (LF), Roe, HLL and HLLEM Riemann solvers and polynomial degree of approximation $k=1, \dotsc , 4$.}
	\label{fig:Couette_RateConvergence}
\end{figure}

Finally, the rates of convergence of the mixed variables in the last mesh refinement, $r_{\diamond}$, are examined with respect to the Reynolds number.
In particular, whereas for $\Rey = 1$ the four Riemann solvers show similar rates of convergence of $k+1$, as displayed in figure~\ref{fig:Couette_Error}, figure~\ref{fig:Couette_RateConvergence} illustrates the decreasing tendency of such convergence rates as the problem turns convection-dominated.
HLLEM and Roe Riemann solvers display an increased accuracy with respect to Lax-Friedrichs and HLL, keeping optimal rates of convergence even for $\Rey = 1000$.
On the contrary, Lax-Friedrichs and HLL exhibit a steeper drop in accuracy, experiencing a suboptimal behaviour as the Reynolds number increases.

\section{Numerical benchmarks}
\label{sc:Benchmarks}

A set of numerical examples is presented in this section to evaluate the performance and accuracy of the different Riemann solvers for inviscid and viscous compressible flows in the context of the high-order HDG method. Different cases, listed in table~\ref{tb:Examples}, are considered, ranging from viscous laminar to inviscid flows, both in subsonic, transonic and supersonic regimes.

\begin{table} [htbp]
	\caption{List of examples.}
	\centering
	\makebox[\linewidth]{
		\begin{tabular}{ L{2cm}  L{10cm}}
			\toprule
			\multicolumn{2}{l}{Inviscid examples} \\
			\midrule
			\ref{ssc:InviscidCylinder} & Subsonic flow past a circular cylinder\\
			\ref{ssc:transonicNACA} & Transonic flow over a NACA 0012 aerofoil\\
			\ref{ssc:supersonicNACA} & Supersonic flow over a NACA 0012 aerofoil\\
			\midrule
			\multicolumn{2}{l}{Viscous examples} \\ 
			\midrule
			\ref{ssc:flatplate} & Subsonic laminar flow over a flat plate\\
			\ref{ssc:SWBLI} & Shock wave/boundary layer interaction\\
			\ref{ssc:CompressionCorner} & Supersonic flow over a compression corner\\
			\bottomrule
		\end{tabular}
	}
	\label{tb:Examples}
\end{table}

\subsection{Entropy production due to geometrical error: subsonic flow past a circular cylinder}
\label{ssc:InviscidCylinder}

The subsonic flow around a circular cylinder at free-stream Mach number $M_\infty = 0.3$ is considered to assess the numerical dissipation introduced by the different Riemann solvers in the context of HDG methods. 

In particular, it is known that the geometrical error introduced by low-order descriptions of curved boundaries is responsible for a substantial nonphysical entropy production~\cite{Bassi-BR:97}. Possible solutions involve the modification of the wall boundary condition~\cite{Krivodonova-KB:2006} or the incorporation of the exact boundary representation~\cite{SevillaEuler-SFH:08}. As mentioned earlier, isoparametric approximations are considered in this work. Therefore, only approximations of degree at least $k = 2$ are reported, preventing the geometrical error from dominating over the dissipative behaviour of the Riemann solvers under analysis.

\begin{figure}[htbp]
	\centering
	\subfloat[Mesh 1 \label{fig:Cylinder_mesh1}]{\includegraphics[width=0.31\textwidth]{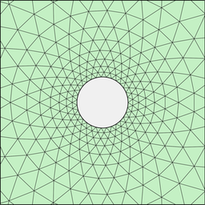}} \quad
	\subfloat[Mesh 2 \label{fig:Cylinder_mesh2}]{\includegraphics[width=0.31\textwidth]{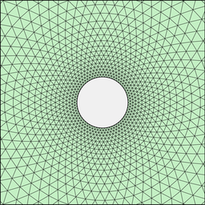}}
	\caption{Subsonic flow around a cylinder - Detail of the meshes near the 2D cylinder, featuring (a) 32 and (b) 64 subdivisions on the circular boundary.}
	\label{fig:Cylinder_meshes}
\end{figure}

Two meshes are considered for this example. The coarsest mesh consists of $1,104$ triangles with $32$ elements to discretise the circle, whereas the finest mesh has $4,635$ elements and $64$ subdivisions on the circle. A detailed view of the corresponding meshes near the cylinder is depicted in figure~\ref{fig:Cylinder_meshes}.
The far-field boundary is placed at 15 diameters from the circle and inviscid wall conditions are set on the cylinder boundary.

For isentropic subsonic flows, entropy production is a measure of the numerical dissipation introduced by the spatial discretisation. The nonphysical entropy production is computed via the so-called \emph{entropy error}, namely
\begin{equation} \label{eq:entropy}
\varepsilon_{\text{ent}} = \frac{p}{p_{\infty}} \left( \frac{\rho_{\infty}}{\rho} \right)^\gamma - 1,
\end{equation}
measuring the relative error of the total pressure with respect to the undisturbed flow in an isentropic process.

\begin{figure}[htbp]
	\subfloat[$k=2$, Mach \label{fig:Cylinder_Mach_isolines_HLL_P2}]{\includegraphics[width=0.32\textwidth]{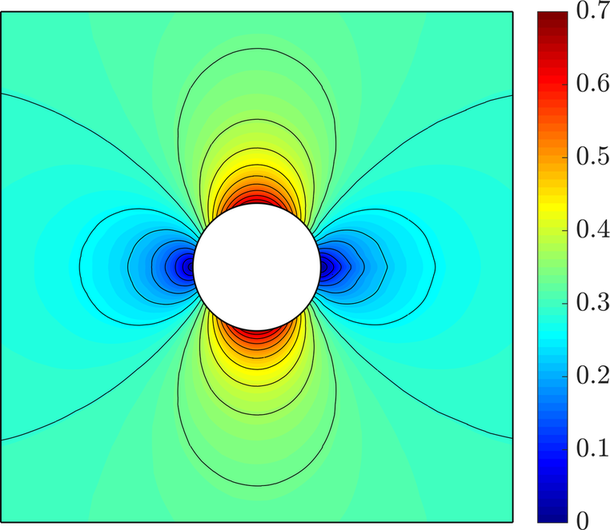}} \hfill
	\subfloat[$k=3$, Mach \label{fig:Cylinder_Mach_isolines_HLL_P3}]{\includegraphics[width=0.32\textwidth]{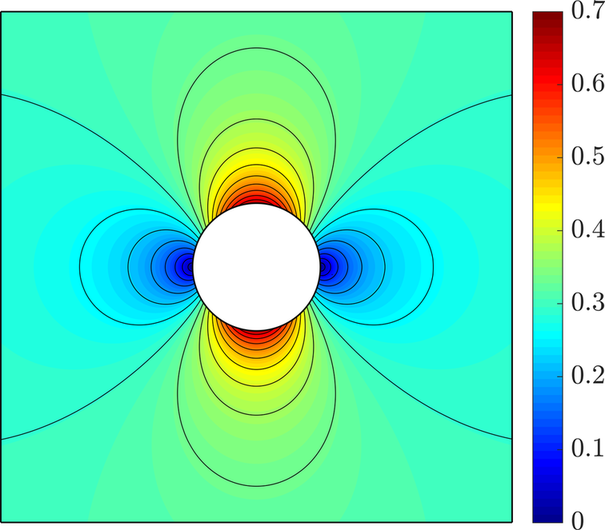}} \hfill
	\subfloat[$k=4$, Mach \label{fig:Cylinder_Mach_isolines_HLL_P4}]{\includegraphics[width=0.32\textwidth]{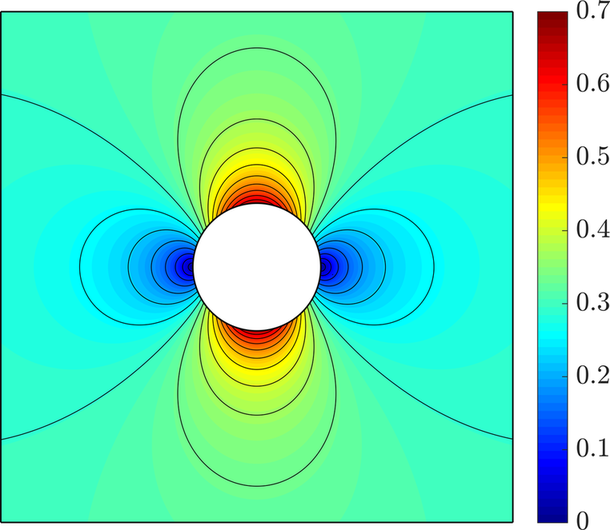}} \\
	\subfloat[$k=2$, entropy error \label{fig:Cylinder_HLL_entropy_P2}]{\includegraphics[width=0.32\textwidth]{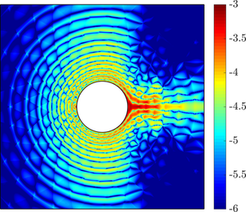}}  \hfill
	\subfloat[$k=3$, entropy error \label{fig:Cylinder_HLL_entropy_P3}]{\includegraphics[width=0.32\textwidth]{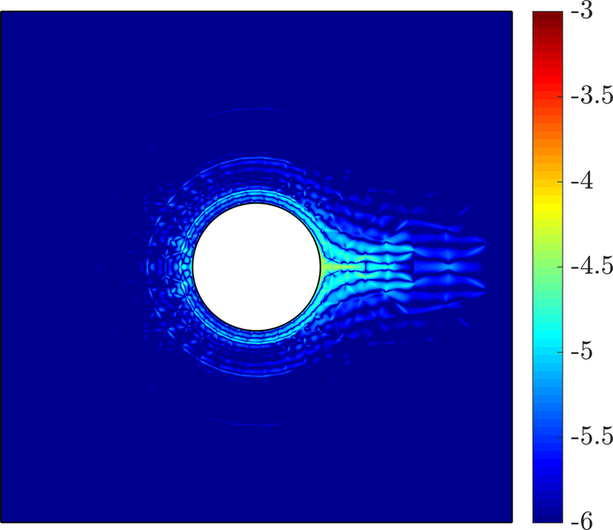}}  \hfill
	\subfloat[$k=4$, entropy error \label{fig:Cylinder_HLL_entropy_P4}]{\includegraphics[width=0.32\textwidth]{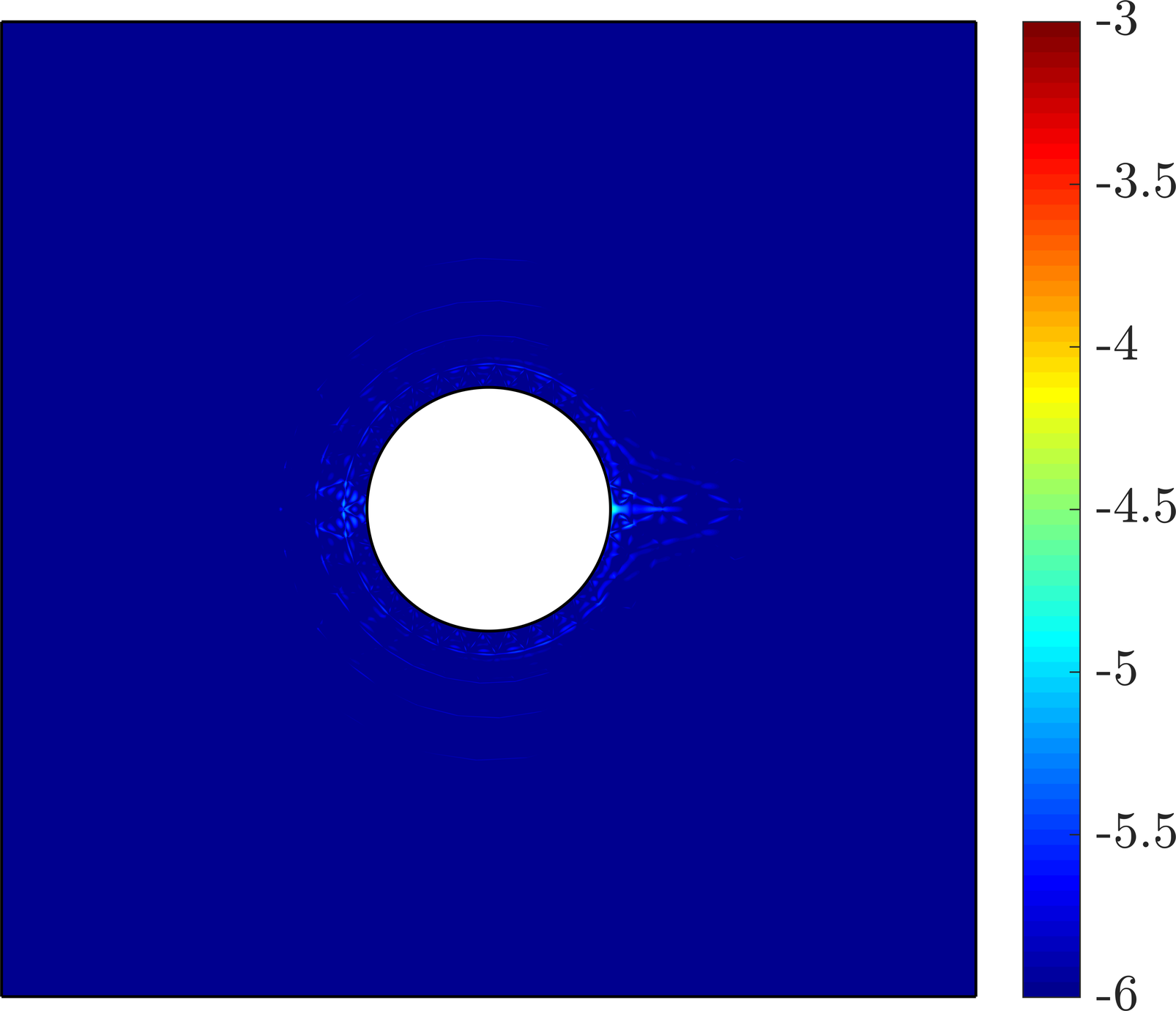}} 
	\caption{Subsonic flow around a cylinder - Mach number distribution and isolines (top) and entropy error in logarithmic scale (bottom) computed on the first mesh using the HLL Riemann solver with $k=2$ (left), $k=3$ (middle) and $k=4$ (right).}
	\label{fig:Cylinder_MachEntropy}
\end{figure}

Figure~\ref{fig:Cylinder_MachEntropy} (top) shows the Mach number distribution and isolines of the numerical solution computed on the first mesh with $k = 2,\ldots,4$, using the HLL Riemann solver. Although the computed distribution of the Mach number is comparable in the three settings, the superiority of high-order approximations becomes evident when the corresponding entropy errors are compared (figure~\ref{fig:Cylinder_MachEntropy}, bottom). The results clearly display that, increasing the polynomial degree of discretisation, the numerical dissipation introduced by the method is localised in the vicinity of the cylinder and its overall amount is reduced.

\begin{SCfigure}[][htbp]
	\includegraphics*[width=0.6\textwidth]{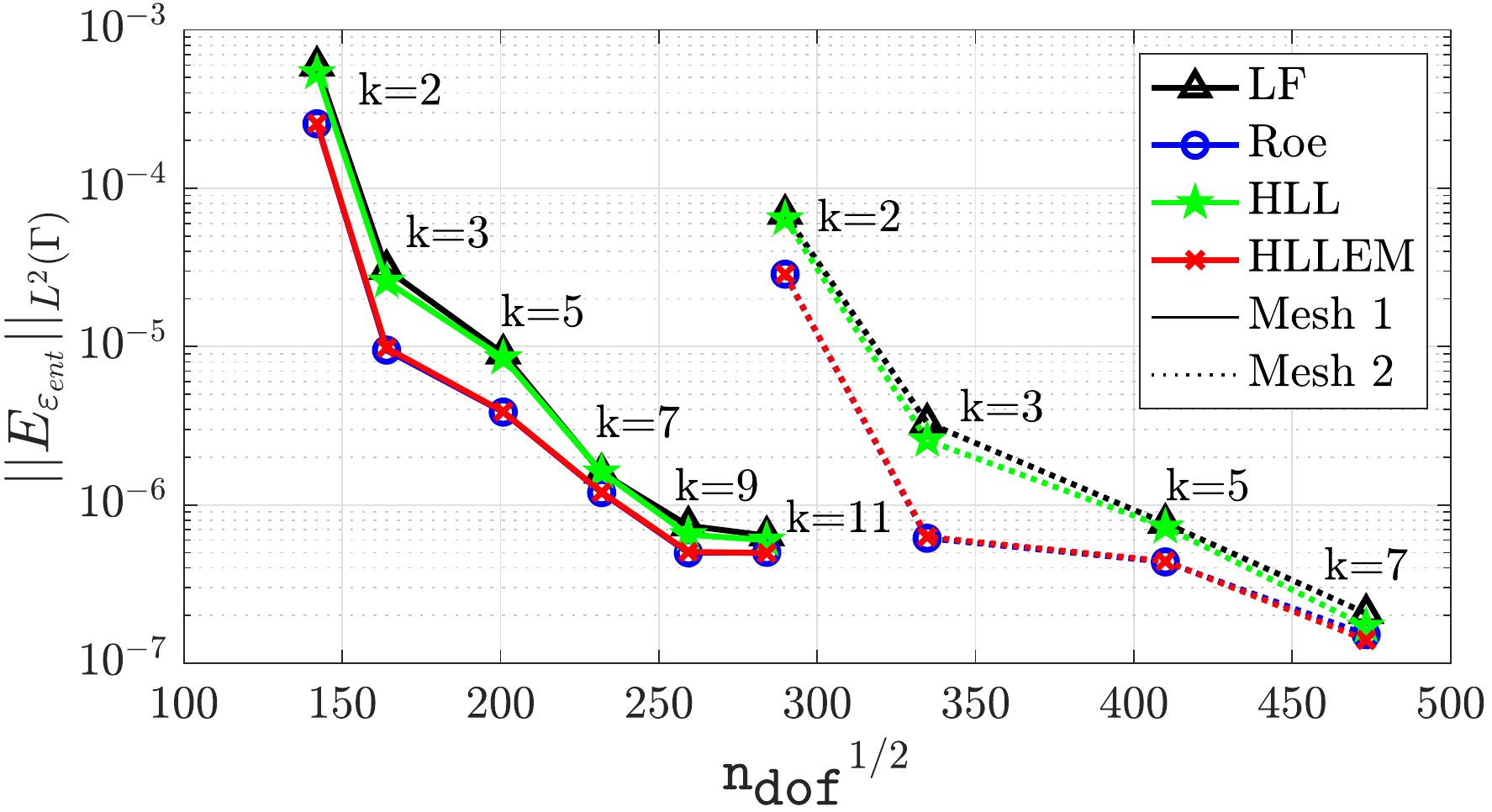}
	\caption{Subsonic flow around a cylinder - Entropy error on the cylinder surface for different meshes and different degrees of polynomial approximation.}
	\label{fig:Cylinder_entropy}
\end{SCfigure}

To quantify the differences between the four Riemann solvers, the nonphysical entropy production is compared through the $\eltwo$ norm of the entropy error, measured on the surface of the cylinder. Figure~\ref{fig:Cylinder_entropy} displays the quantity~\eqref{eq:entropy} as a function of the number of degrees of freedom of the global problem, for the two meshes under analysis and an increasing value of the polynomial degree used to approximate the solution.
The results show that the entropy production of the HLL Riemann solver is almost identical when compared to the Lax-Friedrichs Riemann solver, whereas HLLEM matches the entropy production by the Roe numerical flux. Moreover, as expected for a subsonic flow, the entropy production is slightly lower for the HLLEM and Roe Riemann solvers.

It is worth noting that the differences among the Riemann solvers are less important as the polynomial degree of the approximation increases. This confirms the observation above on the reduced amount of numerical dissipation introduced by the method as the degree of the discretisation increases and the consequent extra accuracy provided by high-order approximations.

\subsection{Boundary layer resolution: subsonic viscous laminar flow over a flat plate}
\label{ssc:flatplate}

The next example consists of the subsonic laminar flow over a flat plate at zero angle of attack.
This classical benchmark follows from the analytical study of boundary layers by Blasius for incompressible flows~\cite{Blasius1908} and has been commonly used to test laminar flow solvers in resolving boundary layers \cite{Schlichting-BL:2016}.

This problem is used to evaluate the numerical diffusion introduced by the different Riemann solvers in the approximation of shear layers and its effect over the boundary layer description.

The example considers a nearly incompressible flow ($\Minf = 0.1$) at a high Reynolds number ($\Rey = 10^5$) while preserving a laminar behaviour of the solution along the flat plate.

\begin{figure}[htbp]
	\centering
	\includegraphics*[width=0.8\textwidth]{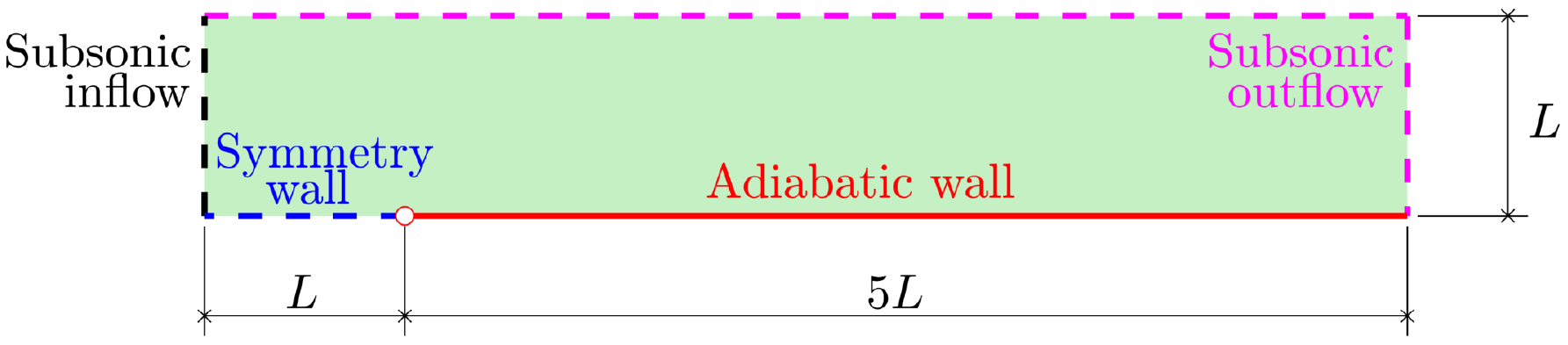}
	\caption{Laminar flow over a flat plate - Sketch of the geometry and boundary conditions.}
	\label{fig:flatplate_Domain}
\end{figure}

The computational domain consists of a flat plate of length $5L$, being $L$ the characteristic length of the problem, embedded in a rectangular domain, as shown in figure~\ref{fig:flatplate_Domain}. Adiabatic wall conditions are imposed along the plate, whereas symmetry wall conditions are imposed upstream of the leading edge. Subsonic inflow and outflow conditions are imposed at the outer boundaries.
The pressure at the outflow is set to $p_\infty$, forcing a zero pressure drop.

Uniform mesh refinement of the boundary layer is performed in order to analyse the convergence of the solution. Details of the refinement are reported in table~\ref{tb:flatplate_meshRefinement}.
In particular, for each level of refinement, the number of layers of elements in the boundary layer, $\nlayers$, and the number of subdivisions along the flat plate, $\ndiv$, are doubled and the height $h_0$ of the first layer is halved. Additionally, $h_0$ is chosen according to the relation $h_0/k \sim \Rey^{-0.75} L$.
Finally, the geometric growth rate of the boundary layer, $r$,  is determined in order for the height of the boundary layer mesh to be $h/L = 0.1$.
\begin{table} [htbp]
	\centering
	\caption{Laminar flow over a flat plate - Mesh refinement details for the convergence study.}
	\begin{tabular}{ C{2.5cm} C{2.25cm} C{2.25cm} C{2.25cm} C{2.25cm} C{2.25cm} }
		\toprule
		Refinement & $\nlayers$ & $\ndiv$ & $h_0/L$  & $r$ & $\numel$ \\
		\midrule
		1 & 4 & 10 & 8 $\cdot 10^{-4}$ & 4 & 501 \\
		2 & 8 & 20 & 4 $\cdot 10^{-4}$ & 2 & 1,154 \\
		3 & 16 & 40 & 2 $\cdot 10^{-4}$ & 1.4 & 3,512 \\
		\bottomrule
	\end{tabular}
	\label{tb:flatplate_meshRefinement}
\end{table}

\begin{figure}[htbp]
	\centering
	\subfloat[Mesh 1. \label{fig:flatplate_meshH1}]{\includegraphics[width=0.65\textwidth]{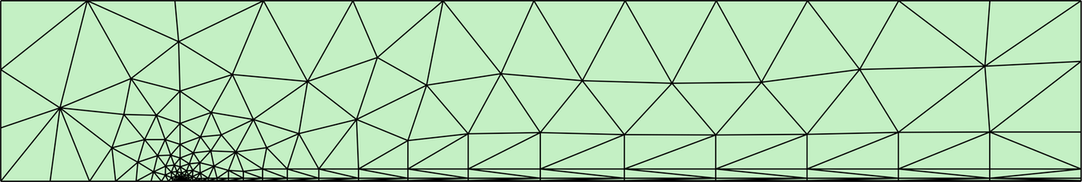}} \\
	\subfloat[Mesh 2. \label{fig:flatplate_meshH2}]{\includegraphics[width=0.65\textwidth]{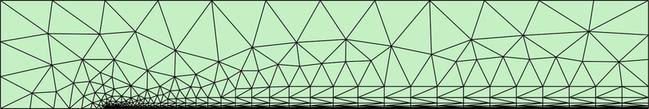}} \\ 
	\subfloat[Mesh 3. \label{fig:flatplate_meshH3}]{\includegraphics[width=0.65\textwidth]{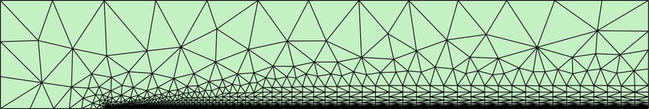}} 
	\caption{Laminar flow over a flat plate - Meshes used for the convergence study.}
	\label{fig:flatplate_meshes}
\end{figure}

The three mesh refinements used for this study are displayed in figure~\ref{fig:flatplate_meshes}. Because of the explicit embedding of the flat plate on the lower boundary of the domain, a singularity is introduced at the leading edge \cite{Tuck1991}. To alleviate its numerical effects, the mesh is further refined at this location.

\begin{figure}[htbp]
	\subfloat[Mesh 1, $k = 1$ \label{fig:flatplate_Cf_H1K1}]{\includegraphics[width=0.48\textwidth]{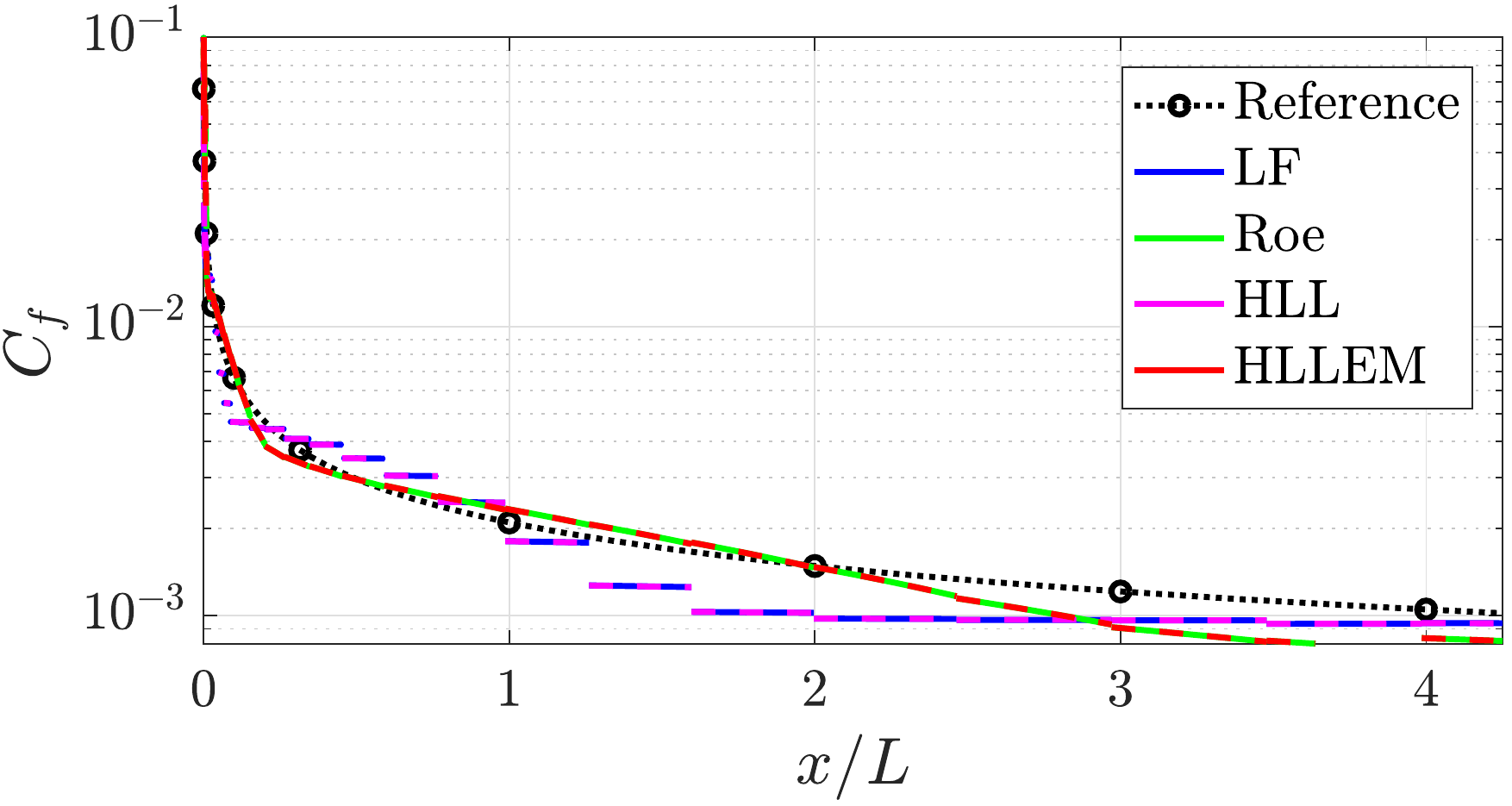}} \hfill
	\subfloat[Mesh 1, $k = 3$ \label{fig:flatplate_Cf_H1K3}]{\includegraphics[width=0.48\textwidth]{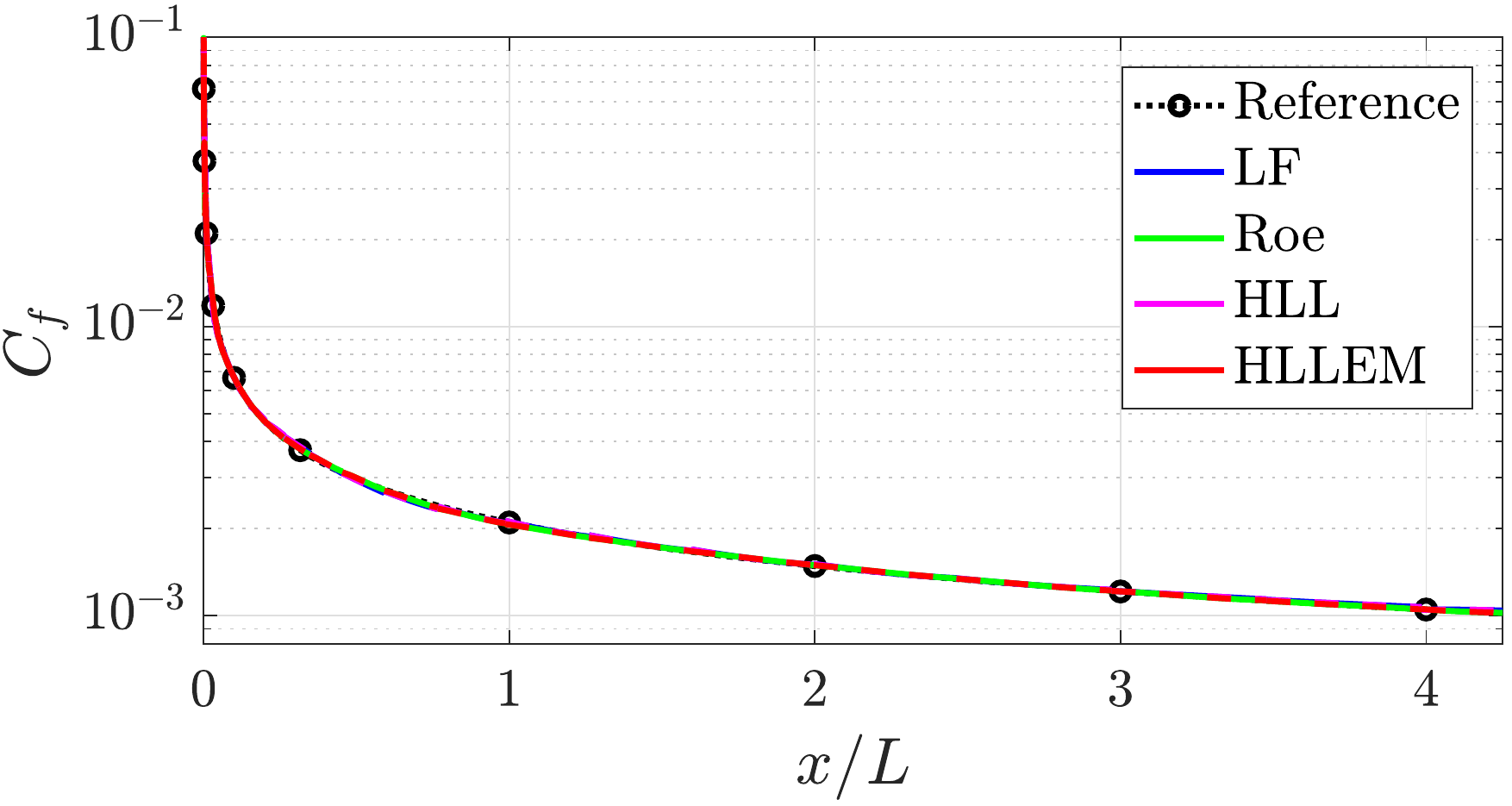}}
	\caption{Laminar flow over a flat plate - Friction coefficient along the flat plate for different polynomial degrees of approximation in the coarsest mesh, using the Lax-Friedrichs (LF), Roe, HLL and HLLEM Riemann solvers. The reference solution is obtained using the HLLEM Riemann solver on the third mesh, with $k=4$.}
	\label{fig:flatplate_Cf}
\end{figure}

The skin friction coefficient computed in the first level of refinement using degree of approximation $k=1$ and $k=3$ is depicted in figure~\ref{fig:flatplate_Cf} for the different Riemann solvers.
The superiority of Roe and HLLEM Riemann solvers with respect to classical Lax-Friedrichs and HLL is clearly displayed in figure~\ref{fig:flatplate_Cf_H1K1}: for low-order approximations, HLLEM and Roe achieve a better accuracy due to their ability to capture contact wave-type phenomena and consequently, boundary layer effects by introducing a lower amount of numerical dissipation. Of course, such difference is reduced when high-order approximations are considered, as the numerical dissipation of the method decreases, see figure~\ref{fig:flatplate_Cf_H1K3}.

\begin{figure}[htbp]
	\subfloat[Mesh 1, $k = 1$ \label{fig:flatplate_VelocityProfiles_H1K1}]{\includegraphics[width=0.32\textwidth]{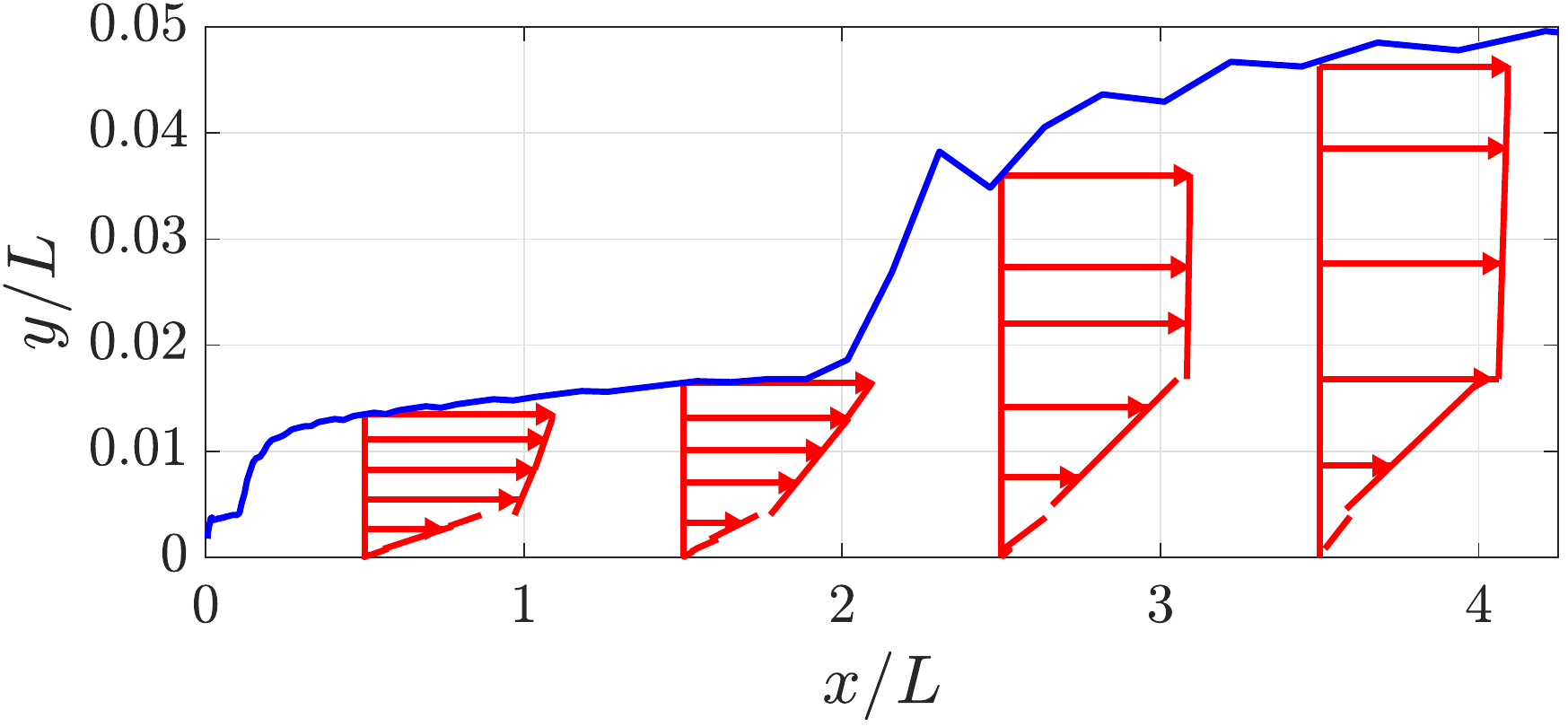}} \hfill
	\subfloat[Mesh 2, $k = 1$ \label{fig:flatplate_VelocityProfiles_H2K1}]{\includegraphics[width=0.32\textwidth]{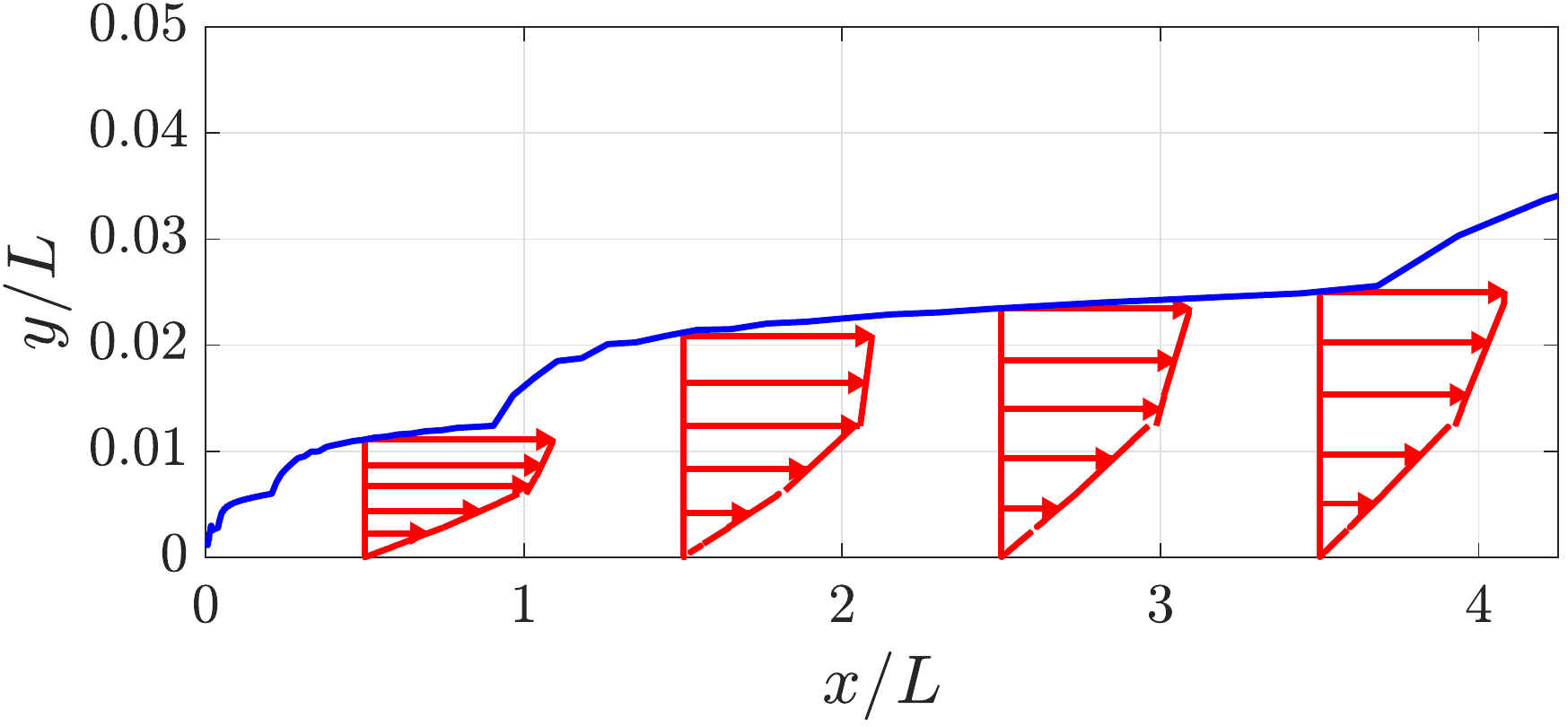}}  \hfill
	\subfloat[Mesh 3, $k = 1$ \label{fig:flatplate_VelocityProfiles_H3K1}]{\includegraphics[width=0.32\textwidth]{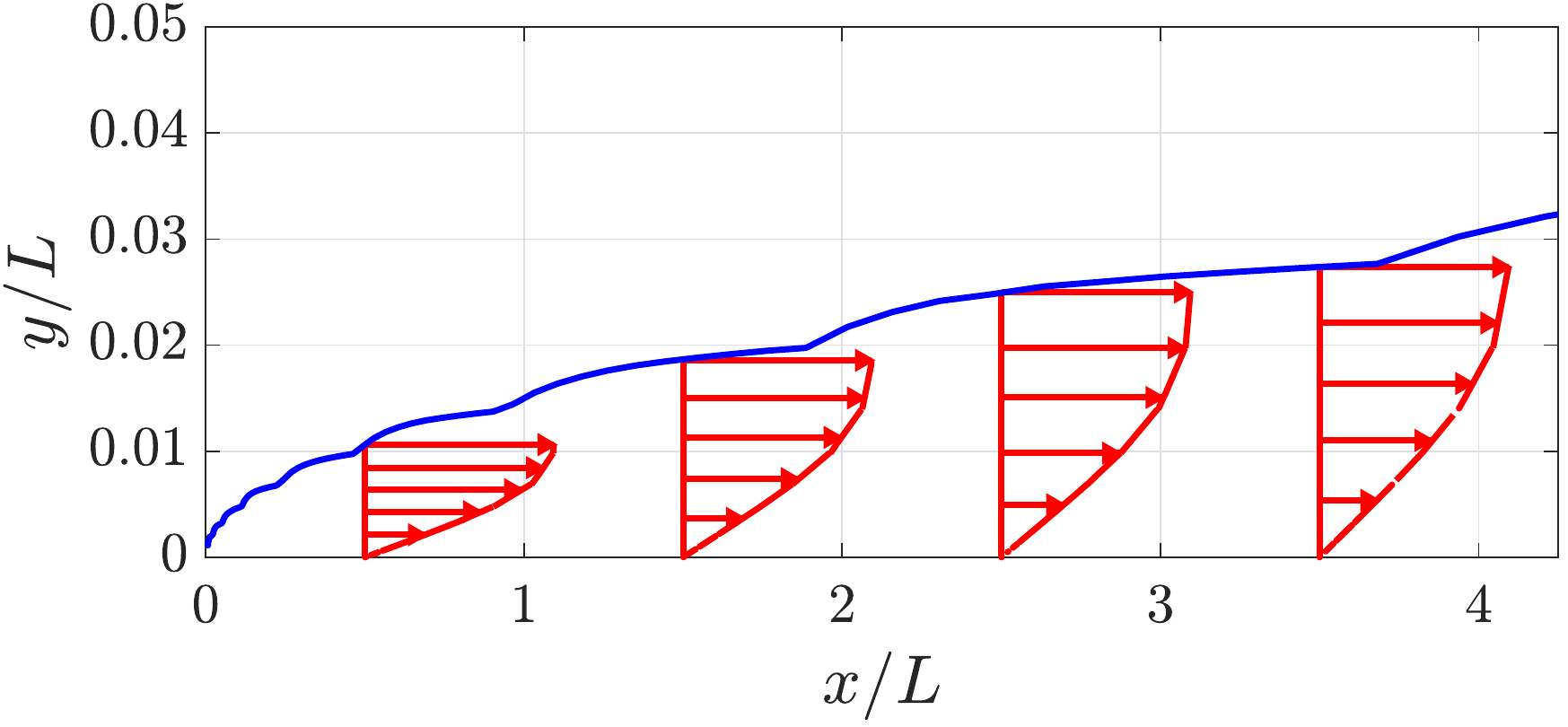}}\\
	\subfloat[Mesh 1, $k = 2$ \label{fig:flatplate_VelocityProfiles_H1K2a}]{\includegraphics[width=0.32\textwidth]{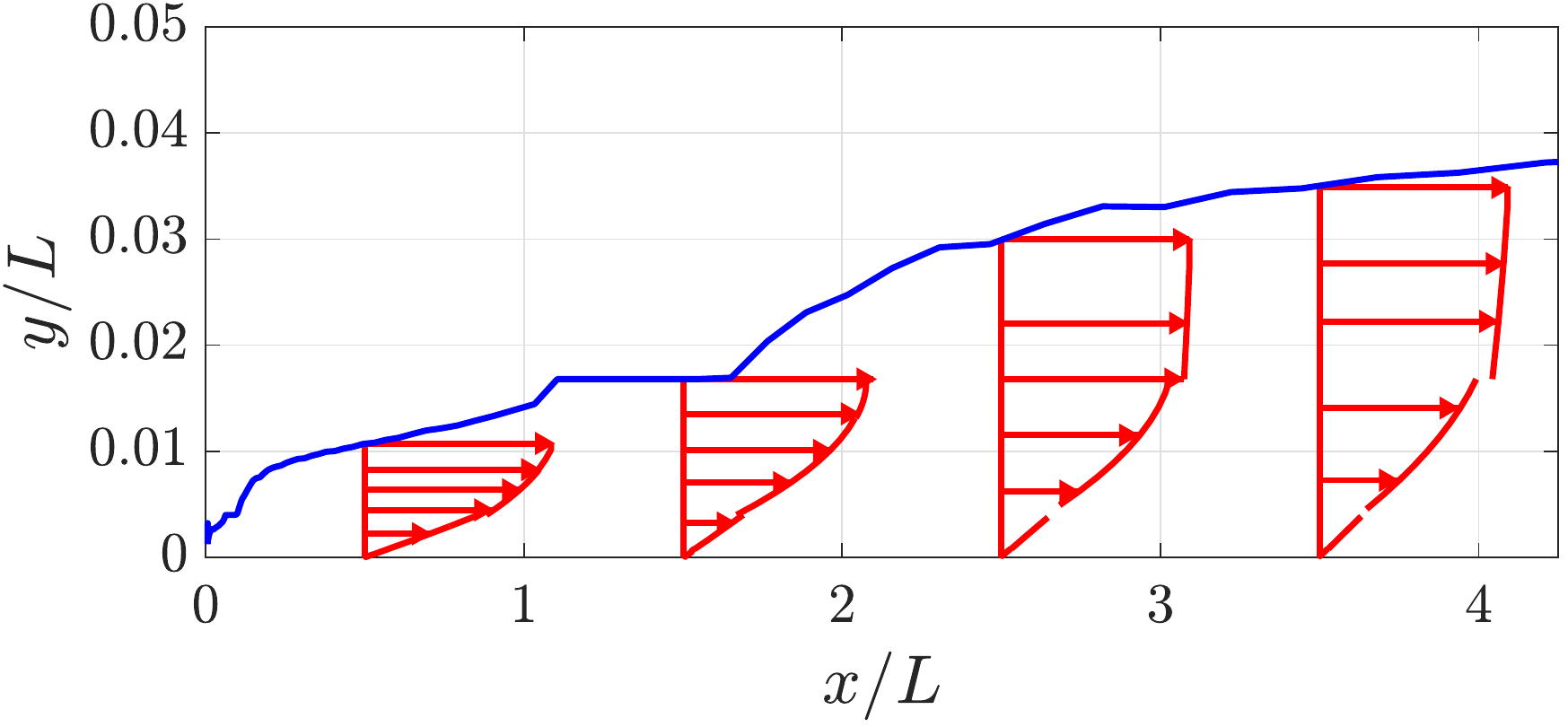}} \hfill
	\subfloat[Mesh 1, $k = 3$ \label{fig:flatplate_VelocityProfiles_H1K3a}]{\includegraphics[width=0.32\textwidth]{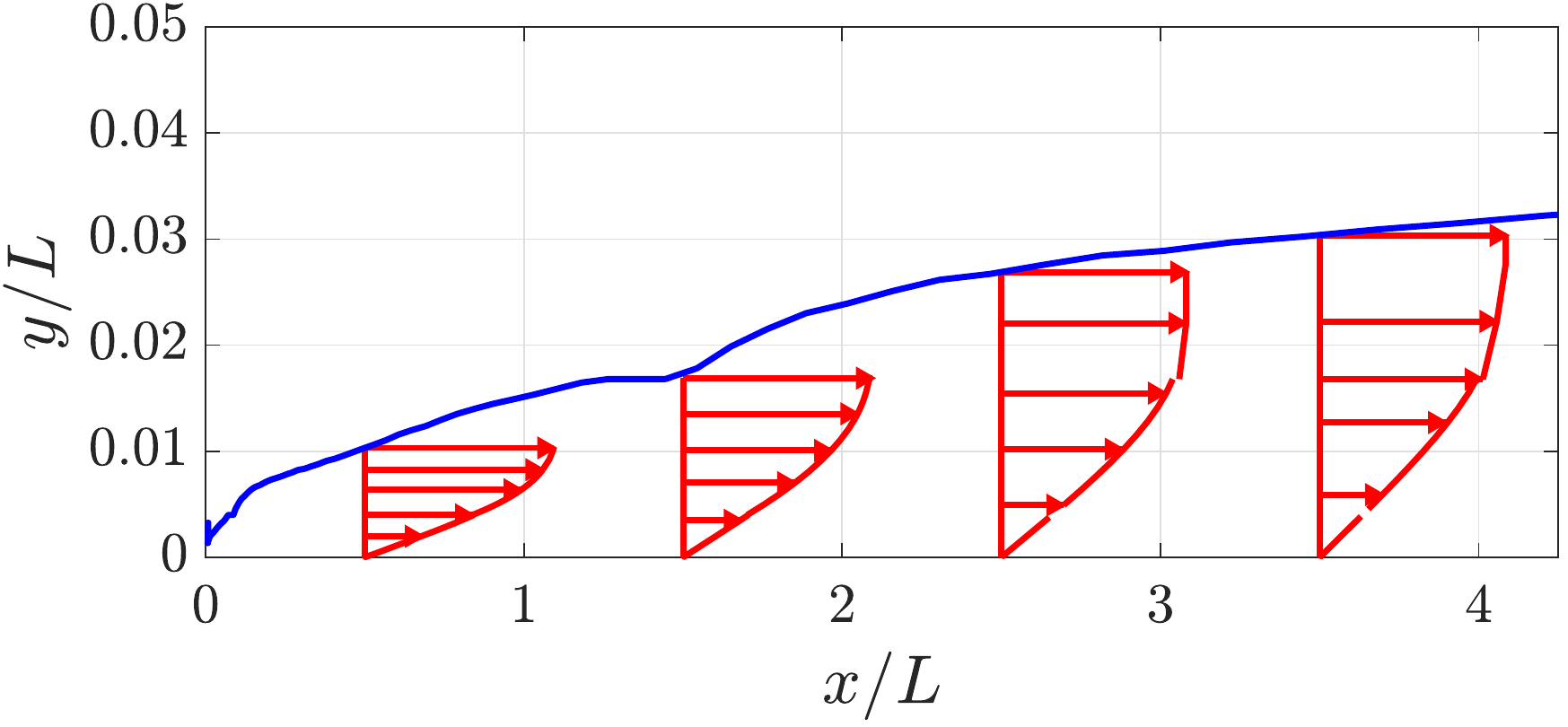}}  \hfill
	\subfloat[Mesh 1, $k = 4$ \label{fig:flatplate_VelocityProfiles_H1K4a}]{\includegraphics[width=0.32\textwidth]{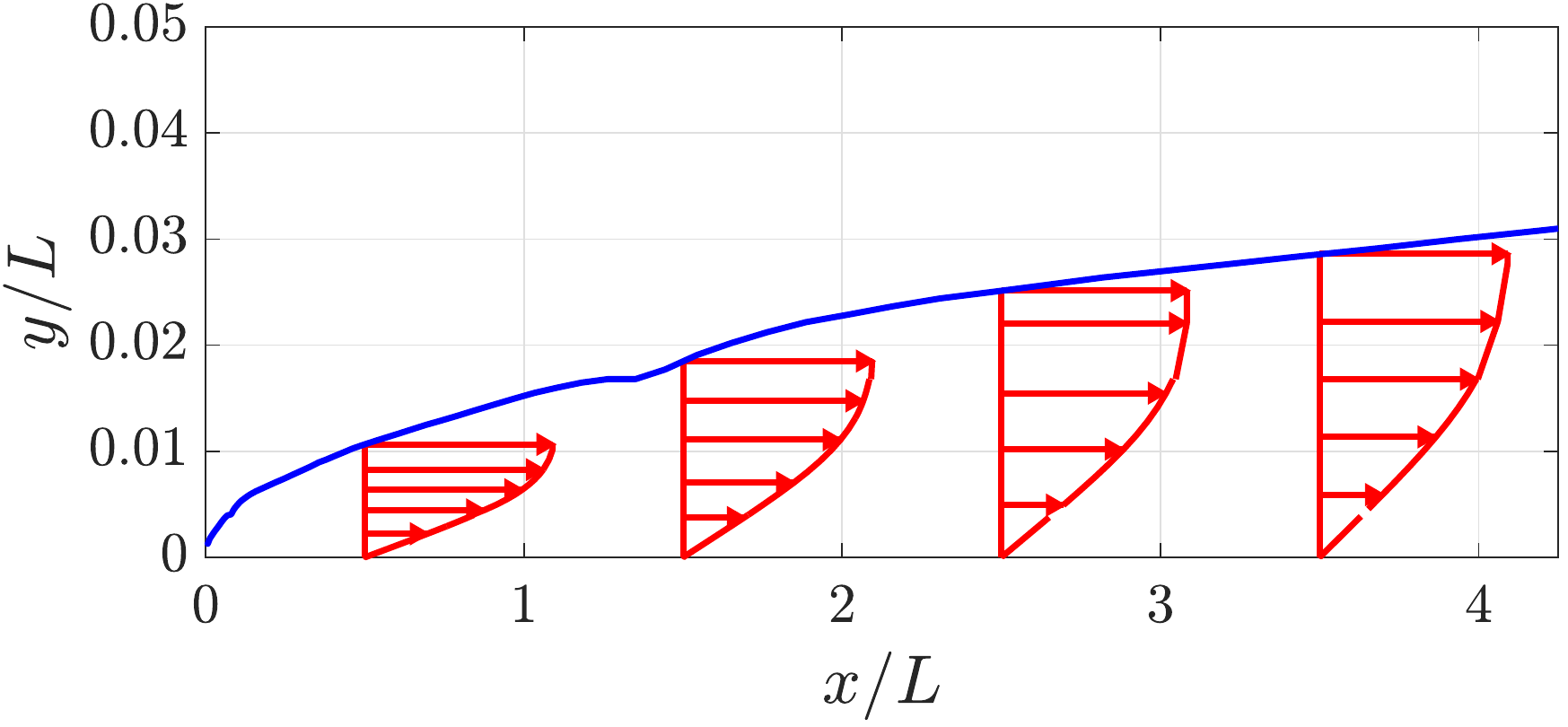}}
	\caption{Laminar flow over a flat plate - Velocity profiles along the flat plate and boundary layer thickness for the different meshes and polynomial degrees of approximation, using an HLLEM Riemann solver.}
	\label{fig:flatplate_VelocityProfiles}
\end{figure}

In a similar fashion, velocity profiles along the flat plate and detail of the boundary layer thickness are sketched in figure~\ref{fig:flatplate_VelocityProfiles} for different degrees of approximation in the different mesh refinements, computed with the HLLEM Riemann solver.
The solution is noticeably improved with mesh refinement (figure~\ref{fig:flatplate_VelocityProfiles}, top). It is worth noting that accurate approximations are achieved on the coarsest mesh using high-order polynomial approximation (figure~\ref{fig:flatplate_VelocityProfiles}, bottom). 

\begin{figure}[htbp]
	\subfloat[Boundary layer thickness, $\delta$ \label{fig:FlatPlate_errorDelta_All}]{\includegraphics[width=0.48\textwidth]{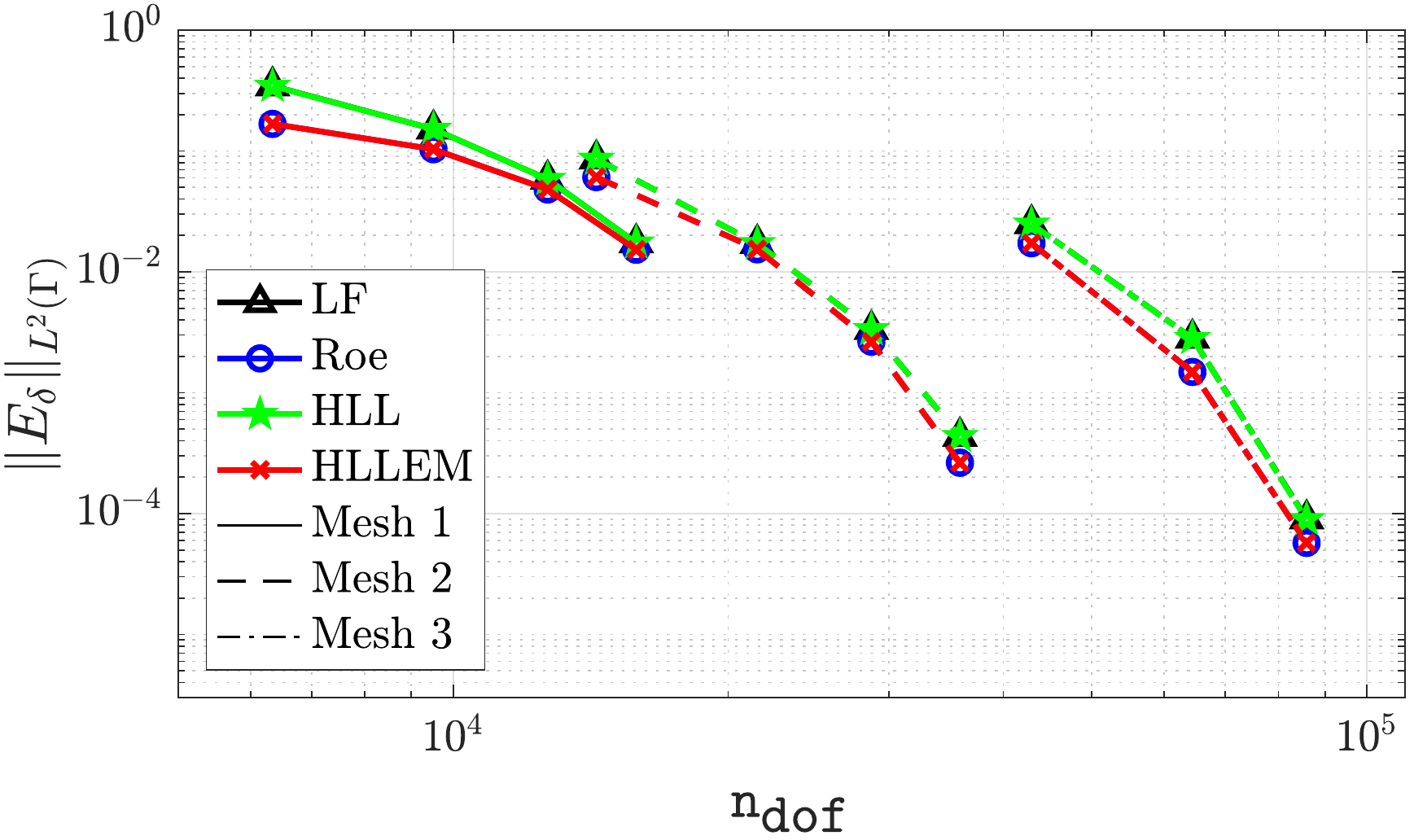}} \hfill
	\subfloat[Friction coefficient, $C_f$ \label{fig:FlatPlate_errorCf_All}]{\includegraphics[width=0.48\textwidth]{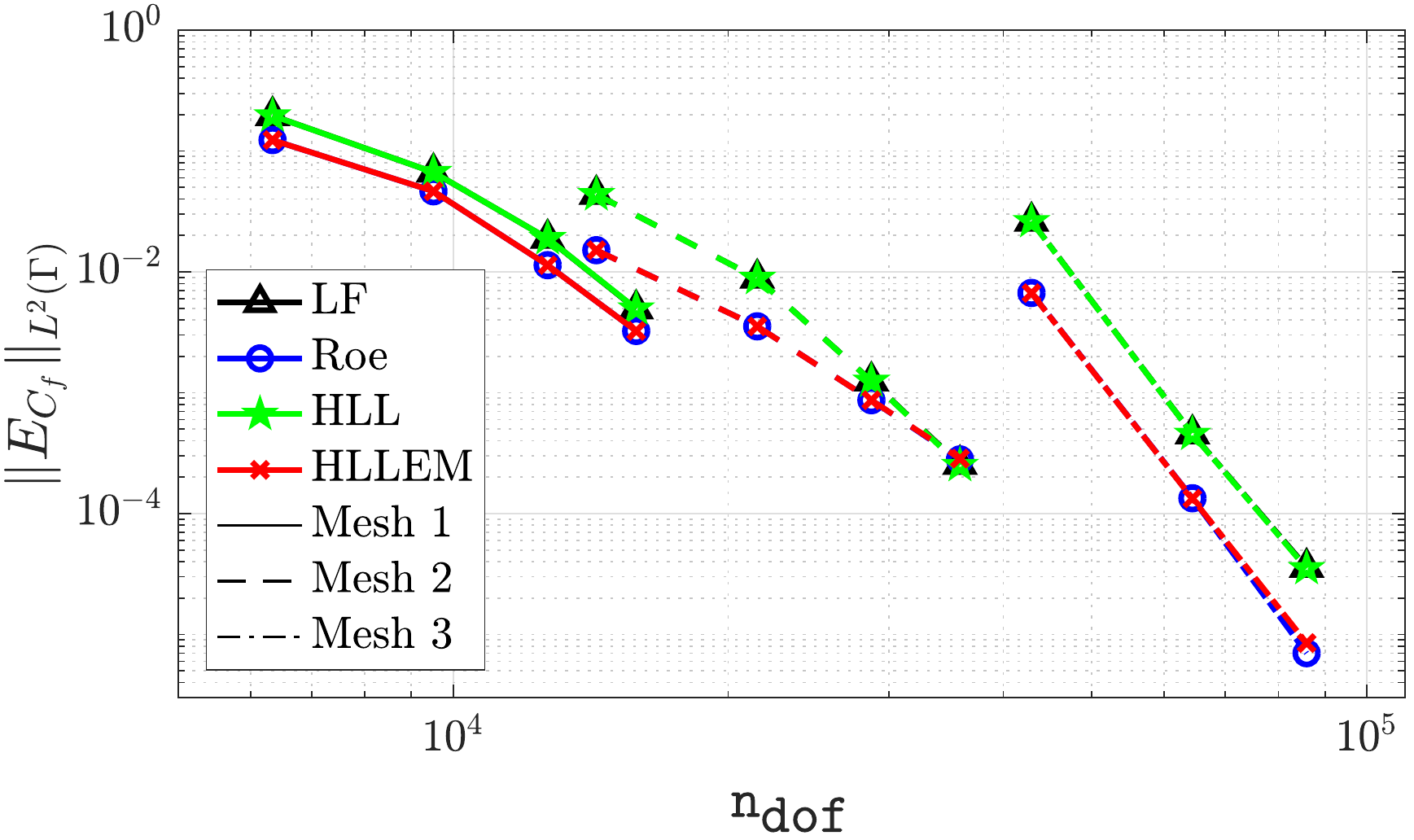}} 
	\caption{Laminar flow over a flat plate - Convergence of the relative $\eltwo$ error of the (a) boundary layer thickness and (b) friction coefficient, using Lax-Friedrichs (LF), Roe, HLL and HLLEM Riemann solvers under $k$-refinement ($k=1, \dotsc , 4$) using three different meshes.}
	\label{fig:flatplate_errors}
\end{figure}

In order to quantify the effect of the numerical dissipation introduced by the Riemann solver in the quality of the approximate solution, the $\eltwo$ error of the boundary layer thickness and of the friction coefficient is measured along the flat plate. The convergence study, shown in figure~\ref{fig:flatplate_errors}, reports the evolution of the error as a function of the number of degrees of freedom, obtained for each mesh by increasing the polynomial degree of approximation from $k=1$ up to $k=4$.

The HLLEM solution on mesh 3 using fourth-order polynomials is taken as reference solution for comparison.
The results display that Lax-Friedrichs and HLL solutions introduce higher levels of error than HLLEM and Roe. These differences are more remarkable in low order approximations, being the choice of Riemann solver a critical issue for the accuracy of the computation.
Furthermore, it is worth noticing that high-order approximations on coarse meshes provide higher accuracy than lower-order solutions with similar number of degrees of freedom, emphasising the interest for increasing the polynomial degree of approximation.

\begin{figure}[htbp]
	\subfloat[$k=1$ \label{fig:FlatPlate_Drag_k1}]{\includegraphics[width=0.33\textwidth]{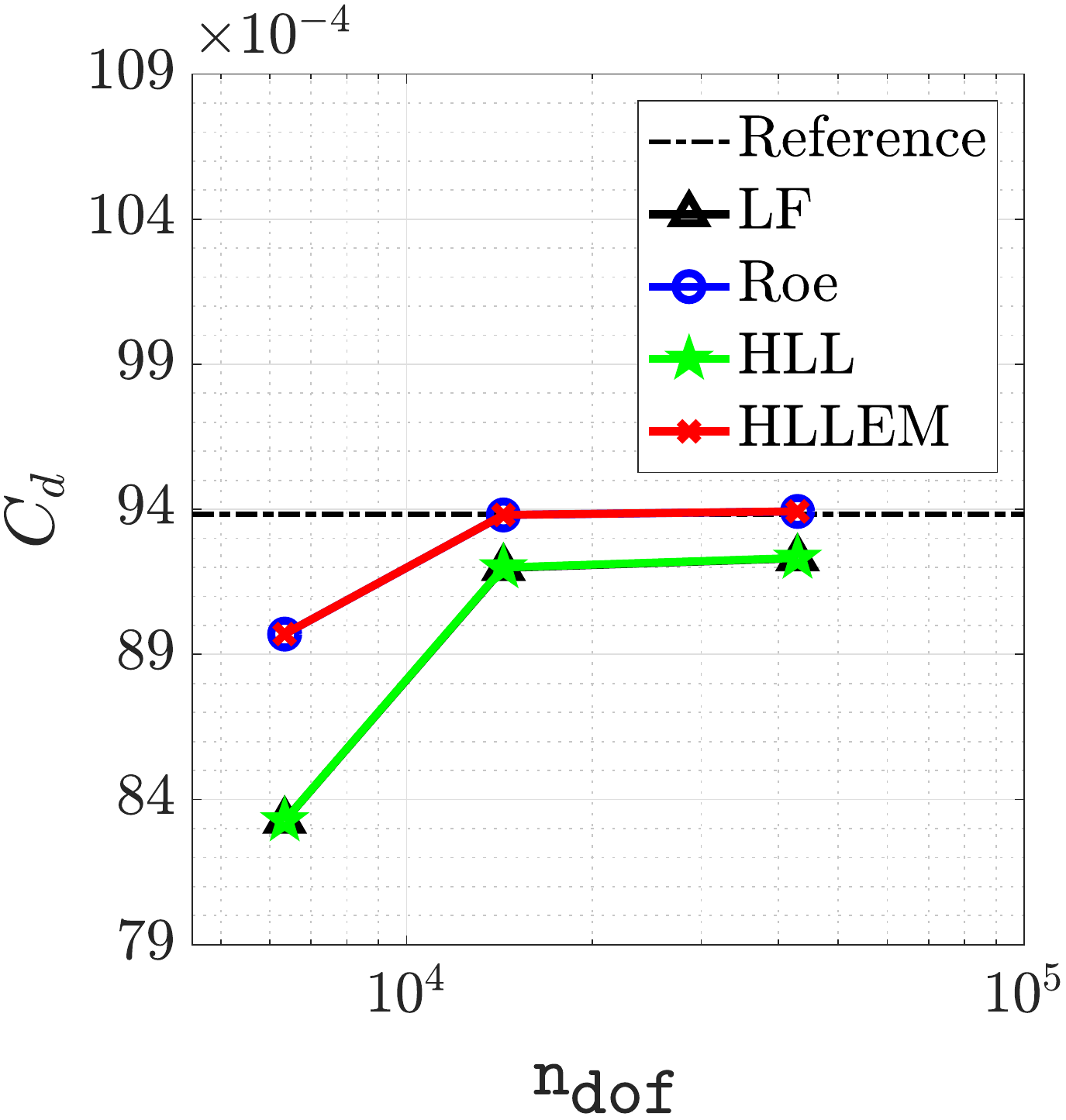}}  \hfill
	\subfloat[$k=2$ \label{fig:FlatPlate_Drag_k2}]{\includegraphics[width=0.33\textwidth]{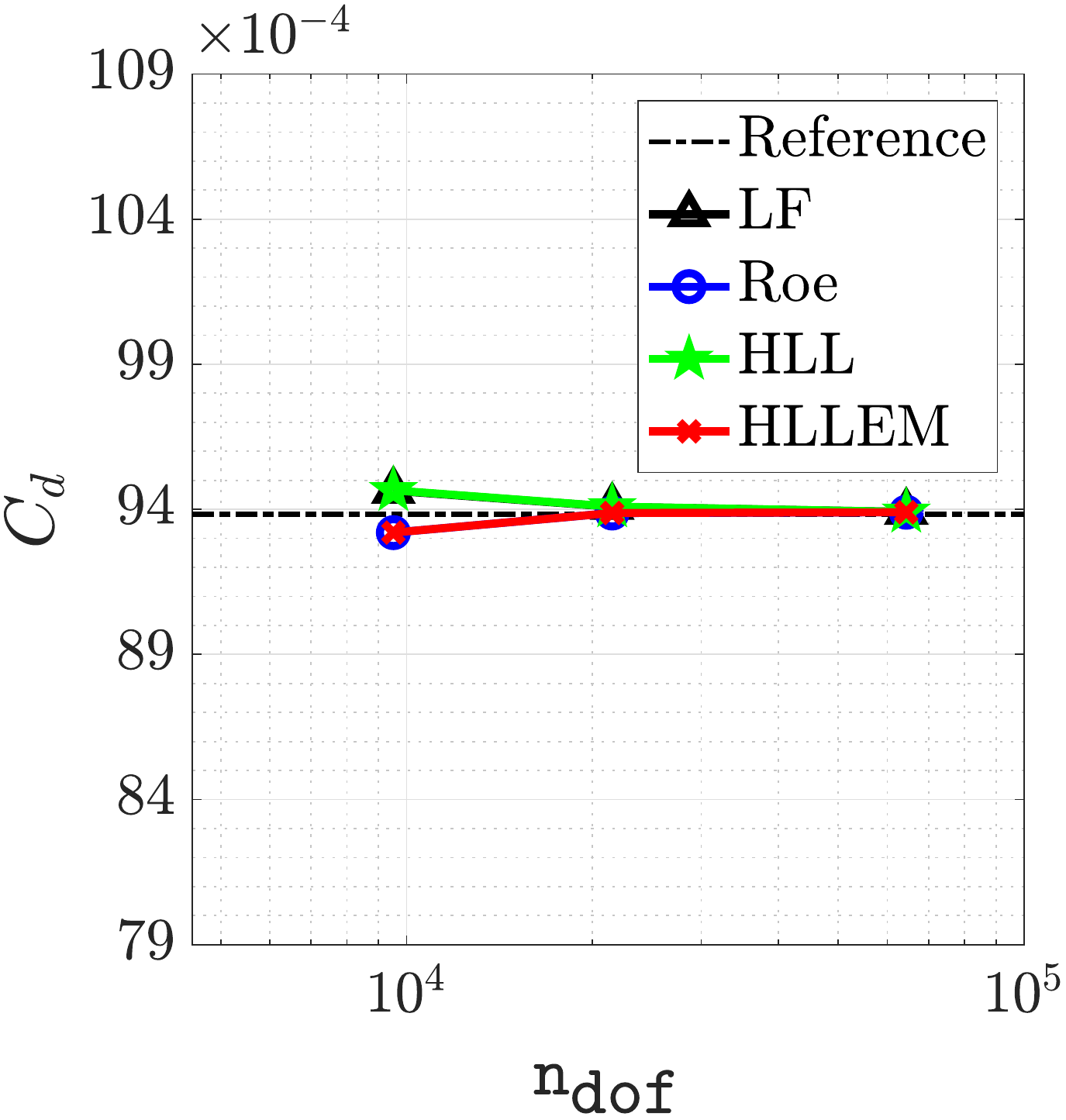}}  \hfill
	\subfloat[$k=3$ \label{fig:FlatPlate_Drag_k3}]{\includegraphics[width=0.33\textwidth]{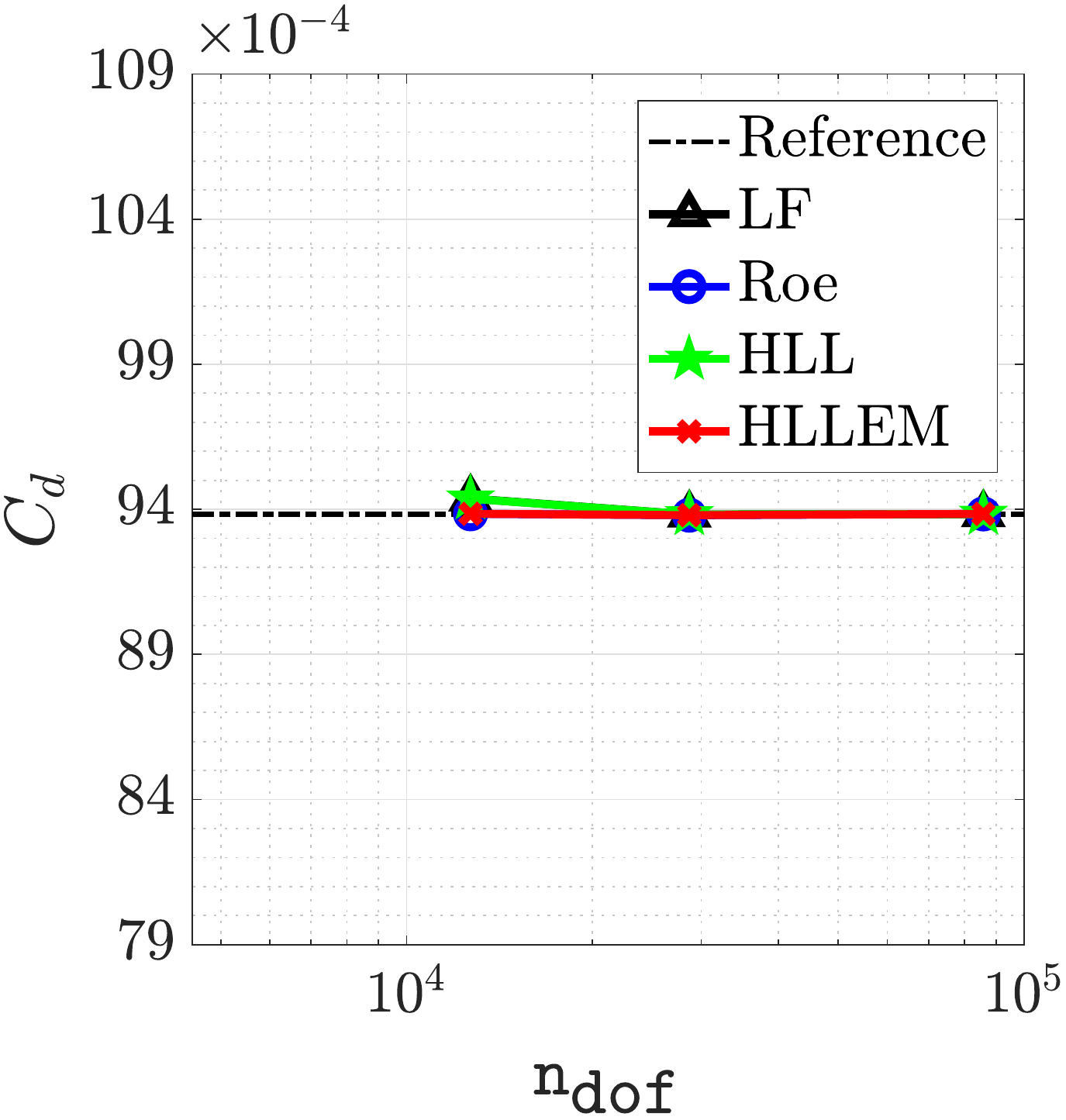}} 
	\caption{Laminar flow over a flat plate - Convergence of the drag coefficient, $C_d$, using Lax-Friedrichs (LF), Roe, HLL and HLLEM Riemann solvers under $h$-refinement using three different polynomial degrees of approximation.}
	\label{fig:flatplate_Drag}
\end{figure}

Similarly, the convergence of the drag coefficient is reported in figure~\ref{fig:flatplate_Drag}. It is confirmed that HLLEM and Roe Riemann solvers display an increased accuracy with respect to Lax-Friedrichs and HLL, which is especially evident for $k=1$. In this case, even in the coarsest mesh, the drag coefficient computed with HLLEM and Roe solutions lies within the admissible error of five drag counts, contrary to HLL and Lax-Friedrichs. As the degree of approximation increases, differences among Riemann solvers are notably reduced, due to the lower numerical dissipation introduced by HDG.

Hence, Roe and HLLEM Riemann solvers have proved to be able to resolve the flow solutions in thin boundary layers exhibiting an increased accuracy when low-order approximations are constructed. More precisely, the numerical diffusion introduced by Riemann solvers misrepresenting middle waves (i.e. Lax-Friedrichs and HLL) results critical for an accurate approximation of the solution in the boundary layer and its derived quantities. As the resolution increases, either by mesh refinement or by increasing the polynomial order of approximation, such numerical diffusion is reduced and the differences among Riemann solvers become negligible.

Henceforth, and in order to fully exploit the advantages of the presented HDG solver with the different Riemann solvers, as proved in the previous examples, only high-order approximations are considered.

\subsection{Shock treatment in inviscid flows: transonic inviscid flow over a NACA 0012 aerofoil}
\label{ssc:transonicNACA}

The transonic inviscid flow over a NACA 0012 aerofoil, at free-stream conditions $M_\infty = 0.8$ and angle of attack $\alpha = 1.25\degree$, is presented to assess the performance of the shock capturing method for inviscid flows.
This example is a classical benchmark used to verify numerical inviscid codes and implementations of shock capturing techniques, see for instance~\cite{Thibert-TGO:1979,AGARD:1985,Sevilla-SHM:2013} or the test case MTC2 in~\cite{Kroll2009}.

\begin{SCfigure}[][htbp] 
	\centering
	\includegraphics[width=0.6\textwidth]{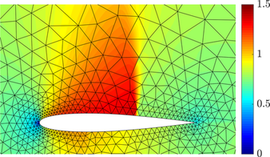}
	\caption{Transonic flow over a NACA 0012 aerofoil - Mach number distribution computed using HLL Riemann solver with polynomial degree of approximation $k=4$.}
	\label{fig:NACA_solution}
\end{SCfigure}

The steady state problem is solved via a relaxation approach with a time step $\Delta t = 10^{-1}$ such that the Courant number is $C = 22$. Convergence to the steady state is achieved when the residual of the steady terms of the continuity equation reaches $10^{-6}$ or is decreased by three orders of magnitude from its maximum value.

All Riemann solvers are equipped with the Laplacian-based shock capturing technique described in section~\ref{ssc:LaplacianBasedShockCapturing} and the value $\varepsilon_0 = 0.4$ is selected. In the case under analysis, no entropy fix is required by the Roe flux since the artificial viscosity introduced by the shock capturing strategy allows the Riemann solver to fulfill the entropy conditions. Nonetheless, it is worth remarking that the need of an entropy fix is not known \emph{a priori} and the value of the corresponding parameter $\delta$ depends upon the problem and requires to be appropriately tuned by the user.
More details will be provided in section~\ref{ssc:supersonicNACA} for the case of a supersonic flow over the NACA 0012 aerofoil.

A mesh with $1,877$ triangular elements, without any specific refinement in the shock region, is used and an approximation degree $k=4$ is considered. The far-field boundary is placed 10 chord units away from the aerofoil.

Figure~\ref{fig:NACA_solution} displays the Mach number distribution computed using the HLL Riemann solver. An accurate description of the flow around the aerofoil is obtained and the shock is precisely captured with a coarse mesh, owing to the high-order polynomial approximation constructed using the HDG framework and the shock capturing term introduced. The resolution of the shock is clearly related to the local mesh size and sharper representations may be obtained by performing local mesh refinement in the shock region, as described in~\cite{Nguyen-NP:2011}. Comparable results, not reported here for brevity, were obtained by the proposed HDG method with Lax-Friedrichs, Roe and HLLEM Riemann solvers.

The accuracy of the different numerical fluxes is thus evaluated comparing the pressure coefficient, given by
\begin{equation}
C_{p}  = \frac{p-p_{\infty}}{0.5 \rho_{\infty} v_{\infty}^2},	
\end{equation}
over the aerofoil profile.
\begin{figure}[htbp]
	\centering
	\includegraphics[width=\textwidth]{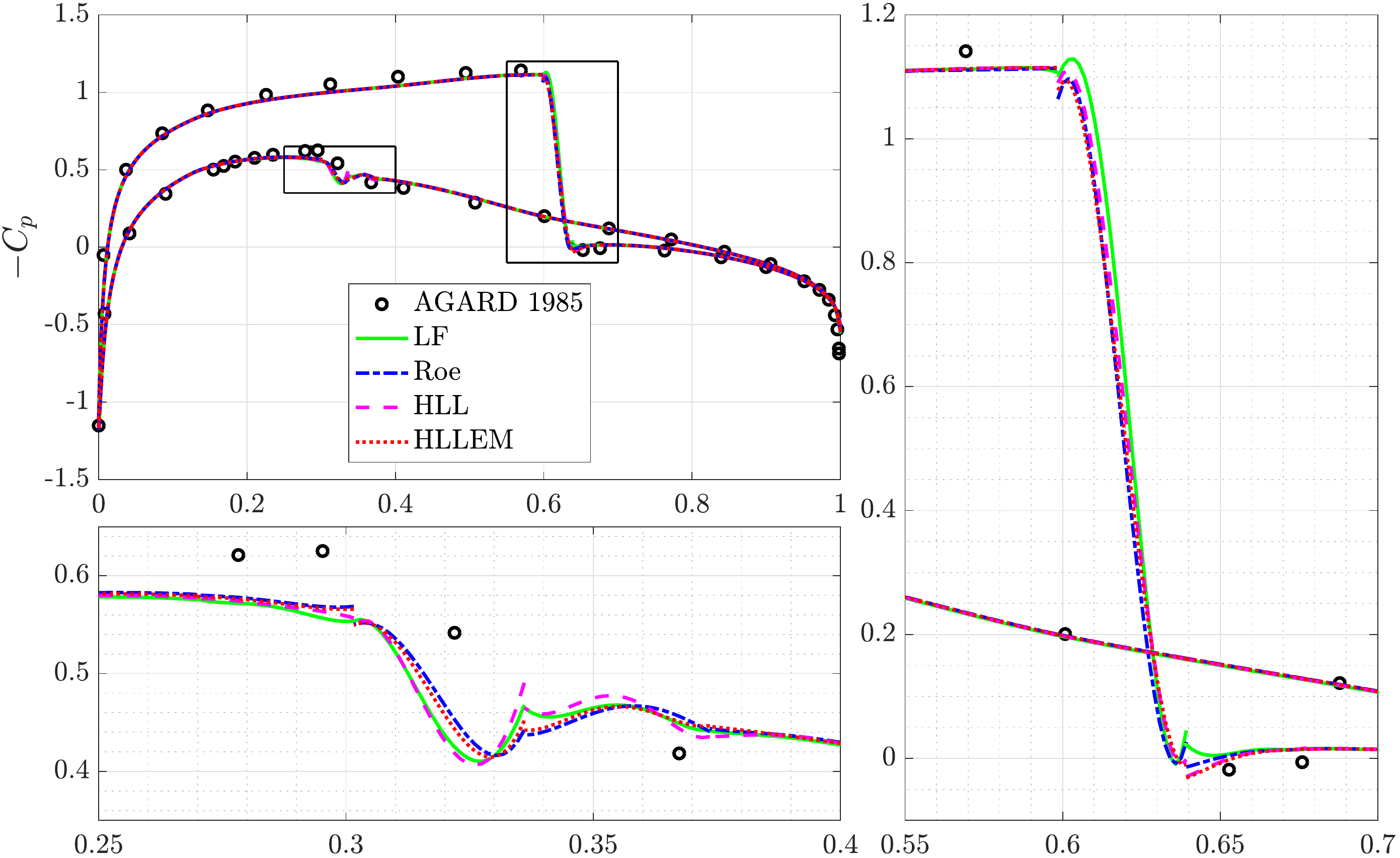}
	\caption{Transonic flow over a NACA 0012 aerofoil - Pressure coefficient around the aerofoil surface computed using different Riemann solvers with polynomial degree of approximation $k=4$ and detailed views of the lower (left) and upper (right) shocks.}
	\label{fig:NACA_Cp}
\end{figure}

A well resolved solution, in agreement with experimental data from~\cite{AGARD:1985}, is obtained using all Riemann solvers. 
The results in figure~\ref{fig:NACA_Cp} display that HLL and HLLEM Riemann solvers provide an approximation without oscillations and with accuracy similar to the one of the Roe numerical flux near the upper, stronger shock. 
It is worth noting that the jumps appearing at the extrema of the shock region are due to the discontinuous nature of the HDG approximation.
The lower, weaker shock, is reproduced less precisely by the four Riemann solvers. In this case, HLL presents a behaviour closer to the Lax-Friedrichs solution, whereas HLLEM and Roe produce a similar approximation.

Accordingly, the lift and drag coefficients reported in table~\ref{tb:transonicNACA_LiftDrag} allow to quantify very little differences among Riemann solvers. 
\begin{table} [htbp]
	\centering
	\caption{Transonic inviscid flow over a NACA 0012 aerofoil - Lift and drag coefficients for different Riemann solvers using a polynomial degree of approximation $k=4$.}
	\begin{tabular}{ L{1cm} C{3cm} C{3cm} C{2.25cm} C{2.25cm}  }
		\toprule
		 &Lax-Friedrichs & Roe & HLL  & HLLEM \\
		\midrule
		$C_l$  & $0.320$ & $0.314$ & $0.317$ & $0.315$ \\
		\midrule
		$C_d$ & $0.0193$ & $0.0190$ & $0.0192$ & $0.0191$ \\
		\bottomrule
	\end{tabular}
	\label{tb:transonicNACA_LiftDrag}
\end{table}

The obtained values lie between 25 and 35 lift and drag counts with respect to typical reference values~\cite{Thibert-TGO:1979}.
Note that the precision of the aerodynamic coefficients is strongly dependent on the location of the far-field boundary, as reported by Yano and Darmofal~\cite{YanoDarmofal2012,Wang2013}.
In particular, for such kind of comparisons, far-field boundaries are tipically located at distances from 50 up to $10^4$ chord lengths from the  aerofoil~\cite{Balan-BWM:2015,May-WBMS-14,YanoDarmofal2012,Wang2013,Ekelschot2017}.

Finally, the entropy production is considered for this non-isentropic case. In this context, such quantity allows to estimate the numerical dissipation introduced in the upstream region before the shock and the entropy produced by the artificial viscosity.
\begin{figure}[htbp]
	\centering
	\subfloat[HLL, sensor activation \label{fig:NACA_HLL_activation}]{\includegraphics[width=0.4\textwidth]{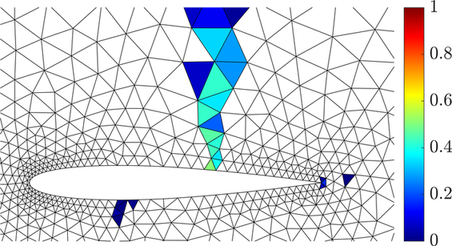}} \quad
	\subfloat[HLL, entropy production \label{fig:NACA_HLL_entropy}]{\includegraphics[width=0.4\textwidth]{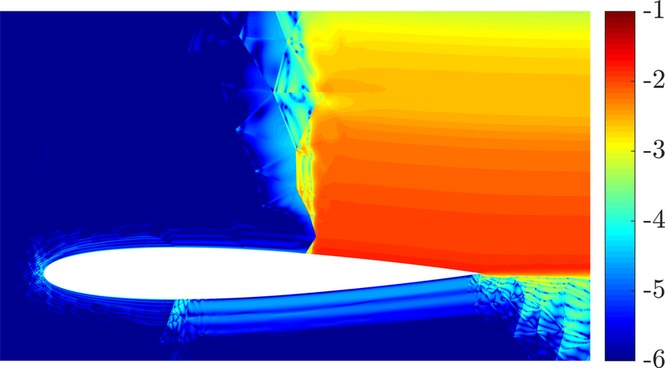}}  \\
	\subfloat[HLLEM, sensor activation \label{fig:NACA_HLLEM_activation}]{\includegraphics[width=0.4\textwidth]{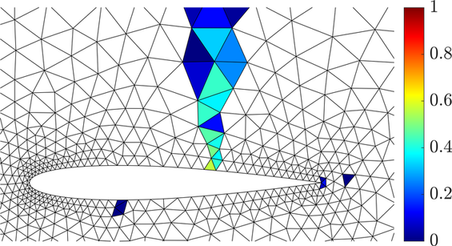}} \quad
	\subfloat[HLLEM, entropy production \label{fig:NACA_HLLEM_entropy}]{\includegraphics[width=0.4\textwidth]{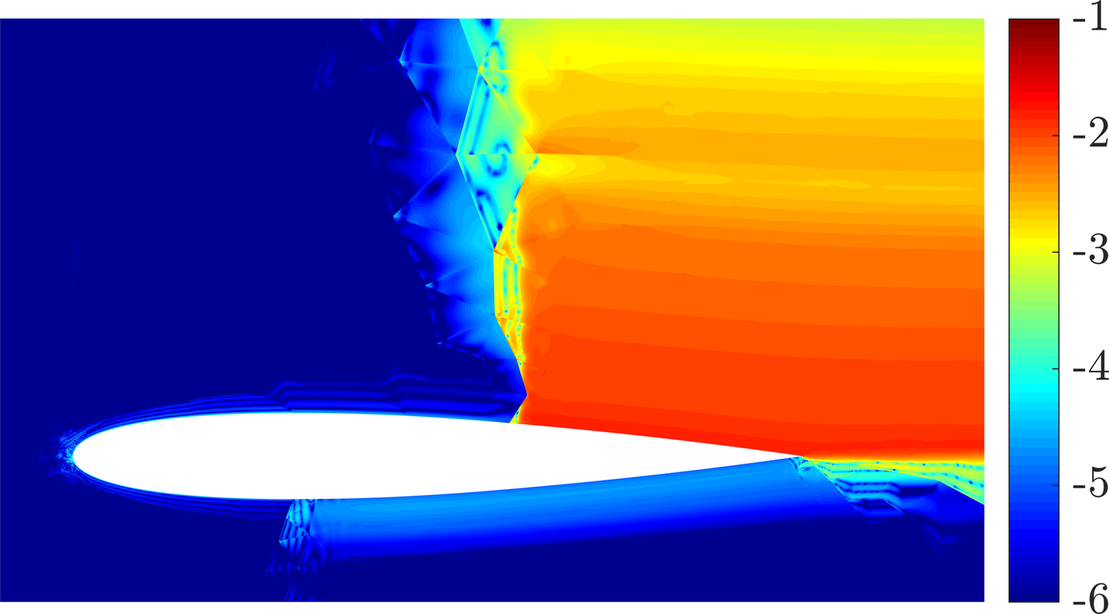}}  \\
	\subfloat[LF, sensor activation \label{fig:NACA_LF_activation}]{\includegraphics[width=0.4\textwidth]{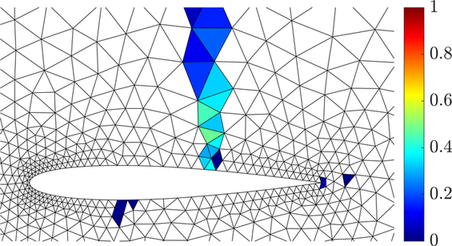}} \quad
	\subfloat[LF, entropy production \label{fig:NACA_LF_entropy}]{\includegraphics[width=0.4\textwidth]{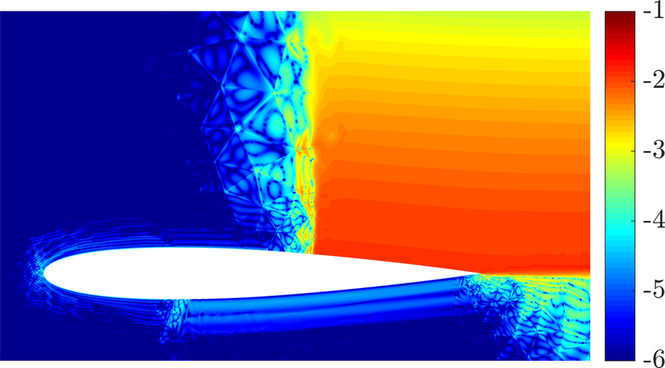}}  \\	
	\subfloat[Roe, sensor activation \label{fig:NACA_Roe_activation}]{\includegraphics[width=0.4\textwidth]{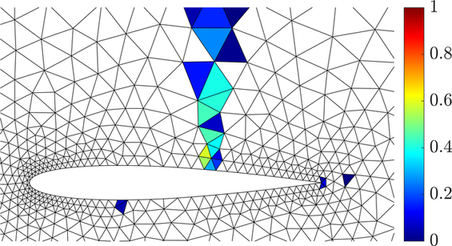}} \quad
	\subfloat[Roe, entropy production \label{fig:NACA_Roe_entropy}]{\includegraphics[width=0.4\textwidth]{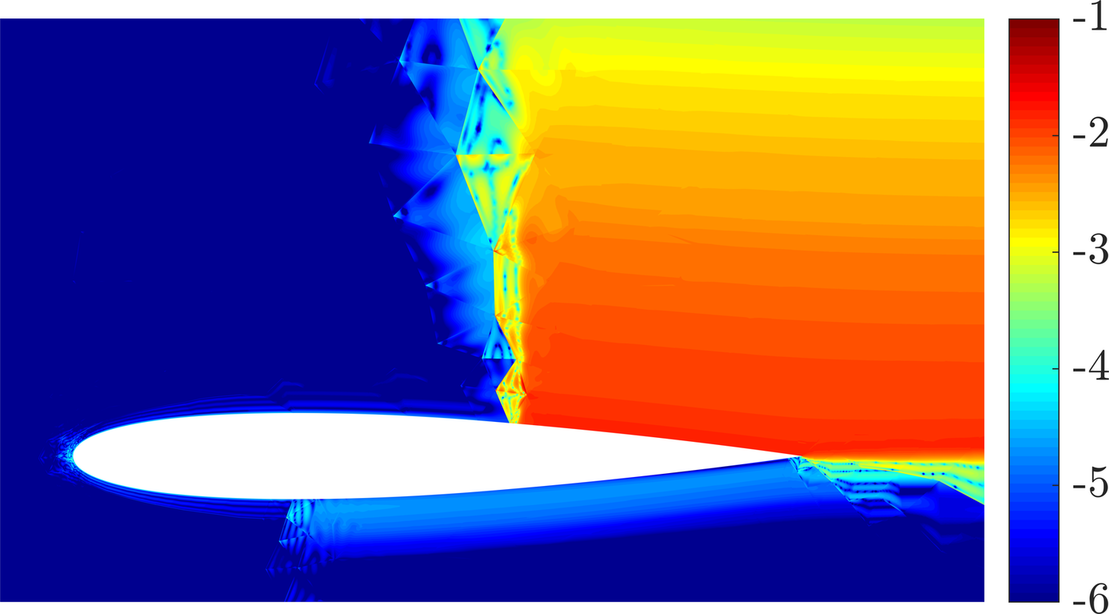}} 
	\caption{Transonic flow over a NACA 0012 aerofoil - Regions of activation of the shock sensor (left) and entropy production in logarithmic scale (right) for HLL (a-b), HLLEM(c-d), Lax-Friedrichs (LF, e-f) and Roe (g-h) Riemann solvers using a polynomial degree of approximation $k=4$.}
	\label{fig:NACA_SensorEntropy}
\end{figure}

On the one hand, the results in figure~\ref{fig:NACA_SensorEntropy} show that the regions of activation of the sensor are almost identical for the four Riemann solvers.
On the other hand, the different amount of numerical dissipation introduced by the numerical fluxes is responsible for the production of entropy. 
As observed in figure~\ref{fig:NACA_Cp}, HLL, HLLEM and Roe Riemann solvers present a similar behaviour in the vicinity of the upper, stronger shock, where comparable approximations are achieved.
On the contrary, the Lax-Friedrichs numerical flux introduces the largest amount of numerical dissipation in this region, as shown in figure~\ref{fig:NACA_LF_entropy}.
In the vicinity of the weaker shock on the lower part of the aerofoil, the four Riemann solvers show a similar entropy production.
Finally, Roe and HLLEM solvers provide the most accurate results in the region near the trailing edge, where the HLL and the Lax-Friedrichs numerical fluxes introduce extra dissipation.

This example demonstrates an overall good performance of the Laplacian-based shock capturing method for inviscid compressible flows. Furthermore, no significant differences are observed among the Riemann solvers using a high-order approximation of order $k=4$. In particular, the four numerical fluxes lead to similar approximate solutions, as reported with the aerodynamic measures of lift, drag and pressure coefficients, while displaying an accurate and positively conservative treatment of the shock waves.

\subsection{Positivity-preserving properties in presence of shocks: supersonic inviscid flow over a NACA 0012 aerofoil}
\label{ssc:supersonicNACA}

The second example of inviscid flow around a NACA 0012 aerofoil consists of a supersonic flow at a free-stream Mach number $\Minf = 1.5$ and zero angle of attack \cite{Persson-PP:2006, Balan-BWM:2015}.

This supersonic test case challenges the performance of the proposed Riemann solvers in HDG in capturing solutions involving shocks and sharp gradients while ensuring positivity-preserving properties using high-order approximations. It is worth noticing that, in such case, Riemann solvers may fail to provide physically admissible solutions, leading to a violation of the positiveness of the approximate density and pressure fields~\cite{Fleischmann-FAHA:2019,Peery1988,Quirk1994}.

The computational mesh described in the previous case~\ref{ssc:transonicNACA}, consisting of 1,877 triangular elements and a far-field boundary placed at 10 chord units away from the aerofoil, is employed for the simulation.
A time step $\Delta t = 8 \times 10^{-2}$ is considered to advance in time and the corresponding Courant number is $C = 20$. Convergence to the steady state is achieved when the residual of the steady terms of the continuity equation reaches $10^{-6}$ or is decreased by three orders of magnitude from its maximum value. 
The shock treatment is handled by means of the Laplacian-based technique discussed in section~\ref{ssc:LaplacianBasedShockCapturing}, with a maximum value of artificial viscosity $\varepsilon_0 = 1$.

\begin{SCfigure}[][htbp] 
	\centering
	\includegraphics[width=0.6\textwidth]{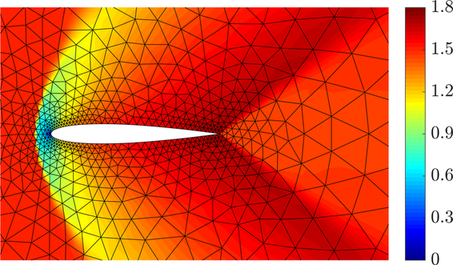}
	\caption{Supersonic flow over a NACA 0012 aerofoil - Mach number distribution computed using an HLL Riemann solver with polynomial degree of approximation $k=4$.}
	\label{fig:NACAsupersonic}
\end{SCfigure}
The Mach number distribution computed using the HLL Riemann solver with a polynomial degree of approximation $k = 4$ is presented in figure~\ref{fig:NACAsupersonic}. The method is able to accurately capture the physics of the problem, even on a coarse mesh, owing to the high-order functional discretisation introduced by the HDG scheme.

\begin{SCfigure}[][htbp] 
		\includegraphics[width=0.65\textwidth]{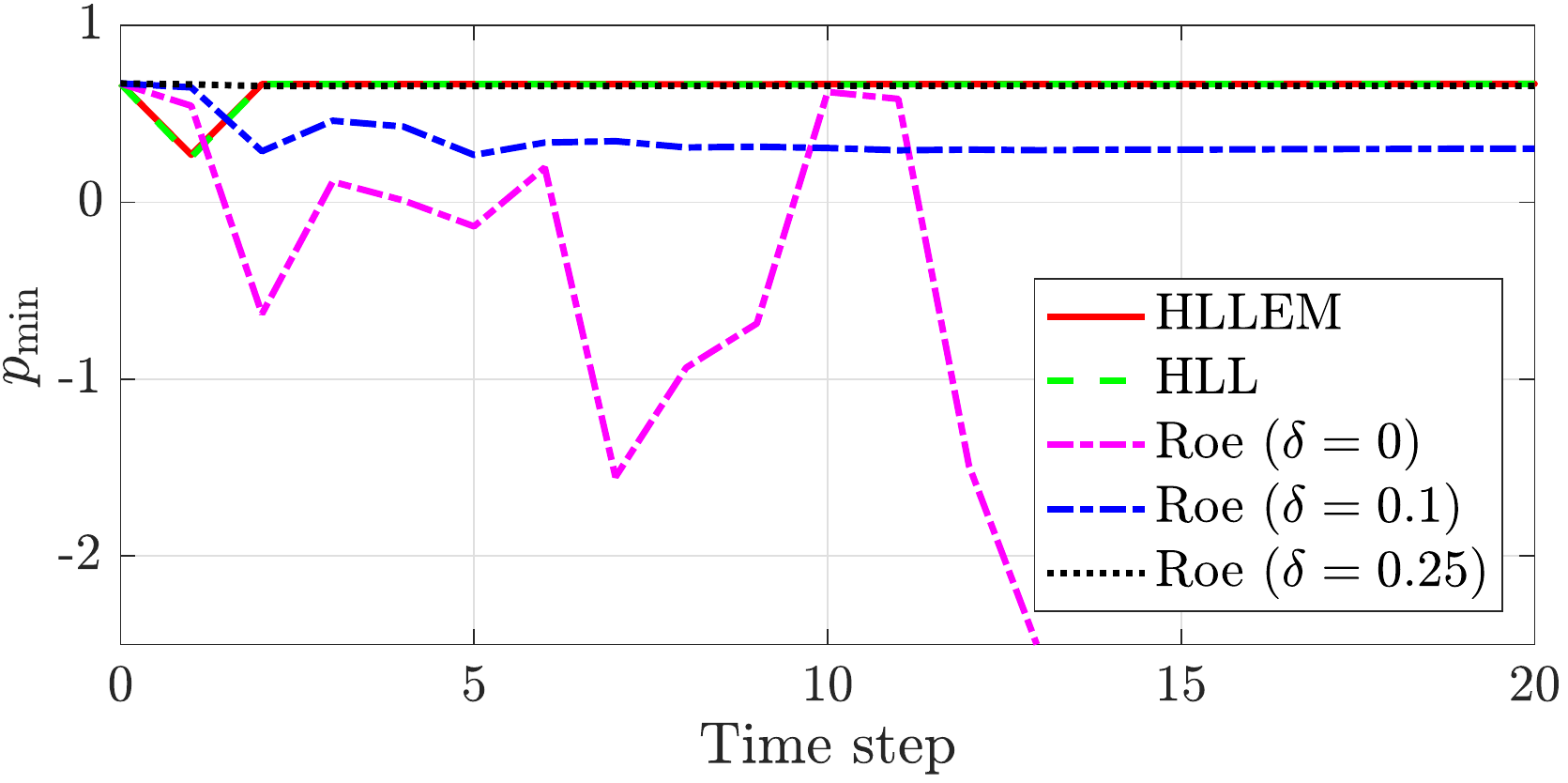}
		\caption{Supersonic flow over a NACA 0012 aerofoil - Minimum nodal value of the pressure computed using the different Riemann solvers with polynomial degree of approximation $k=4$.}
		\label{fig:pressureNonPhys}
\end{SCfigure}

This supersonic problem is especially challenging since it features an abrupt shock in front of the aerofoil and allows to test the positivity properties of the approximate solution.
For this purpose, the performance of the Roe Riemann solver is compared with those of the HLL family. Figure~\ref{fig:pressureNonPhys} shows the minimum nodal value of the pressure computed using the Roe numerical flux with different values of the HH entropy fix as well as with HLL and HLLEM.

In the case with no entropy fix ($\delta = 0$), the Roe solver displays an insufficient numerical dissipation. After few iterations, negative values of the pressure are computed, leading to a nonphysical solution. This error is amplified from one time step to the following ones and rapidly leads to the divergence of the Newton-Raphson algorithm employed to solve the nonlinear problem. 
To remedy this issue, inherent to the Roe Riemann solver, an HH entropy fix with an empirically tuned value of the threshold parameter $\delta$ is considered. It is worth emphasising that the tuning of such parameter is problem-dependent. With a setting of $\delta = 0.1$, the HDG method with Roe Riemann solver converges to a steady state solution including some nonphysical undershoots in the pressure and density fields, giving rise to overshoots in the Mach distribution.

\begin{figure}[htbp]
	\centering
	\subfloat[HLL, Mach \label{fig:NACAsup_HLL_Mach_detail}]{\includegraphics[width=0.4\textwidth]{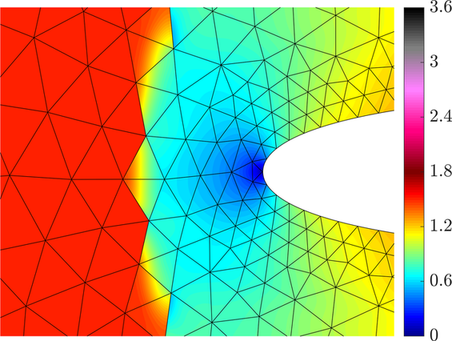}} \quad
	\subfloat[HLL, artificial viscosity \label{fig:NACAsup_HLL_viscosity_detail}]{\includegraphics[width=0.415\textwidth]{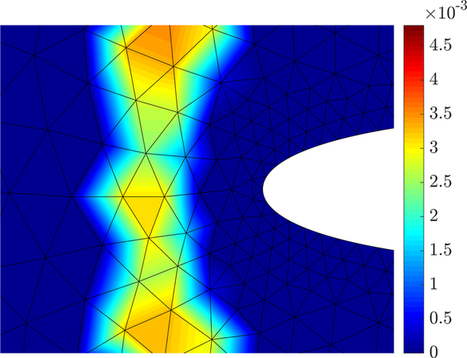}} \\
	\subfloat[Roe HH-EF $\delta=0.1$, Mach \label{fig:NACAsup_Roe1_Mach_detail}]{\includegraphics[width=0.4\textwidth]{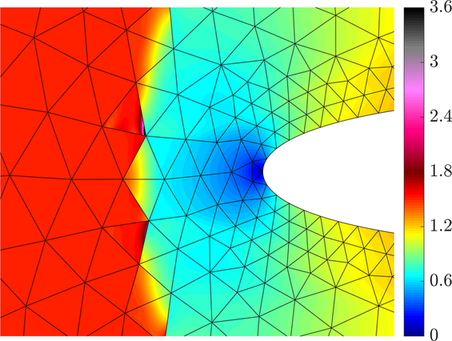}} \quad
	\subfloat[Roe HH-EF $\delta=0.1$, artificial viscosity \label{fig:NACAsup_Roe1_viscosity_detail}]{\includegraphics[width=0.415\textwidth]{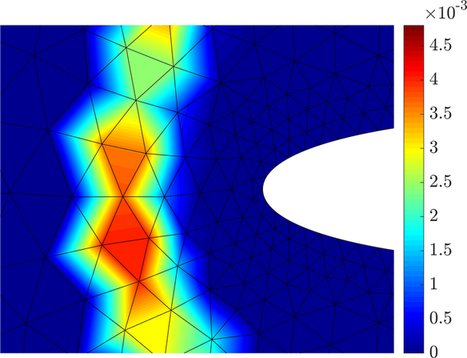}}
	\caption{Supersonic flow over a NACA 0012 aerofoil - Detail of the Mach number distribution (left) and corresponding artificial viscosity (right) in the front shock near the leading edge computed using HLL (top) and Roe Riemann solver with HH entropy fix with threshold parameter $\delta = 0.1$ (bottom) with polynomial degree of approximation $k=4$.}
	\label{fig:NACAsupersonicRoe}
\end{figure}

Precisely, the corresponding Mach number distribution computed using the Roe numerical flux with entropy fix parameter $\delta = 0.1$ is reported in figure~\ref{fig:NACAsupersonicRoe} to illustrate such spurious oscillations appearing in the region in front of the shock (Fig.~\ref{fig:NACAsup_Roe1_Mach_detail}). Such oscillations appear despite the artificial viscosity introduced in the corresponding elements, as displayed in figure~\ref{fig:NACAsup_Roe1_viscosity_detail}. Hence, this value of the HH entropy fix parameter leads to insufficient stabilisation and a higher threshold needs to be introduced. 


\begin{remark}
	It is worth noting that the colour scale of figure~\ref{fig:NACAsupersonicRoe} keeps the same gradation of colours of figure~\ref{fig:NACAsupersonic} for the interval $M \in [0, 1.8]$ but extends up to $M = 3.6$ to visualise the peak values achieved by the overshoots in the Roe solution.  
\end{remark}

Such numerical issues are fixed by increasing the threshold value $\delta$ of the HH entropy fix. Numerical results showed that a value $\delta = 0.25$ or larger allows the high-order HDG solver to achieve a physically admissible solution with no overshoots, as reported in figure~\ref{fig:pressureNonPhys}.
Nonetheless, in case of exceeding the threshold value of the entropy fix, the associated numerical dissipation of the Roe Riemann solver is increased, turning the solver overdiffusive. In the limit, $\delta \to \lamax$, the Lax-Friedrichs Riemann solver is obtained.
On the contrary, HLL and HLLEM numerical fluxes provide a robust approximation with no oscillations without the need of any user-defined entropy fix.

\begin{figure}[htbp]
	\centering
	\subfloat[HLL, sensor activation \label{fig:NACAsup_HLL_sensor}]{\includegraphics[width=0.33\textwidth]{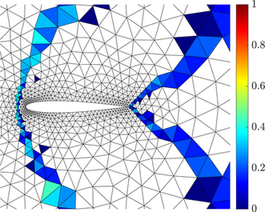}} \hfill
	\subfloat[HLLEM, sensor activation \label{fig:NACAsup_HLLEM_sensor}]{\includegraphics[width=0.33\textwidth]{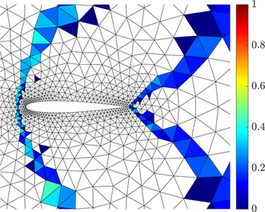}} \hfill
	\subfloat[LF, sensor activation \label{fig:NACAsup_LF_sensor}]{\includegraphics[width=0.33\textwidth]{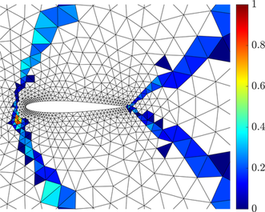}} \\
	\subfloat[HLL, entropy production \label{fig:NACAsup_HLL_entropy}]{\includegraphics[width=0.33\textwidth]{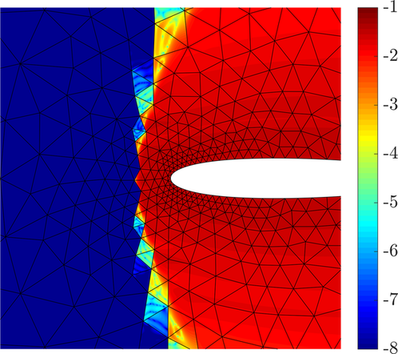}} \hfill
	\subfloat[HLLEM, entropy production \label{fig:NACAsup_HLLEM_entropy}]{\includegraphics[width=0.33\textwidth]{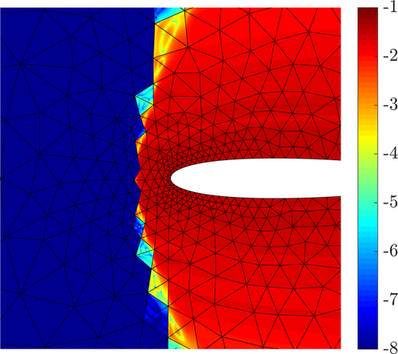}} \hfill
	\subfloat[LF, entropy production \label{fig:NACAsup_LF_entropy}]{\includegraphics[width=0.33\textwidth]{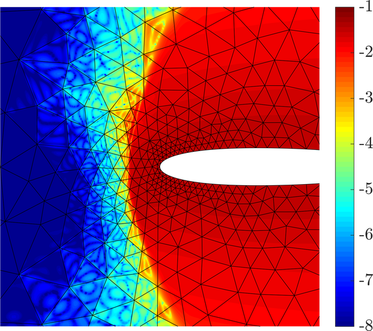}}
	\caption{Supersonic flow over a NACA 0012 aerofoil - Regions of activation of the shock sensor (top) and entropy production in logarithmic scale (bottom) for HLL (left), HLLEM (middle) and Lax-Friedrichs (LF, right) Riemann solvers using polynomial degree of approximation $k=4$.}
	\label{fig:NACAsupersonic_entropy}
\end{figure}

The entropy production is then examined for this non-isentropic case. In this context, such quantity allows to estimate the numerical dissipation introduced in the upstream region before the shock and the entropy produced by the artificial viscosity.
The map of the entropy production is reported in figure~\ref{fig:NACAsupersonic_entropy} for the HLL, HLLEM and Lax-Friedrichs numerical fluxes.
The results display that HLL-type Riemann solvers introduce a limited amount of numerical dissipation in the vicinity of the front shock.
On the contrary, the Lax-Friedrichs solver is responsible for a large entropy production in the shock region, confirming its over-diffusive nature also in supersonic problems.
Figure~\ref{fig:NACAsupersonic_entropy} also confirms that the shock-capturing sensor is activated in the same regions independently on the Riemann solver considered.

Finally, the accuracy of the approximate solutions corresponding to the different Riemann solvers is quantitatively evaluated with respect to the error in the lift coefficient. It is well-known that a symmetric aerofoil subject to a flow at zero angle of attack produces no lift force. 
Table~\ref{tb:NACA_LiftError} gathers the lift coefficient computed with the different Riemann solvers.
The HLL-type numerical fluxes, i.e., HLL and HLLEM, are the most accurate in such computation, with a lift coefficient laying at 5 and 6 lift counts from the reference value, respectively.
Both the lift coefficient computed by Roe with an entropy fix $\delta = 0.25$ and by Lax-Friedrichs (LF) feature a higher error of 8 lift counts with respect to the reference value.

\begin{table} [htbp]
	\centering
	\caption{Supersonic inviscid flow over a NACA 0012 aerofoil - Lift coefficient for different Riemann solvers using a polynomial degree of approximation $k=4$.}
	\begin{tabular}{ L{1cm} C{2cm} C{3cm} C{3cm} C{2.25cm} C{2.25cm}  }
		\toprule
		& Reference &Lax-Friedrichs & Roe ($\delta = 0.25$) & HLL  & HLLEM \\
		\midrule
		$C_l$ & $0$ & $-0.008$ & $-0.008$ & $-0.005$ & $-0.006$ \\
		\bottomrule
	\end{tabular}
	\label{tb:NACA_LiftError}
\end{table}

This example involving a strong shock wave illustrates the ability of HLL-type Riemann solvers, such as HLL and HLLEM, of guaranteeing positivity and thus producing physically admissible solutions in a robust and parameter-free strategy, in contrast with Roe Riemann solver.

\subsection{Shock wave/boundary layer interaction}
\label{ssc:SWBLI}

The next example considers the strong interaction between a shock wave and a laminar boundary layer. Such interaction is a basic phenomenon of viscous-inviscid interaction that happens when a shock impinges on the boundary layer producing separation in it. In such a case, the shock, instead of reflecting off the wall, turns into a combination of an expansion fan  at the edge of the boundary layer plus two compression waves around the separation and reattachment points \cite{Hakkinen1959,Katzer1989}.

\begin{figure}[htbp]
	\centering
	\subfloat[ \label{fig:SWBLI_Domain1}]{\includegraphics[width=0.73\textwidth]{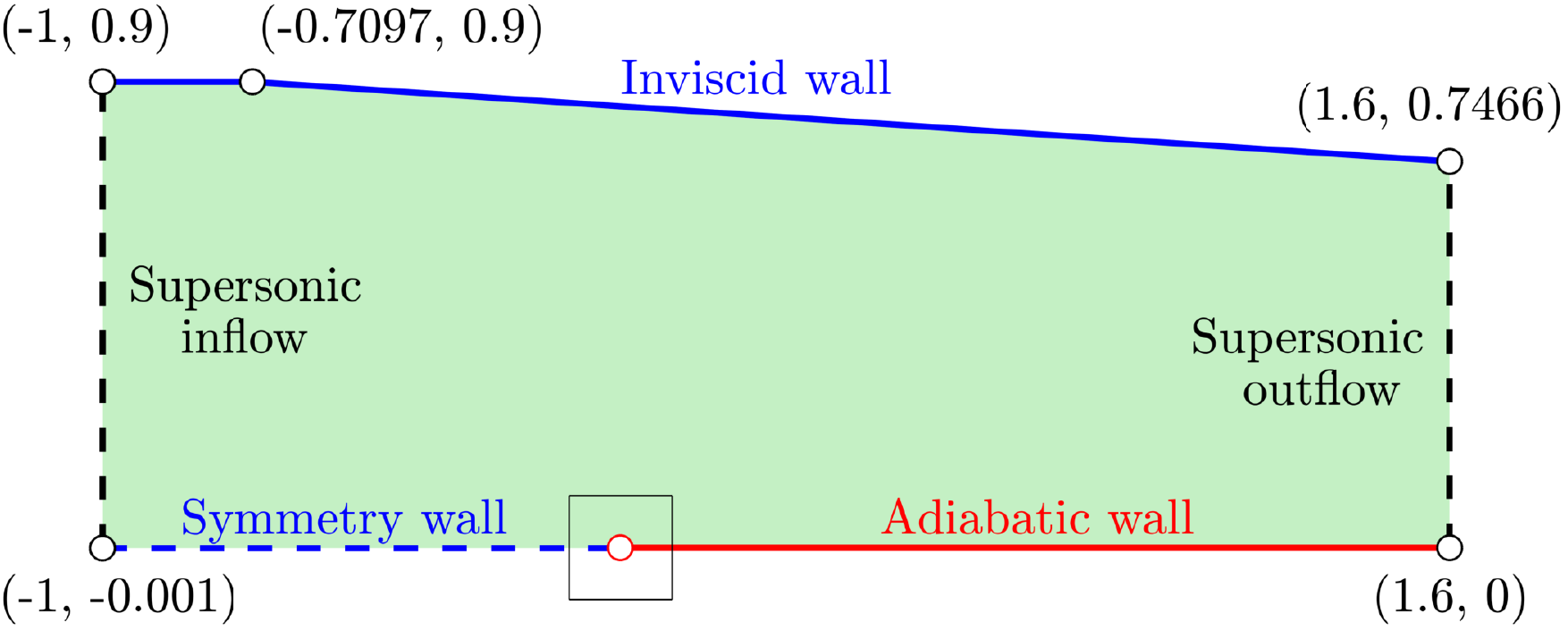}} \hfill
	\subfloat[ \label{fig:SWBLI_Domain2}]{\includegraphics[width=0.26\textwidth]{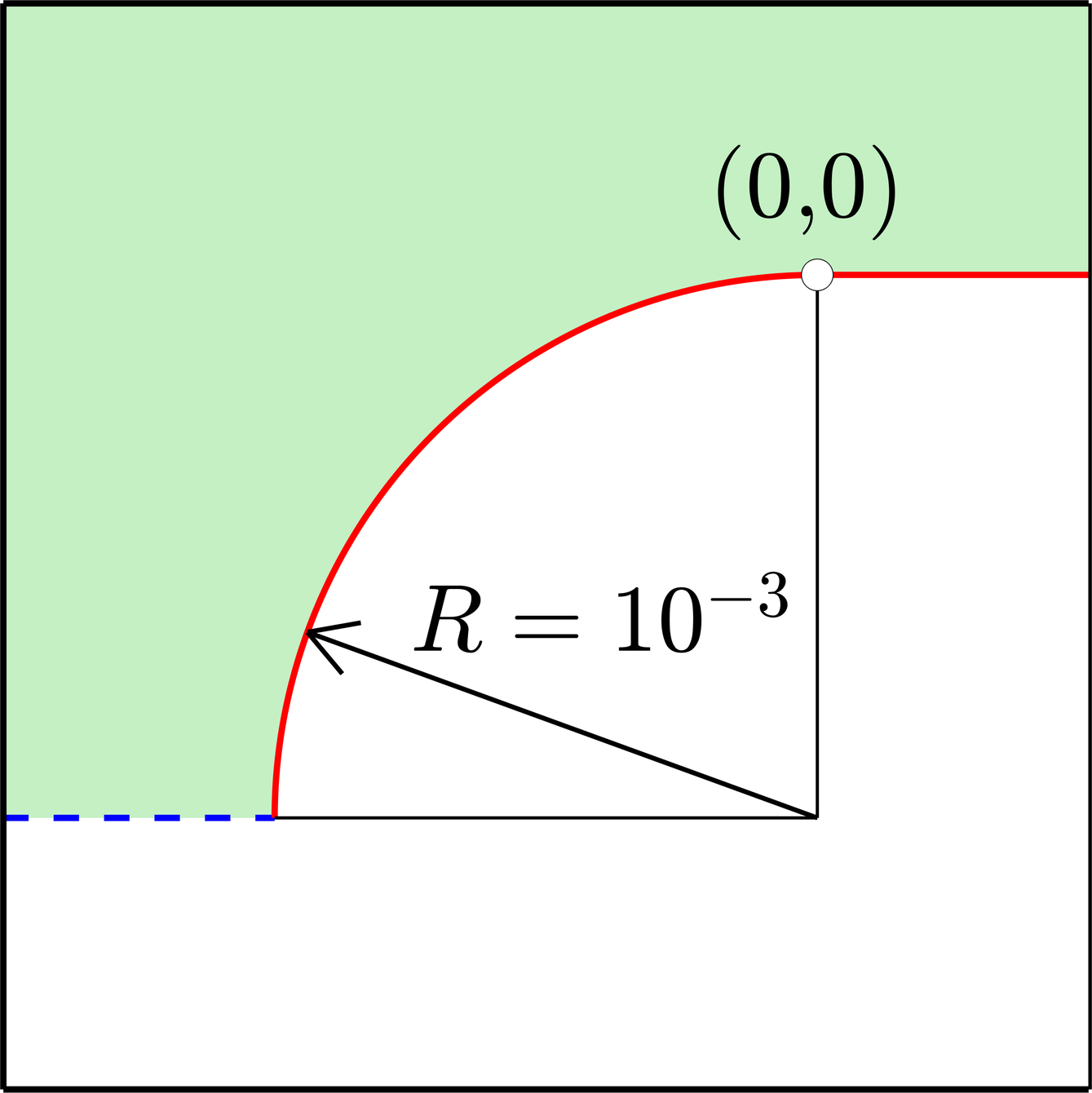}}
	\caption{Shock wave/boundary layer interaction - (a) Geometry and boundary conditions and (b) detail of the fillet at the leading edge.}
	\label{fig:SWBLI_Domain}
\end{figure}

The setup of this test case replicates the one introduced by Degrez et al. \cite{Degrez1987} and later reproduced by Moro et al. \cite{Moro2017} using a high-order HDG discretisation with $k=3$. It consists of a flat plate and a shock generator mounted inside a stream at $\Minf = 2.15$ and $\Rey = 10^5$.
A sketch of the geometry for a characteristic length of $L=1$ and the corresponding boundary conditions are detailed in figure~\ref{fig:SWBLI_Domain}.
It is worth noticing that a fillet is introduced at the leading edge in order to avoid the singularity at this location. 

\begin{figure}[htbp]
	\centering
	\includegraphics*[width=0.6\textwidth]{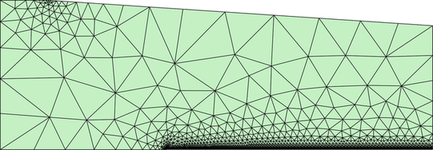}
	\caption{Shock wave/boundary layer interaction - Computational mesh.}
	\label{fig:SWBLI_Mesh}
\end{figure}

The computational mesh, depicted in figure~\ref{fig:SWBLI_Mesh}, is composed of 3,379 triangular elements of degree $k=3$. The boundary layer mesh consists of $\nlayers = 12$ layers of elements with a growing rate $r=1.4$ and the first layer located at a height of $h_0/L = 2.5\cdot 10^{-4}$. In addition, the mesh is refined at the leading edge and $\ndiv = 80$ divisions are defined along the plate.

\begin{figure}[htbp]
	\subfloat[\label{fig:SWBLI_Mach}]{\includegraphics[width=0.55\textwidth]{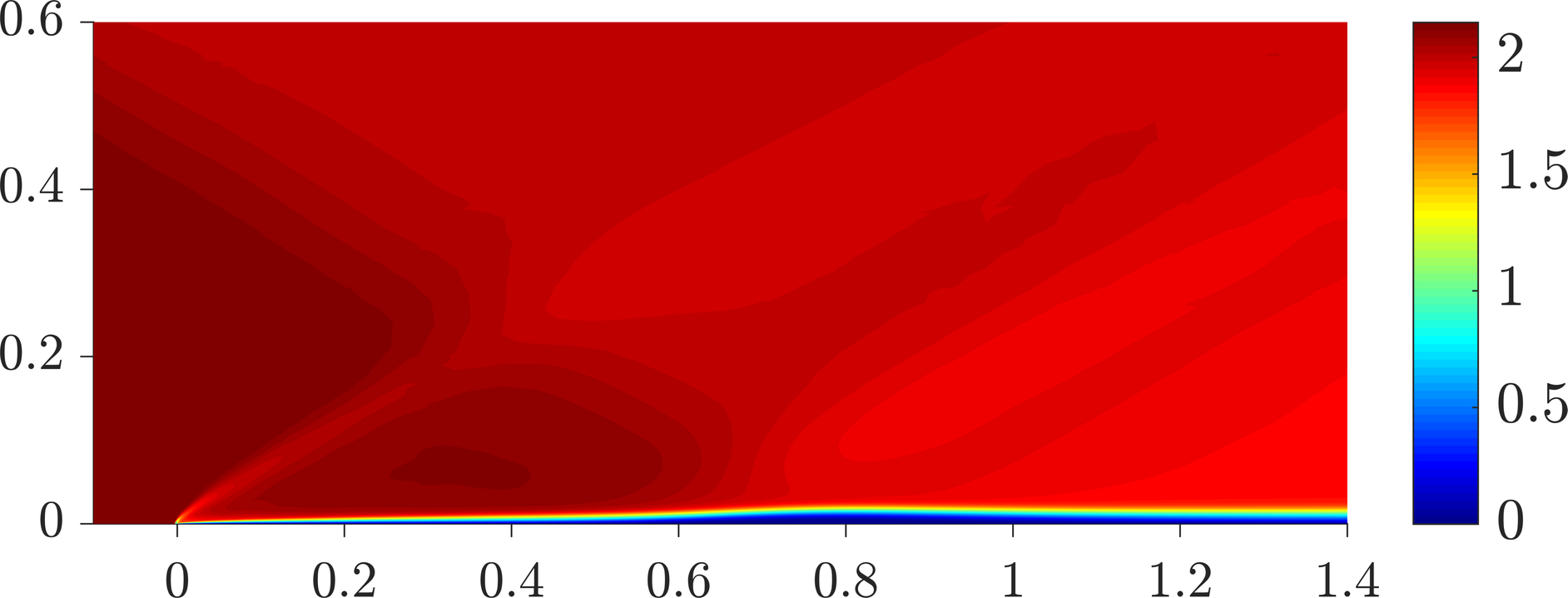}} \hfill
	\subfloat[\label{fig:SWBLI_MachZoom}]{\includegraphics[width=0.42\textwidth]{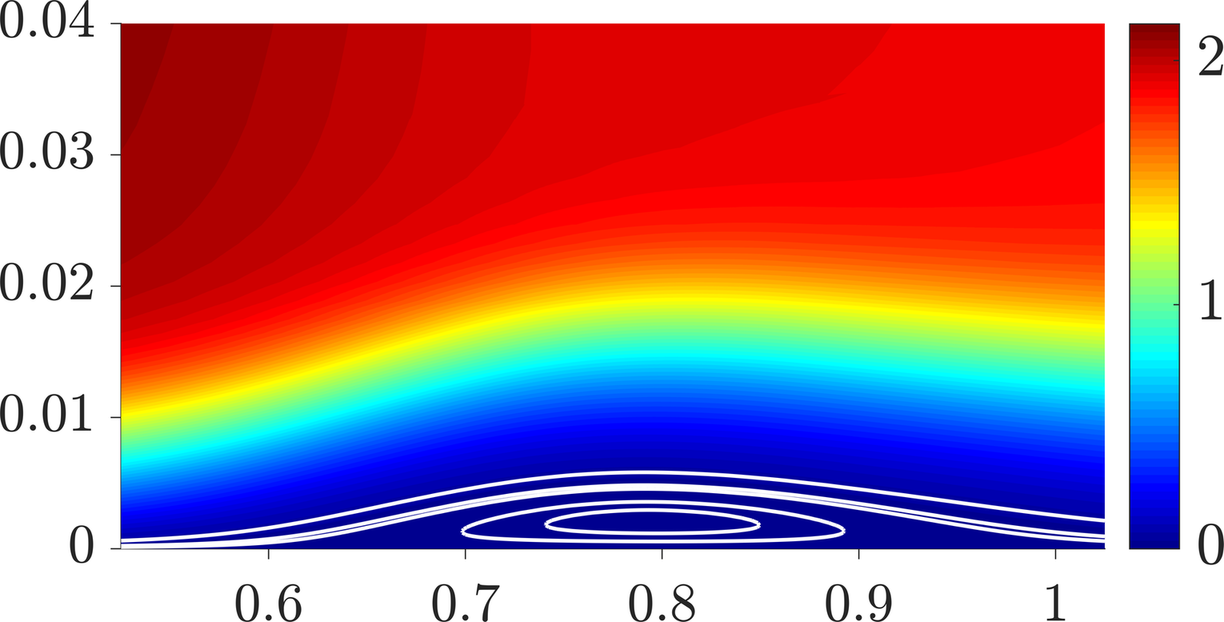}}
	\caption{Shock wave/boundary layer interaction - Mach number distribution obtained with the HLLEM Riemann solver and polynomial degree of approximation $k=3$. In (b), detail of the shock-induced separation bubble, where isolines of the Mach are drawn in white.}
	\label{fig:SWBLI_MachGlobal}
\end{figure}

The simulation is performed using the HLLEM Riemann solver due to its positivity-preserving properties in presence of shocks, contrary to Roe, and its superiority with respect to HLL or Lax-Friedrichs in resolving boundary layers.
The physics-based shock treatment involving an artificial bulk viscosity described in section~\ref{ssc:PhysicsBasedShockCapturing} is employed. The resulting flowfield is depicted in figure~\ref{fig:SWBLI_Mach}. The presence of shocks generated at different locations as well as the effect of the strong shock wave/boundary layer interaction can be observed.
Detail of the impingement region showing the separation bubble induced by the interaction between the reflecting shock wave and the boundary layer is illustrated in figure~\ref{fig:SWBLI_MachZoom}.

\begin{figure}[htbp]
	\subfloat[Pressure coefficient \label{fig:SWBLI_Cp}]{\includegraphics[width=0.49\textwidth]{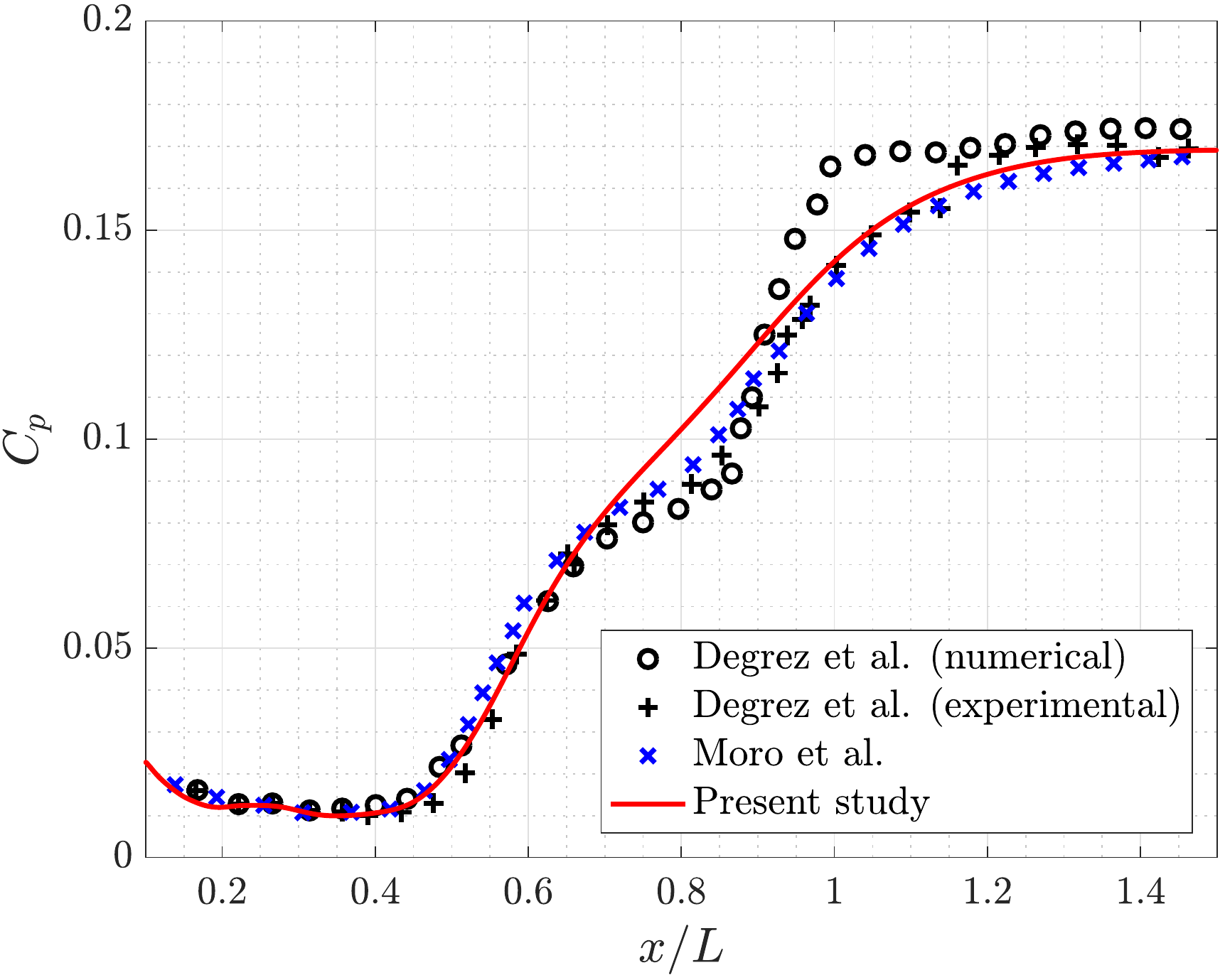}} \hfill
	\subfloat[Friction coefficient \label{fig:SWBLI_Cf}]{\includegraphics[width=0.474\textwidth]{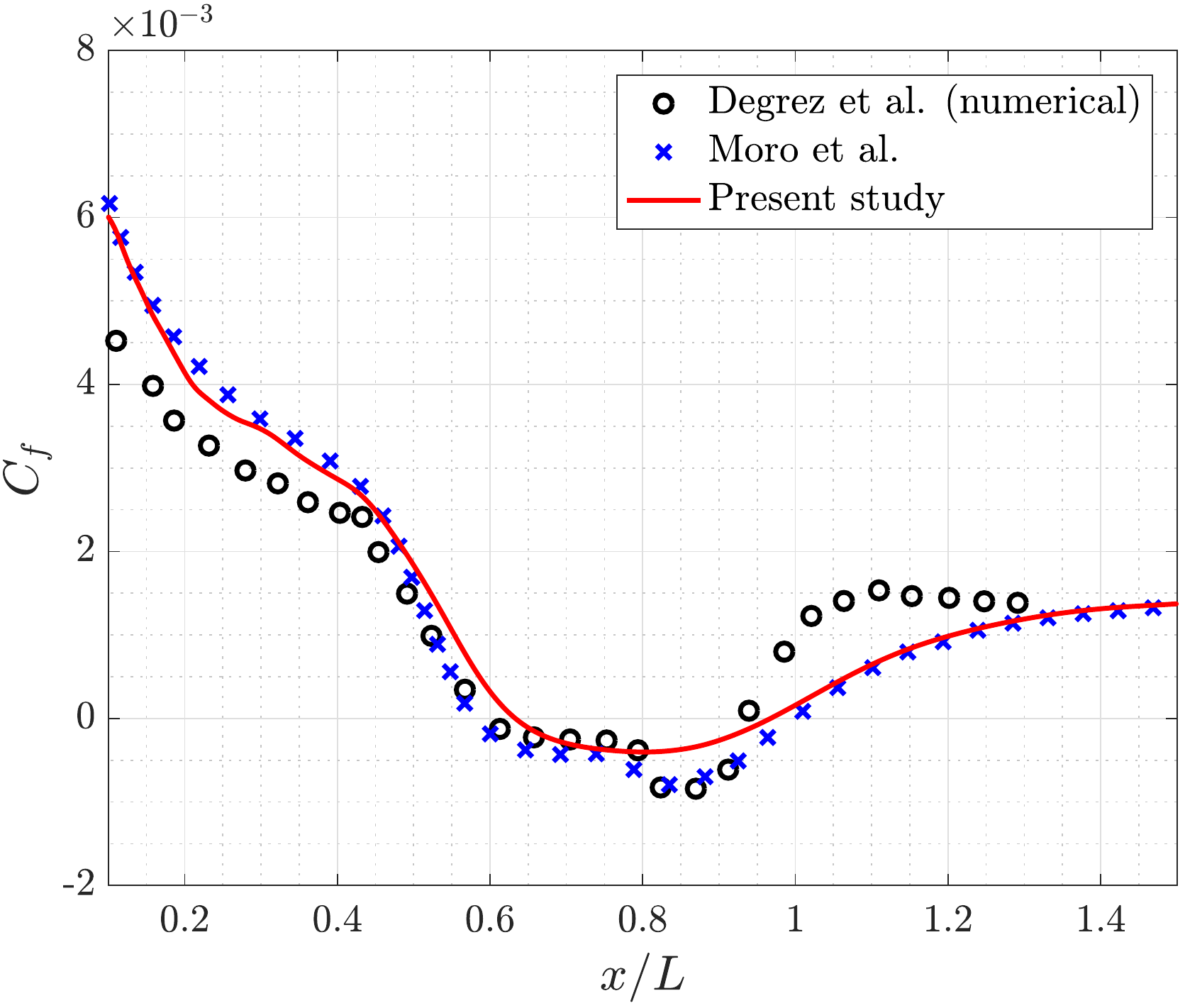}} 
	\caption{Shock wave/boundary layer interaction - Pressure (a) and friction (b) coefficients along the flat plate using the HLLEM Riemann solver and order of polynomial approximation $k=3$.}
	\label{fig:SWBLI_Coefficients}
\end{figure}

Finally, figure~\ref{fig:SWBLI_Coefficients} displays a comparison of the pressure coefficient and the skin friction coefficient computed using the HLLEM Riemann solver with benchmark results in \cite{Degrez1987,Moro2017}.

The computed pressure and skin friction coefficients show excellent agreement with both the experimental curve by Degrez et al. and the numerical solution by Moro et al., whereas the numerical curve by Degrez et al. deviates from the rest of results, especially downstream of the separation bubble.
The HLLEM computed solution lies on top of the reference results except for the region of shock impingement, where the highly anisotropic adapted meshes by Moro et al. outperform the presented results.
It is worth recalling that the HLLEM simulation in this study is performed on a mesh with no a priori refinement except for the boundary layer regions and the leading edge point.

This test case demonstrates a good behaviour of the HLLEM Riemann solver not only in the resolution of the boundary layer or in the treatment of shock waves in high-order but also in the strong interaction of these two flow features which challenges the performance of Riemann solvers.

\subsection{Supersonic flow over a compression corner}
\label{ssc:CompressionCorner}

The last case presented in this study considers the $\Ma_{\! \infty} = 3$ supersonic flow over a $10^\circ$ compression corner. This example represents a classical benchmark for viscous laminar compressible flow, first introduced by Carter~\cite{Carter1972} and later reproduced by several authors, see for example~\cite{Shakib1991,Aliabadi1993,Qamar2006,Mittal2001,Kotteda2014,Hung1976}.

\begin{figure}[htbp]
	\centering
	\includegraphics*[width=0.7\textwidth]{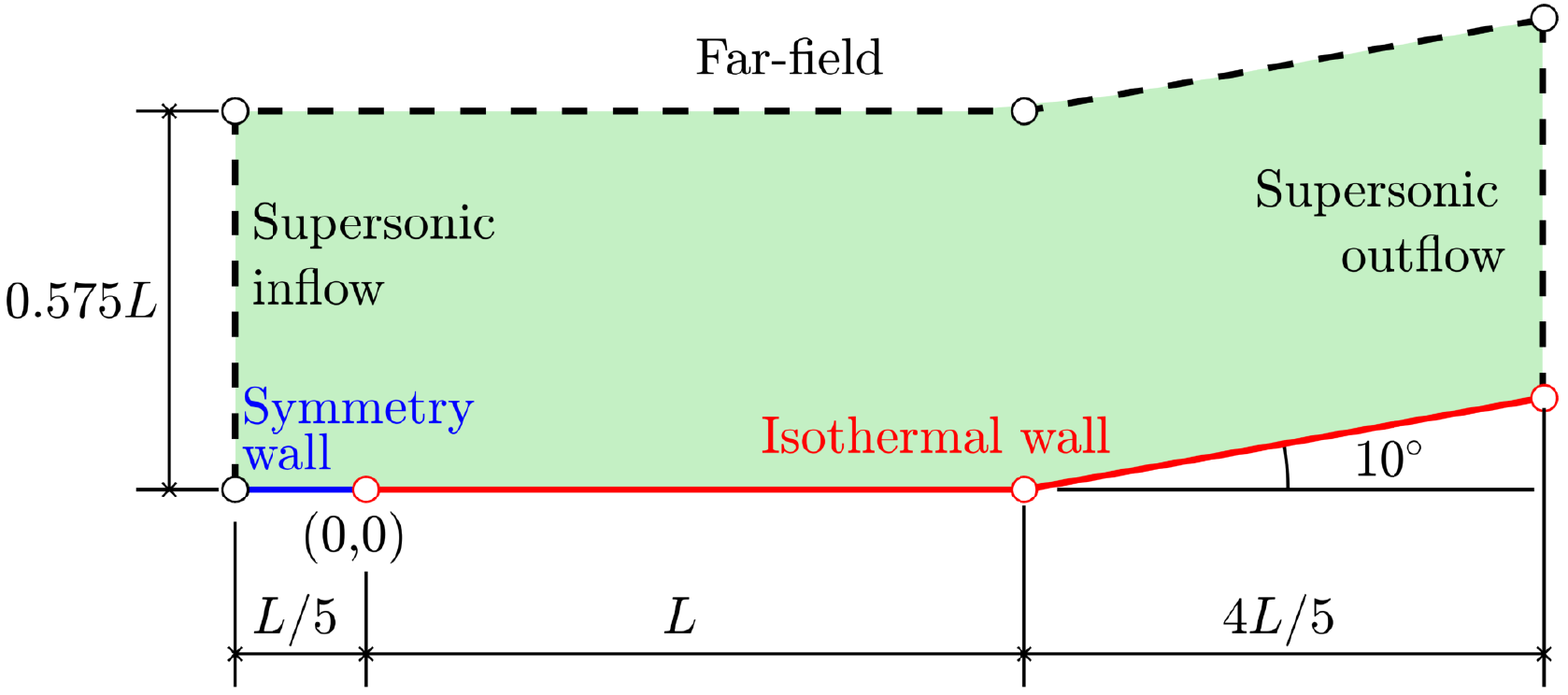}
	\caption{Supersonic flow over a compression corner - Sketch of the geometry and boundary conditions.}
	\label{fig:CompressionCorner_Domain}
\end{figure}

The setup of this problem consists of a laminar flow at $\Rey = 16,800$ over an isothermal flat plate of length $L$ (the characteristic length of the problem) ended with a $10^\circ$ wedge.
The isothermal surface is kept at the free-stream stagnation temperature, namely 
\begin{equation}
T_w = T_{\infty,0} = \frac{1}{(\gamma - 1) \Ma_{\! \infty}^2} \left(1 + \frac{\gamma - 1}{2} \Ma_{\! \infty}^2 \right).
\end{equation}
A sketch of the geometry and detail of the corresponding boundary conditions is depicted in figure~\ref{fig:CompressionCorner_Domain}.

\begin{figure}[htbp]
	\centering
	\subfloat[ \label{fig:compressioncorner_Mesh}]{\includegraphics[width=0.58\textwidth,valign=c]{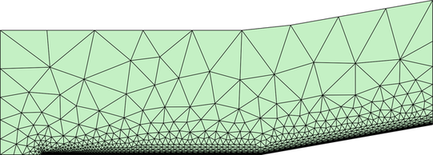}
	\vphantom{\includegraphics[width=0.38\textwidth,valign=c]{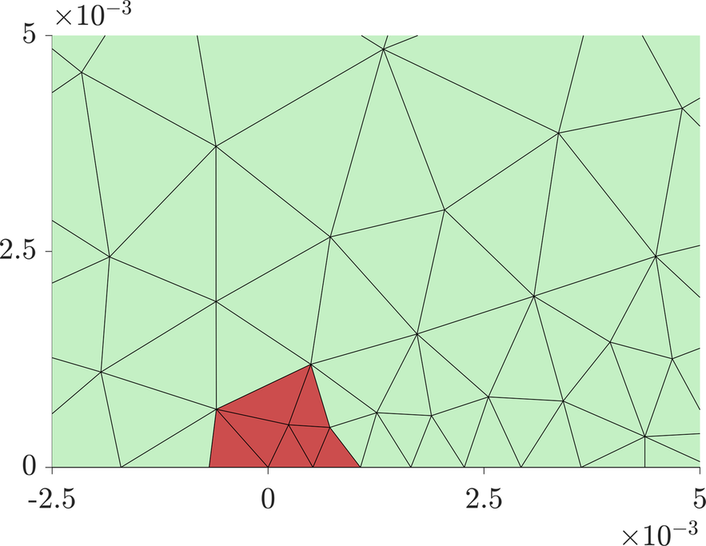}}} \hfill
	\subfloat[ \label{fig:compressioncorner_MeshPmap}]{\includegraphics[width=0.38\textwidth,valign=c]{compressioncorner_MeshPmap}}
	\caption{Supersonic flow over a compression corner - (a) Computational mesh and (b) detail of the leading edge, showing in red the elements employing a lower degree of approximation, $k=2$.}
	\label{fig:compressioncorner_MeshAll}
\end{figure}

The computational domain is discretised with 2,773 triangular elements of degree $k=3$, as illustrated in~\ref{fig:compressioncorner_MeshAll}.
In contrast to the shock wave/boundary layer interaction example, the leading edge of the flat plate is not rounded by means of a fillet, thus introducing a singularity. Such singular behaviour is alleviated by means of further refinement and by reducing the order of polynomial approximation to $k=2$ in the elements surrounding the singularity, as depicted in red in~\ref{fig:compressioncorner_MeshPmap}.

The boundary layer mesh consists of $\nlayers = 12$ layers of elements with the first layer located at a height of $h_0/L = 5\cdot 10^{-4}$ and a growing rate of $r=1.4$. The isothermal wall is divided into $\ndiv = 72$ elements.

\begin{figure}[htbp]
	\centering
	\subfloat[Density \label{fig:compressioncorner_Density}]{\includegraphics[width=0.49\textwidth]{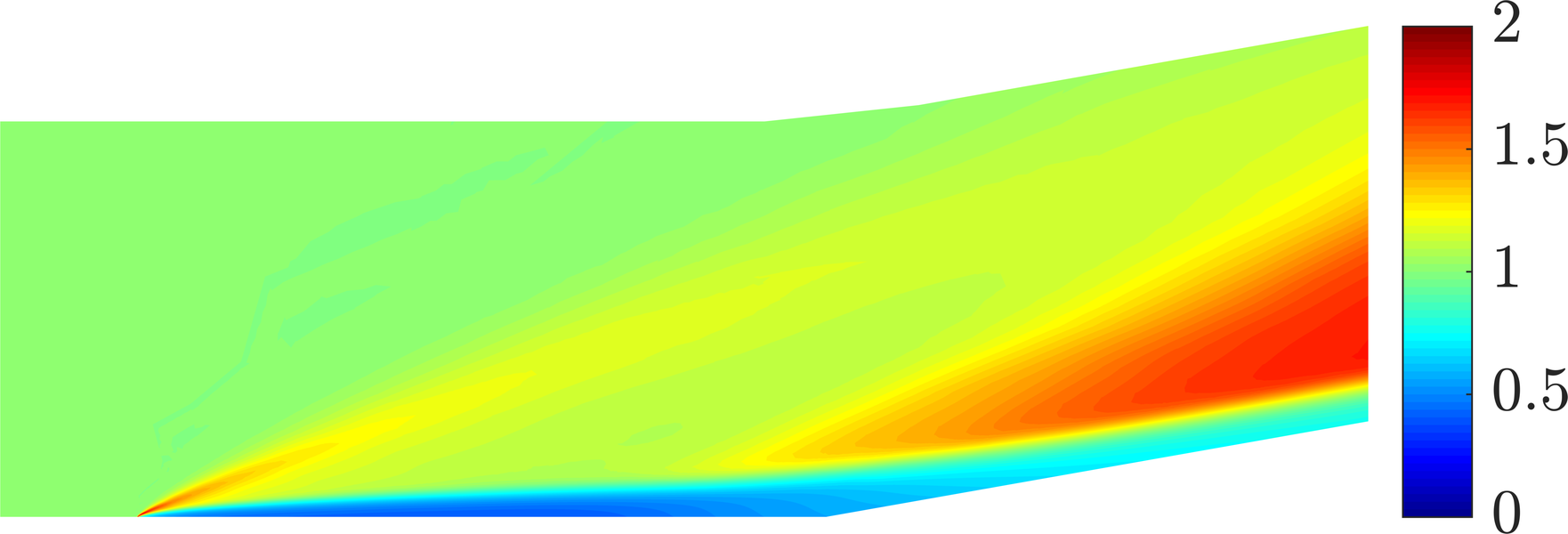}} \hfill
	\subfloat[Mach \label{fig:compressioncorner_Mach}]{\includegraphics[width=0.48\textwidth]{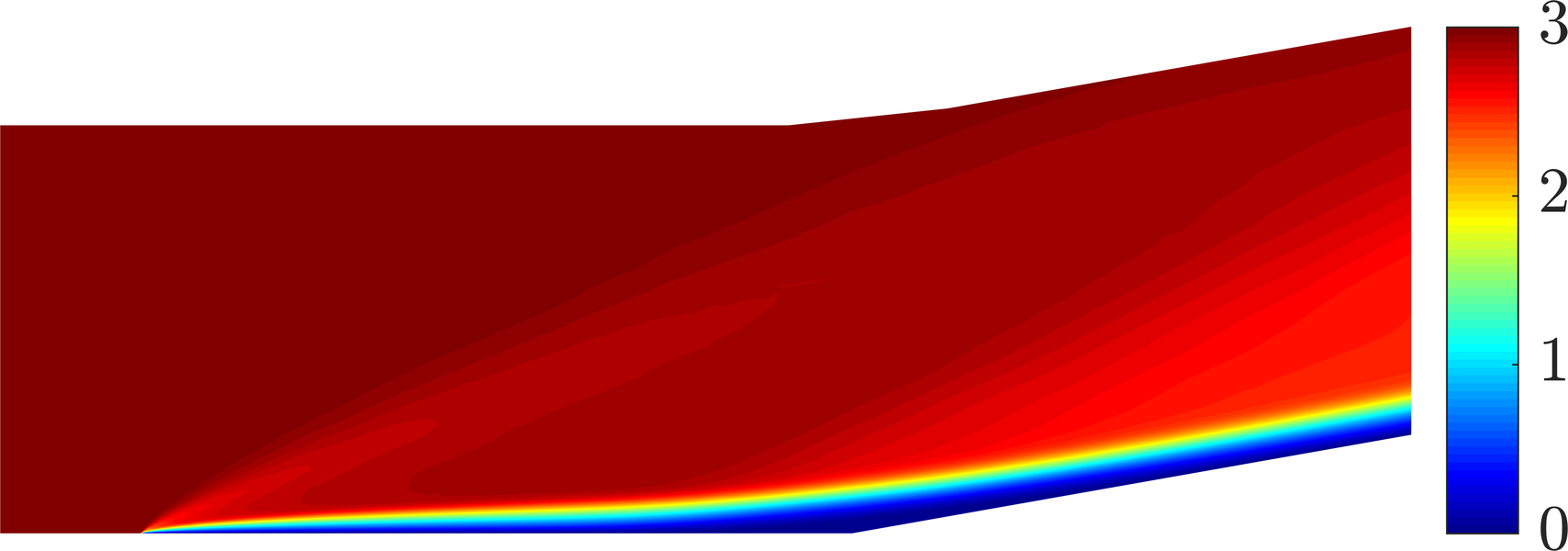}}
	\caption{Supersonic flow over a compression corner - Density (a) and Mach number (b) distributions using the HLLEM Riemann solver with a combined polynomial degree of approximation $k=2$ and $k=3$.}
	\label{fig:compressionCorner_Solution}
\end{figure}

The physics-based shock capturing based on artificial bulk viscosity described in~\ref{ssc:PhysicsBasedShockCapturing} is employed for the simulation.
The resulting flowfield obtained with the HLLEM Riemann solver is presented in figure~\ref{fig:compressionCorner_Solution}. The density field in figure~\ref{fig:compressioncorner_Density} illustrates the regions of high compression, namely the shock wave generated at the leading edge and the compression fan induced by the wedge.

\begin{SCfigure}[][htbp] 
	\centering
	\includegraphics[width=0.6\textwidth]{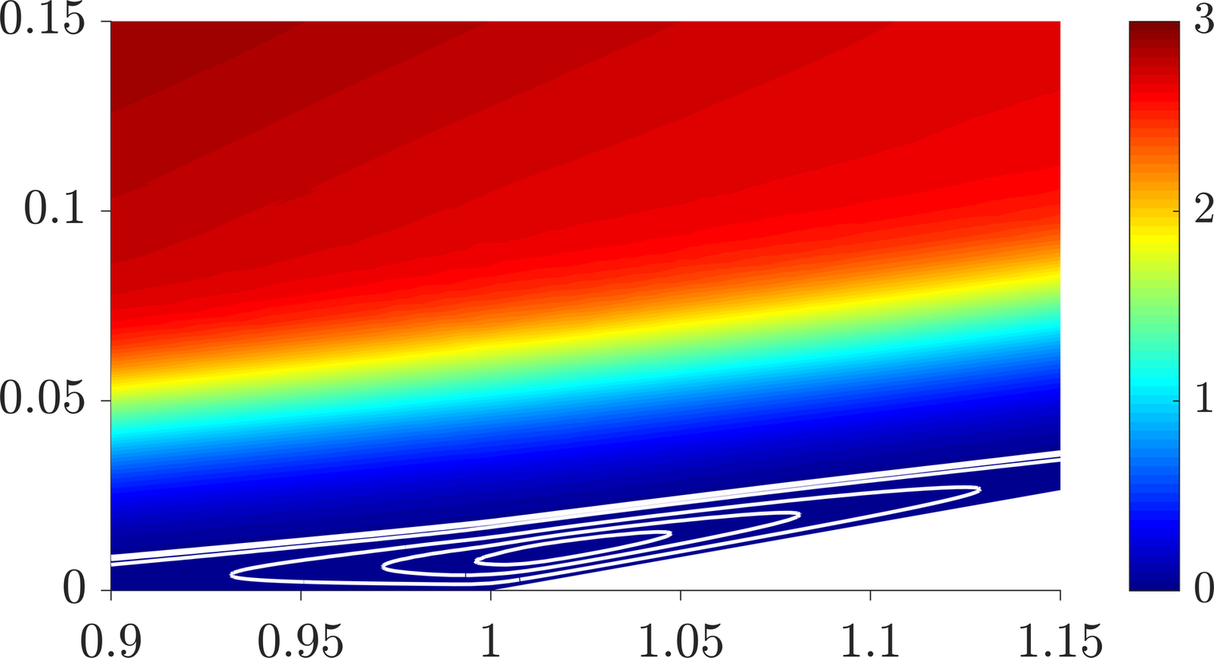}
	\caption{Supersonic flow over a compression corner - Detail of the Mach number distribution around around the corner, using the HLLEM Riemann solver with a combined polynomial degree of approximation $k=2$ and $k=3$. Isolines of the Mach are drawn in white.}
	\label{fig:compressionCorner_MachDetail}
\end{SCfigure}

Good resolution of the flow solution can be observed in figure~\ref{fig:compressionCorner_MachDetail}, where the separation bubble induced by the corner is depicted.

\begin{figure}[htbp]
	\subfloat[Pressure \label{fig:compressioncorner_Cp}]{\includegraphics[width=0.49\textwidth]{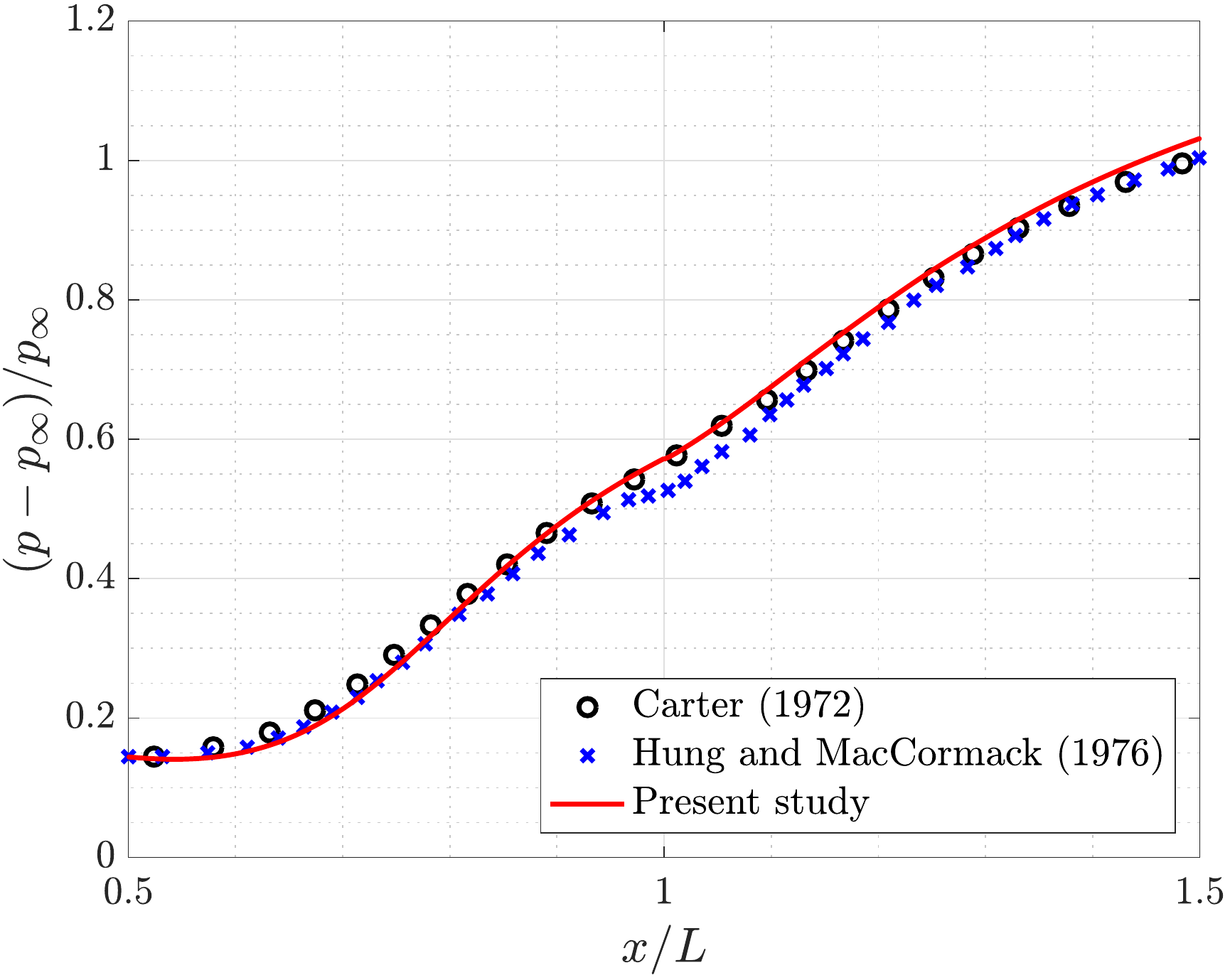}} \hfill
	\subfloat[Friction coefficient \label{fig:compressioncorner_Cf}]{\includegraphics[width=0.49\textwidth]{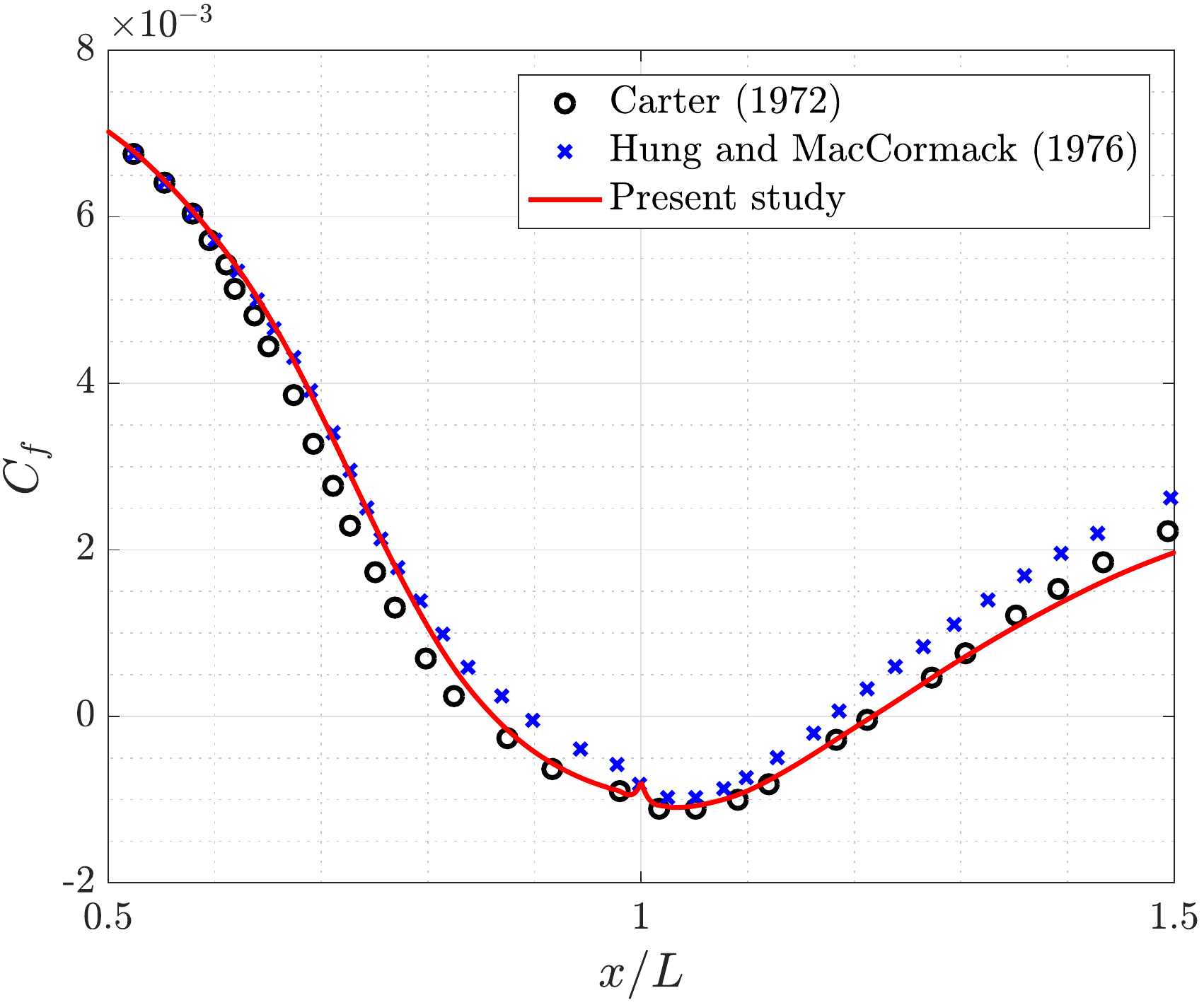}} 
	\caption{Supersonic flow over a compression corner - Pressure (a) and friction coefficient (b) along the flat plate using the HLLEM Riemann solver with a combined polynomial degree of approximation $k=2$ and $k=3$.}
	\label{fig:compressionCorner_coefficients}
\end{figure}

A qualitative comparison of the obtained results is carried out through the wall pressure and the skin friction coefficient. Figure~\ref{fig:compressionCorner_coefficients} compares such quantities with respect to the reference results by Carter~\cite{Carter1972} and Hung and MacCormack~\cite{Hung1976}, showing an excellent agreement.
Additional numerical results available in the literature such as~\cite{Shakib1991,Aliabadi1993,Qamar2006,Mittal2001,Kotteda2014} are not included in the comparison for the sake of readability because of the similarity among them.

Finally, the position of the separation, $x_s$, and reattachment, $x_r$, points, gathered in table~\ref{tb:compressioncorner_reattachment}, allows a quantitative assessment of the computed solution.

\begin{table} [htbp]
	\centering
	\caption{Supersonic flow over a compression corner - Position of the separation, $x_s$, and reattachment, $x_r$, points around the wall.}
	\begin{tabular}{ L{6.5cm} C{2cm} C{2cm}}
		\toprule
		Reference & $x_s/x_c$ & $x_r/x_c$ \\
		\midrule
		Present study & 0.86 & 1.20 \\				
		\midrule
		Carter~\cite{Carter1972} & 0.84 & 1.22 \\
		Hung and MacCormack~\cite{Hung1976} & 0.89 & 1.18 \\
		Shakib et al.~\cite{Shakib1991} & 0.88 & 1.17 \\
		Mittal and Yadav~\cite{Mittal2001} & 0.89 & 1.13 \\
		Kotteda and Mittal~\cite{Mittal2001} & 0.88 & 1.17 \\
		\bottomrule
	\end{tabular}
	\label{tb:compressioncorner_reattachment}
\end{table}

The obtained results show a strong consistency with respect to those available in the literature, proving the good performance of the high-order HDG solver.

\section{Concluding remarks}
\label{sc:Conclusions}

This paper presents a review of the formulation of inviscid and viscous compressible flows, i.e. the Euler and the compressible Navier-Stokes equations, in the context of high-order hybridised discontinuous Galerkin methods.
Moreover, it introduces a unified framework for the derivation of traditional Riemann solvers, namely Lax-Friedrichs and Roe, already formulated in HDG, and HLL and HLLEM Riemann solvers, which are devised for the first time for hybridised discretisations.
According to the HDG rationale, the intermediate state utilised to evaluate the numerical fluxes is constructed by means of the HDG hybrid variable and the information of the Riemann solver itself is encapsulated in the HDG stabilisation matrix.
In addition, the present formulation of the compressible Navier-Stokes equations introduces a new choice for the mixed variables employed to describe the viscous flux tensor, namely the deviatoric strain rate tensor and the temperature gradient. Such election for the mixed variables allows to impose pointwise the symmetry of the stress tensor with reduced computational cost, while retrieving optimal accuracy.

Optimal convergence properties of the HDG discretisation have been verified using Lax-Friedrichs, Roe, HLL and HLLEM Riemann solvers both for inviscid and viscous cases and for a wide range of the Reynolds number.
HDG demonstrates its ability to approximate the conserved quantities as well as the viscous stress and the heat flux with optimal order of convergence, $k+1$. Whereas the role of the Riemann solver shows little effect in the precision of the approximate primal variables, significant differences are noticed in the precision of the approximated mixed variables. In particular, HLLEM and Roe Riemann solvers yield a gain in accuracy in the approximation of the heat flux and viscous stress, specially as the Reynolds number increases.

Then, a set of 2D numerical benchmarks has been presented to show the advantages of high-order approximations for compressible flow problems and the capabilities of the novel HLL and HLLEM Riemann solvers in different flow regimes, from subsonic to supersonic, with special attention to its comparison with well-established Lax-Friedrichs and Roe Riemann solvers in the context of HDG.

In particular, HLL-type Riemann solvers exhibit a superior performance in supersonic cases, illustrating their positivity preserving properties. This allows a robust and parameter-free strategy in the solution of supersonic flows involving shock waves, contrary to Roe Riemann solver, which may fail to produce physically admissible solutions because of a lack of dissipation.
Furthermore, HLLEM Riemann solver demonstrates its ability to preserve contact or shear layers, likewise Roe, producing results that introduce less numerical dissipation than HLL and Lax-Friedrichs, and displaying a major advantage in the approximation of boundary layers.

Finally, a couple of benchmarks involving the interaction of boundary layers and shock waves demonstrate the overall good performance of the high-order HDG method, equipped with a shock-capturing technique based on artificial viscosity and an HLLEM Riemann solver, in the resolution of problems with such combination of viscous and inviscid-type phenomena.

\section*{Acknowledgements}

This work was partially supported by the Spanish Ministry of Economy and Competitiveness (Grant number: DPI2017-85139-C2-2-R). J.V.P. was supported by the Spanish Ministry of Economy and Competitiveness, through the Mar\'ia de Maeztu programme for units of excellence in R\&D (Grant number: MDM-2014-0445). M.G. and A.H. are also grateful for the support provided by the Spanish Ministry of Economy and Competitiveness through the Severo Ochoa programme for centres of excellence in RTD (Grant number: CEX2018-000797-S) and the Generalitat de Catalunya (Grant number: 2017-SGR-1278). R.S. also acknowledges the support of the Engineering and Physical Sciences Research Council (Grant number: EP/P033997/1).

\bibliographystyle{elsarticle-num}
\bibliography{CompressibleHDG}   

\end{document}